%% file: breakdown.tex
\documentclass[12pt]{amsart}
\usepackage{amssymb,latexsym,amsthm, amsmath, comment, bbm, tikz-cd, fullpage, multirow, lscape, subcaption}
\usepackage[all]{xy}
\usepackage{graphicx, comment, libertine, wrapfig}
\usepackage[utf8]{inputenc}
\usepackage[T2A]{fontenc}

\ifx\pdfoutput\undefined
\usepackage{hyperref}
\else
\usepackage[pdftex,colorlinks=true,linkcolor=blue,urlcolor=blue]{hyperref}
\fi

\usepackage[frozencache]{minted}
\usepackage{anyfontsize}
\usetikzlibrary{cd,patterns,decorations.pathmorphing,shapes}
\DeclareMathOperator{\coker}{coker}

\theoremstyle{plain}

\theoremstyle{definition}
\newtheorem{definition}{Definition}[section]
\newtheorem{theorem}{Theorem}[section]

\newtheorem{corollary}{Corollary}[section]

\newtheorem{conjecture}{Conjecture}[section]
\newtheorem{remark}{Remark}[section]
\newtheorem*{prop*}{Proposition}

\newtheorem*{ack}{Acknowledgements}
\newtheorem{example}{Example}[section]

\newcommand{\bC}{\mathbb{C}}

\newcommand{\bP}{\mathbb{P}}
\newcommand{\bQ}{\mathbb{Q}}
\newcommand{\bR}{\mathbb{R}}
\newcommand{\bT}{\mathbb{T}}
\newcommand{\bZ}{\mathbb{Z}}
\newcommand{\cA}{\mathcal{A}}
\newcommand{\cT}{\mathcal{T}}
\newcommand{\fp}{\mathfrak{p}}
\newcommand{\co}{\colon\,}
\newcommand{\ch}{\mathrm{ch}}
\renewcommand{\Col}{\mathrm{Col}}

\begin{document}

\title{New results on tilings via cup products and Chern characters on tiling spaces}
\author{Jianlong Liu}
\address{University of Texas at Austin}
\email{jlliu@utexas.edu}
\author{Jonathan Rosenberg}
\address{University of Maryland, College Park}
\email{jmr@umd.edu}
\author{Rodrigo Trevi\~no}
\address{University of Maryland, College Park}
\email{rodrigo@trevino.cat}
\date{\today}
\begin{abstract}
  We study the cohomology rings of tiling spaces $\Omega$ given by cubical
  substitutions.  While there have been many calculations before of
  cohomology \emph{groups} of such tiling spaces, the innovation
  here is that we use computer-assisted methods to compute the
  cup-product structure.  This leads to examples of substitution
  tilings with isomorphic cohomology groups but different
  cohomology rings.  Part of the interest in studying the cup
  product comes from Bellissard's \emph{gap-labeling conjecture},
  which is known to hold in dimensions $\le 3$, but where a
  proof is known in dimensions $\ge 4$ only when the Chern character
  from $K^0(\Omega)$ to $H^*(\Omega,\bQ)$ lands in $H^*(\Omega,\bZ)$.
  Computation of the cup product on cohomology often makes it
  possible to compute the Chern character.
  We introduce a natural generalization of the gap-labeling conjecture,
  called the \emph{equivariant gap-labeling conjecture}, which
  applies to tilings with a finite symmetry group.  Again this
  holds in dimensions $\le 3$, but we are able to show that it
  \emph{fails} in general in dimensions $\ge 4$.  This, plus
  some of our cup product calculations, makes it plausible that
  the gap-labeling conjecture might fail in high dimensions.
\end{abstract}
\maketitle

\section{Introduction}
\label{sec:intro}
This paper is about tiling spaces $\Omega$ coming from primitive
substitution tilings of $\bR^d$.  All the tilings in this paper will
be self-similar and have finite local complexity.  Under those
hypotheses, the action of $\bR^d$ on $\Omega$ by translations
is minimal and uniquely ergodic \cite{solomyak:SS}.  We begin with a
review of the basic definitions and literature on tilings in section
\ref{sec:tilingspaces}.

The problems that we will study here were motivated by the
gap-labeling conjecture of Bellissard \cite{bellissard:K}.
His idea was that in a quasi-crystal (an almost periodic physical
system modeled mathematically by a tiling), physical properties of
the material are determined by the spectra of almost periodic
Schr{\"o}dinger operators which live in a $C^*$-algebra
$\cA_p(\Omega)$ determined by the tiling, so that gaps in the spectra
can be determined by computing $K_0(\cA_p(\Omega))$ and the
map from this $K$-group to the reals induced by the trace.  The algebra
$\cA_p(\Omega)$  is highly noncommutative and thus complicated to deal
with, but it contains a natural copy of the commutative algebra
$C(\Omega)$ of functions on the tiling space.  Bellissard noticed that
at least in many cases, the $K$-theory of $\cA_p(\Omega)$ all lies in
the image of the $K$-theory of this subalgebra, which can be computed
purely topologically.  The gap-labeling conjecture asserts that this
should always be the case.  Several papers
\cite{BOO:GLT,BBG:GLT,KP:GLT} claimed to prove the conjecture, but
they all implicitly assumed that the top-dimensional piece of the
Chern character $\text{ch}_d\co K^{-d}(\Omega)\to\check{H}^d(\Omega;\bQ)$
gives an isomorphism
$K^{-d}(\Omega)\xrightarrow{\cong}\check{H}^d(\Omega;\bZ)$.
(See \cite[\S9]{ADRS:bloch} for a detailed discussion.)

Our goal in this paper is to compute the structure of the integral
cohomology ring $H^*(\Omega;\bZ)$ for several examples of cubical
substitution tilings.  This is computation-intensive and requires
computer assistance for the calculation of the cup-product.
Calculations are done in sections \ref{sec:cubical} and
\ref{sec:breakdown}.  In several cases, we are also able to determine
at least an important part of the Chern character map
$\text{ch}_d\co K^{-d}(\Omega)\to\check{H}^d(\Omega;\bQ)$.
Among the major results are:
\begin{enumerate}
\item  There are cases of tiling spaces with isomorphic cohomology as
  groups but not as \emph{rings} (section \ref{sec:distinguishing}).
\item There are cases where the Chern character does \emph{not} give an
  isomorphism $K^{-d}(\Omega)\to\check{H}^d(\Omega;\bZ)$ (section
  \ref{subsec:ex}). 
\end{enumerate}

Unfortunately we are not able with the present techniques to produce a
counterexample to the gap-labeling conjecture in high dimensions,
though we would not be surprised if such a counterexample exists.

However, we show that there is a natural generalization of the
gap-labeling conjecture, which we call the \emph{equivariant
  gap-labeling conjecture}, that can be formulated for tiling spaces
with a finite symmetry group.  Roughly speaking, this involves
replacing $K$-theory by equivariant $K$-theory.  This conjecture is
formulated and discussed in section \ref{subsec:equivgap} and
\ref{subsec:equivChern}.  We show that this equivariant conjecture
fails in general when $d\ge 4$ (see Theorem \ref{thm:equivlabeling}
and section \ref{subsec:equivCounter}), which provides additional evidence
that the original gap-labeling conjecture probably fails in high
dimensions. 

\begin{ack}
J.L. was supported in part by NSF grant DMS-1937215. R.T. was
supported by the NSF award DMS-2143133 Career. We are grateful to
Michael Baake, Franz G\"ahler, and Lorenzo Sadun, for insightful
conversations during the preparation of this manuscript, and for
suggesting the quotient of the squiral substitution. Thanks to
Lorenzo Sadun for pointing out a mistake in the first pair of
examples, and to the anonymous referee for useful comments.
\end{ack}

\section{Tiling spaces and their cohomology}
\label{sec:tilingspaces}
This section covers the necessary background on tiling spaces and
their topological invariants.

\subsection{Basics of tilings}
To define a tiling, we start by assuming we are given a finite
set of polytopes in $\bR^d$ called \textbf{prototiles}.
A \textbf{tile} is a translate $t+\tau$, $\tau\in \bR^d$, of
a prototile $t$, along with a choice of color or label
$c\in \Col$ chosen from a set $\Col$ of \textbf{colors}.
We will always assume that $\Col$ is finite.
If $\Col$ has only one element,
a common scenario, then one can ignore the color.
A \textbf{tiling} $\cT$ of $\bR^d$ is a decomposition of $\bR^d$
as a union of tiles which only intersect along their faces. We write
$\bR^d = \bigcup_j t_j$, where each $t_j$ is the
\textbf{support of the tile $(t_j,c_j)$}, and $c_j\in \Col$ 
is the color of the tile. In a slight abuse of notation, we will often
fail to distinguish between tiles and their supports,
but this should not cause confusion from
context. A \textbf{patch} of $\mathcal{T}$ is the union of finitely
many tiles. Two tiles $t_1,t_2$ of a
tiling $\mathcal{T}$ are \textbf{translation-equivalent} if there exists a
$\tau\in\mathbb{R}^d$ such that $t_1 = t_2-\tau$, i.e., they are
translations of the same prototile, and their labels are
equal. Two patches $P_1,P_2$ of $\cT$ are translation
equivalent if there is a $\tau\in\bR^d$ such that $P_1 = P_2-\tau$.
The translation $\varphi_\tau$ of tiles and patches by vectors
$\tau\in \bR^d$ extends to all of a tiling $\cT$, and we denote by
$\varphi_\tau(\cT)$ the tiling consisting of the union of
tiles of the form $\varphi_\tau(t)$, where $t$ is a tile of
$\mathcal{T}$. 

Let $R_*>0$ be the smallest number so that any ball of radius $R_*$ in
$\bR^d$ contains a tile of $\mathcal{T}$. For any $R\geq R_*$, an $R$-patch
of $\mathcal{T}$ is the largest patch of $\mathcal{T}$ completely
contained in a fixed ball of radius $R$. The tiling $\mathcal{T}$ has
\textbf{finite local complexity} or
\textbf{FLC} if, given any $R>0$, the set of equivalence
classes (under translation) of $R$-patches of the family of tilings
$\{\varphi_t(\mathcal{T})\}_{t\in\mathbb{R}^d}$ is finite. All tilings
considered in this paper will have FLC. 

Given a tiling $\mathcal{T}$, define the metric $d$ on the set of all
translates $\{\varphi_t(\mathcal{T})\}_{t\in\mathbb{R}^d}$ as
follows. Set
\[
d'(\varphi_t(\cT),\varphi_s(\cT)) = \inf\left\{
\varepsilon>0: \substack{\mbox{there is a $\tau\in B_\varepsilon(0)$
    such that the $\varepsilon^{-1}$-patch of }
  \\ \mbox{$\varphi_{t+\tau}(\cT)$ is equal to the
    $\varepsilon^{-1}$-patch of $\varphi_{s}(\cT)$
}}\right\}
\]
and define
\[
d(\varphi_t(\cT), \varphi_s(\cT)) =
\min\{1,d'(\varphi_t(\cT), \varphi_s(\cT))\}.
\]
This is a metric \cite{solomyak:SS}. The \textbf{(translational) tiling space} $\Omega_\cT$ of
$\mathcal{T}$ is the metric completion of the set of translates of
$\mathcal{T}$ with respect to the metric above, i.e.
\[\Omega_\mathcal{T} :=
\overline{\{\varphi_t(\mathcal{T}):t\in\mathbb{R}^d\}}.
\]
If $\cT$ is FLC, then $\Omega_\cT$ is a compact metric
space, with local product structure $B\times C$, where $B$ is a
Euclidean ball and $C$ is a zero-dimensional compact space
(usually, but not always, a Cantor set). The translation of
$\cT$ by vectors of $\bR^d$ extends to an action of
$\bR^d$ by translations. \\
\indent Let us mention two more topological spaces associated to \(\mathcal{T}\), where one incorporates rotations in addition to translations, that is, \(SE(d)=\mathbb{R}^d\rtimes SO(d)\). For details, see \cite[\S 4.1--4.5]{sadun:book}. The metric is first relaxed to allow \(\varepsilon\)-rotations, then one can define
\[
\Omega_\mathcal{T}^\textnormal{rot}=\overline{\{\phi(\mathcal{T}):\phi\in SE(d)\}}.
\]
It has a canonical transversal
\[
\Omega_\mathcal{T}^0=\Omega_\mathcal{T}^\textnormal{rot}/SO(d)\cong\Omega_\mathcal{T}/\{\textnormal{rotational symmetries}\}.
\]
which we call the \textbf{rotational tiling space} of \(\mathcal{T}\).

\subsection{Substitution tilings}
Let $(t_1,c_1),\dots, (t_m,c_m)\subset \mathbb{R}^d\times \Col$ be
compact, connected sets with non-empty interior, where $\Col$ is a finite
set. If $A\subset GL(d,\bR)$ is an expanding matrix, then a
\textbf{substitution rule on $t_1,\dots, t_m$ with expansion $A$} is a
decomposition of $At_i$ as
\[
At_i = \bigcup_{j=1}^m\bigcup_{\tau\in \Lambda_{ij}} t_j-\tau
\]
for any $i$, where each $\Lambda_{ij}$ is a finite set. The
\textbf{substitution matrix} corresponding to this substitution rule
is the doubly-indexed collection of numbers $S$ defined by $S_{ij} =
|\Lambda_{ij}|$. The substitution is \textbf{primitive} if there is a
$k\in\mathbb{N}$ such that the substitution matrix $S^k$ has all
positive entries. Any substitution rule can be iterated: the set
$A^2t_i$ can be tiled as
\[
A^2t_i = \bigcup_{j=1}^m\bigcup_{\tau\in \Lambda_{ij}} At_j-A\tau,
\]
where $At_j$ can be decomposed as above. A tiling $\mathcal{T}$ is
given by a substitution rule with expansion $A$ if any patch of
$\mathcal{T}$ is translation equivalent to a subpatch of $A^kt_i$ for
some $k\in\mathbb{N}$, $1\leq i\leq m$, where $A^kt_i$ is a tiled set
as above. In that case, a patch which is translation-equivalent to one of the form $A^kt_i$ as above is a \textbf{level-$k$ supertile}.

Let $\mathcal{T}$ be a FLC tiling with prototiles $t_1,\dots,
t_m$. For a tile $t\in\mathcal{T}$, let $\mathcal{T}(t)$ be the patch
consisting of all tiles of $\mathcal{T}$ which intersect $t$. This
type of patch is called a \textbf{collared tile} of $\mathcal{T}$. Now
consider the product $\Omega_\cT\times \bR^d$ with the
product topology, where $\Omega_\cT$ has the discrete topology
and $\bR^d$ the traditional topology. Let $\sim_1$ be the
equivalence relation on this product which declares
$(\mathcal{T}_1,u_1)\sim_1(\mathcal{T}_2,u_2)$ if
$\mathcal{T}_1(t_1)-u_1 = \mathcal{T}_2(t_2)-u_2$ for some tiles
$t_1,t_2$ with $u_i\in t_i\in\mathcal{T}_i$. The
\textbf{Anderson-Putnam} complex of $\Omega_\mathcal{T}$ is the set
$\Gamma_\mathcal{T}:=(\Omega_\mathcal{T}\times\mathbb{R}^d)/\sim_1$. Let
$\pi_1\co\Omega_\cT\times\bR^d\rightarrow \Gamma_\cT$ be the quotient map.
$\Gamma_{\cT}$ is a flat branched manifold. 

If $\mathcal{T}$ is a primitive substitution tiling then there exists a
locally expanding affine map
$\gamma\co\Gamma_\cT\rightarrow  \Gamma_\cT$ such that there is a
homeomorphism of tiling spaces
      \begin{equation}
        \label{eqn:AP}
        \Omega_\mathcal{T} \cong \varprojlim (\Gamma_\mathcal{T},\gamma),
      \end{equation}
which is the seminal result of Anderson and Putnam \cite{AP}.
In this case \cite{solomyak:SS}, the action of $\bR^d$ by translations is
\textbf{uniquely ergodic}: there is a unique $\bR^d$-invariant ergodic
probability measure on $\Omega$, which will be denoted by $\mu$.
By the local product structure of $\Omega$, this measure is locally
of the form $\mathrm{Leb}\times \nu$, where $\mathrm{Leb}$ is the Lebesgue
measure on the local Euclidean component and $\nu$ is a measure on the
local Cantor component and will be discussed further in \S \ref{subsec:RS}.

\begin{remark}
  \label{rem:force}
  The reader may wonder why $\Omega_\mathcal{T}$ may not be recovered
  as the inverse limit obtained by the substitution rule directly on
  tiles and not necessarily on collared tiles, as required for the
  construction of the Anderson-Putnam complex $\Gamma_\mathcal{T}$
  above. Although an inverse limit would be well-defined in such case,
  it may not necessarily be homeomorphic to the tiling space
  $\Omega_\mathcal{T}$. This happens if and only if the substitution
  rule \textbf{forces the border} (see \cite{AP} for a proof of this
  result). To quote \cite{sadun:book}, ``A substitution rule
  \textbf{forces the border} if there exists a positive integer $n$
  such that any two level-$n$ supertiles of the same type have the
  same pattern of neighboring tiles''. Although the concept of forcing
  the border will not be central in this paper, it will be mentioned
  in a couple of the examples below, and we direct the reader to
  \cite{AP} or \cite{sadun:book} to learn why in that case one can use
  uncollared tiles to form the desired inverse limit. 
\end{remark}

\begin{definition}
  A \textbf{cubical substitution} in dimension $d$ with $m$ colors with
  expansion $\lambda\in\mathbb{N}$ is a substitution rule on
  $[0,1]^d \times \{1,\dots, m\}$ with expansion matrix
  $A = \lambda\cdot \mathrm{Id}$.
\end{definition}
Note that there are $m^{\lambda^d+1}$ different cubical substitution rules
in dimension $d$ with $m$ colors and expansion $\lambda$.

\subsection{Topological invariants}
Suppose that $\cT$ is a primitive substitution tiling.
The Anderson-Putnam homeomorphism (\ref{eqn:AP}) allows us to express
the cohomology and $K$-theory of $\Omega = \Omega_\cT$ in digestible terms:
\begin{equation}
\label{eqn:dirLimits}
\check{H}^*(\Omega;\bZ) = \varinjlim
\left(\check{H}^*(\Gamma;\bZ),\gamma^* \right)\hspace{.6in}\mbox{ and}
\hspace{.6in}K^*(\Omega) = \varinjlim \left(K^*(\Gamma),\gamma^* \right).
\end{equation}
Recall that the Chern character
$\ch\co K^*(\Omega)\rightarrow \check{H}^*(\Omega;\bQ)$
is defined by
\begin{equation}
\label{eqn:chern}
\ch([L]) = \sum_{k\geq 0}\frac{c_1([L])^k}{k!},
\end{equation}
for $L$ a line bundle on $\Omega$ and $[L]$ its
class in $K$-theory, where $c_1([L])$ is the first Chern
class of $[L]$, and then extended to all virtual vector bundles using
the ``splitting principle'' (which says that one can pretend all
vector bundles split into direct sums of line bundles).
This is a ring homomorphism, and after tensoring the domain with $\bQ$,
it becomes a rational isomorphism of rings.
(Strictly speaking, we have only defined the Chern character on
even-degree $K$-theory, with values in even-degree cohomology,
but it extends by suspension to odd degrees as well.)
This isomorphism in fact is derived from the isomorphism
$\ch\co K^*(\Gamma)\otimes \bQ\xrightarrow{\cong}\check{H}^*(\Gamma;\bQ)$
at the AP complex level and
(\ref{eqn:dirLimits}). When the dimension of $\Gamma$ is at least four,
then because of the denominators in \eqref{eqn:chern},
$\ch$ does not necessarily come from an integral isomorphism between
$K^*(\Gamma)$ and $H^*(\Gamma;\bZ)$.  (It does give such an
isomorphism for spheres and tori --- see for example
\cite[Proposition 4.3]{VB} --- but not for $\bC\bP^k$, $k\ge 2$
\cite[\S4.1]{VB}.)

In \S \ref{sec:breakdown} we will
give an explicit example of how the Chern character can fail to be an integral
isomorphism for tiling spaces. This will give a negative answer to
a question of Benameur and Mathai \cite[p.\ 8]{MR4030282}.

\subsection{Frequencies, traces and the gap-labeling conjecture}
\label{subsec:RS}
Suppose $\cT$ is a primitive substitution tiling and
$P\subset \cT$ is a patch. Let us make some comments about patches.

First, the patch has a \textbf{frequency}: let
$\mathrm{freq}_R(P,\mathcal{T})$ be the number of copies of $P$ found
(up to equivalence)
as subpatches of $\mathcal{T}$ which are completely contained in a
ball of radius $R$ around the origin. Since the tiling came from a
primitive substitution, then the limit
$R^{-d}\mathrm{freq}_R(P,\mathcal{T})\rightarrow
\mathrm{freq}(P,\mathcal{T})\geq 0$ exists as $R\rightarrow \infty$,
and it is called the \textbf{patch frequency} of $P$. Since there is a
unique $\mathbb{R}^d$-invariant measure on $\Omega_\mathcal{T}$, then
this limit is actually independent of $\mathcal{T}$, meaning that the
patch frequency of the patch $P$ is the same for any tiling in
$\Omega$. This will be denoted by $\mathrm{freq}(P)$. 

Secondly, if $P$ is a patch, then $P$ is given by a closed subset
$\bar{P}\subset \Gamma$ and a $k\geq 0$ in that the set
$\pi^{-1}_k(\bar{P})\subset \Omega$ is isometric to $P\times
\mathcal{C}_P$, where $\mathcal{C}_P$ is a (transversal) Cantor
set. It follows from Birkhoff's ergodic theorem that
$\nu(\mathcal{C}_P) = \mathrm{freq}(P)$, and for this reason $\nu$ is
sometimes referred to as the \textbf{frequency measure}. 

Finally, if $P$ is a patch, then it defines a cohomology class $[P]\in
\check{H}^d(\Omega;\mathbb{Z})$ as follows. Recall that there exists a
closed set $\bar{P}\subset \Gamma$ such that
$\pi^{-1}_k(\bar{P})\subset \Omega$ is isometric to $P\times
\mathcal{C}_P$. Then the interior of $\bar{P}$ in $\Gamma$ is an open
set with a \v Cech cohomology class $[\bar{P}]\in
\check{H}^d(\Gamma;\mathbb{Z})$ which by (\ref{eqn:dirLimits}) defines
a class $[P]\in \check{H}^d(\Omega;\mathbb{Z})$. The map $[P]\mapsto
\mathrm{freq}(P) = \nu(\mathcal{C}_P)$ turns out to give a
homomorphism $C_\mu:\check H^d(\Omega;\mathbb{Z})\rightarrow
\mathbb{R}$ called the \textbf{Schwartzman asymptotic cycle} or the
\textbf{Ruelle-Sullivan current} \cite{KP:RS}. The image
$C_\mu(\check{H}^d(\Omega;\mathbb{Z}))\subset \mathbb{R}$ is called
the \textbf{frequency module}.

Associated with any substitution tiling space there are several
$C^*$-algebras that one can define (see \cite{KP:survey} for
details). The principal one is $\cA_{p}(\Omega)$, the
\textbf{unstable (punctured) algebra}. Another important one is
$\cA_{AF}(\Omega)$, the AF-algebra constructed from the
substitution matrix defining the tiling space $\Omega$. Since
the substitution is assumed to be primitive, it comes with a
unique tracial state
$\tau_{AF}\co K_0(\cA_{AF}(\Omega))\rightarrow \bR$. For
a primitive substitution tiling, there is also a unique tracial state
$\tau\co K_0(\cA_{p}(\Omega))\rightarrow \bR$ which
can be recovered from the unique $\bR^d$-invariant
probability measure $\mu$ on the tiling space $\Omega$. 

How the images $C_\mu(\check{H}^d(\Omega;\mathbb{Z}))$ and
$\tau(K_0(\mathcal{A}_p(\Omega)))$ overlap is a delicate matter.
The conjecture that they are equal was proposed by Bellissard
\cite{bellissard:K}, and went under the name of the
\textbf{gap-labeling conjecture}. There have been three papers claiming to
have proved this equality \cite{BBG:GLT,BOO:GLT,KP:GLT}; we now briefly
review what is at stake. As pointed out in \cite{ADRS:bloch}, several of
these papers invoke arguments which rely on the commutative diagram
\begin{equation}
\label{eq:forreal}
\begin{tikzcd}
K_0(\mathcal{A}_p(\Omega))\otimes \mathbb{Q} \arrow[r, "\mathrm{ch}_d\circ \chi^{-1}"] \arrow[d, "\tau"]
& \check{H}^d(\Omega;\mathbb{Q}) \arrow[d, "C_\mu"] \\
\mathbb{R} \arrow[r,  "\mathrm{Id}"]
& || \mathbb{R}
\end{tikzcd}
\end{equation}
where
$\chi:K^{-d}(\Omega)\rightarrow K_0(\mathcal{A}_p(\Omega))$ is an isomorphism,
$\mathrm{ch}:K^*(\Omega)\otimes \mathbb{Q}\rightarrow
\check{H}^*(\Omega;\mathbb{Z})$ is the isomorphism given by the
Chern character in (\ref{eqn:chern}), and $\mathrm{ch}_d$ is its
projection to $\check{H}^d(\Omega;\mathbb{Q})$. The arguments
used in several papers have the gap that they assume that
$\mathrm{ch}_d$ factors through $\check{H}^d(\Omega;\mathbb{Z})$
to give the diagram 
\begin{equation}
        \label{eqn:wishful}
        \begin{tikzcd}
K_0(\mathcal{A}_p(\Omega)) \arrow[r, "\mathrm{ch}_d\circ \chi^{-1}"] \arrow[d, "\tau"]
& \check{H}^d(\Omega;\mathbb{Z}) \arrow[d, "C_\mu"] \\
\mathbb{R} \arrow[r,  "\mathrm{Id}"]
& || \mathbb{R}
        \end{tikzcd}.
\end{equation}
Although this may hold in low dimensions, in section \S
\ref{sec:breakdown} we give an explicit example where this is not true, and
thus show that any argument to prove the gap-labeling conjecture using this
isomorphism is incomplete.

\input equivlabeling.tex

\section{Cubical substitutions and the cup product}
\label{sec:cubical}
In this section we describe our approach to compute the cohomology rings of
tiling spaces given by cubical substitutions. With an eye of doing
computations in dimension four in \S\ref{sec:breakdown}, we work out examples
in dimension two where we explicitly compute the cohomology ring structure
and other invariants.
We will first start describing in detail how the type of collaring
scheme used to compute the cohomology groups of cubical substitutions
and the ring structure is computed by following \cite{KM:ring}. Table
\ref{tab:1} summarizes the properties of the different
examples\footnote{The table lists groups as $\mathbb{Z}[1/16]$,
  $\mathbb{Z}[1/8]$, $\mathbb{Z}[1/4]$ and $\mathbb{Z}[1/2]$, even
  though these can all be written as $\mathbb{Z}[1/2]$. Following
  \cite[\S 3.5]{sadun:book}, we write them like this to emphasize
  their different scalings under the substitution rule.}. 

The four examples in this section are meant as two pairs of examples
to be compared. Of particular note are examples \ref{ex:2} and
\ref{ex:4}. These are two tiling spaces with isomorphic cohomology
groups over \(\mathbb{Z}\) and the same frequency module. Without the
cup product, they would be indistinguishable. However, the cup product
$\check{H}^1\times \check{H}^1\rightarrow \check{H}^2$ in one case is
surjective, whereas it is not in the other. This is the first example
we are aware of where the cup product is the only distinguishing
feature between two tiling spaces. Examples \ref{ex:1} and \ref{ex:3}
have the feature that they have indistinguishable cohomology spaces
over $\mathbb{Q}$, but can be distinguished either through their
cohomology groups over $\mathbb{Z}$ or through through the ring
structure. 

            \begin{table}[t]        
        \begin{center}
          \begin{tabular}{ |l|l|l|l| } 
            \hline
            Example & Cohomology groups & Cup product & Frequency module \\
            \hline
            \multirow{3}{4em}{ \ref{ex:1}}& $ \check{H}^2(\Omega)=\mathbb{Z}[1/25]\oplus\mathbb{Z}[1/15]\oplus\mathbb{Z}[1/5]\oplus\mathbb{Z}[1/3]$ & &  \\
            &$\check{H}^1(\Omega)=\mathbb{Z}[1/5]^2\oplus\mathbb{Z}[1/3]\oplus\mathbb{Z}$ & Surjective & $\frac{1}{6}\mathbb{Z}[1/5]$  \\
            &$\check{H}^0(\Omega)=\mathbb{Z}$ & & \\
            \hline
            \multirow{3}{4em}{ \ref{ex:2}}& $ \check{H}^2(\Omega)=\mathbb{Z}[1/16]\oplus\mathbb{Z}[1/8]\oplus\mathbb{Z}[1/4]\oplus\mathbb{Z}[1/2]$ & &  \\
            &$\check{H}^1(\Omega)=\mathbb{Z}[1/4]^2\oplus\mathbb{Z}[1/2]\oplus\mathbb{Z}$ & Surjective & $\frac{1}{3}\mathbb{Z}[1/2]$ \\
            &$\check{H}^0(\Omega)=\mathbb{Z}$ & & \\
            \hline
            \multirow{3}{4em}{ \ref{ex:3}}&  & &  \\
            & See Example \ref{ex:3} for details & Non-Surjective & $\frac{1}{84}\mathbb{Z}[1/5]$ \\
            & & & \\
            \hline
            \multirow{3}{4em}{ \ref{ex:4}}& $ \check{H}^2(\Omega)=\mathbb{Z}[1/16]\oplus\mathbb{Z}[1/8]\oplus\mathbb{Z}[1/4]\oplus\mathbb{Z}[1/2]$ & &  \\
            &$\check{H}^1(\Omega)=\mathbb{Z}[1/4]^2\oplus\mathbb{Z}[1/2]\oplus\mathbb{Z}$ & Non-Surjective & $\frac{1}{3}\mathbb{Z}[1/2]$ \\
            &$\check{H}^0(\Omega)=\mathbb{Z}$ & & \\
            \hline
          \end{tabular}
        \end{center}
        \caption{}
        \label{tab:1}
      \end{table}
      \subsection{Computational setup}
\indent Let us first establish the notations and conventions.\\
\indent Recall that a cubical substitution rule \(\varsigma\) on \(m\)
prototiles with (uniform) expansion matrix
\(\lambda=\lambda\cdot\textnormal{Id}_d\) has the set of prototiles 
\([0,1]^d\times\{0,\ldots,m-1\}\). (I.e., the set of colors
is $\Col=\{0,\ldots,m-1\}$.)  We will name each
prototile by its second coordinate, its \textbf{label}. In addition,
recall that each tile is a labeled region
\(t+[0,1]^d\subseteq\mathbb{R}^d\). We interpret this as being formed
by anchoring the point \([0]^d\in[0,1]^d\), its \textbf{puncture}, in
the corresponding prototile to \(t\in\mathbb{R}^d\). We follow the
convention of \cite{KM:ring} where $[x]$ denotes a (degenerate)
interval consisting of the single point $x$.\\ 
\indent We decompose the substitution rule \(\varsigma\) into the
following two steps: the first, denoted \(\lambda\), inflates each
tile by \(\lambda\); the second, denoted \(\sigma\), subdivides the
resulting region into tiles. If the initial tile is of the form
\((t,k)\) with \(t\in\mathbb{Z}^d\), the resulting substituted patch
will have each of the punctures of each of its tiles attached to
\(\mathbb{Z}^d\). Thus, let us put the lexicographic order on
\(\{0,\ldots,\lambda-1\}^d\), and using this order we abbreviate the
substitution rule as \(\varsigma(k)=(k_0,\ldots,k_{\lambda^d-1})\),
where \(0\leq k_i\leq m-1\) is a prototile with its puncture attached
at the coordinate that is the \(i\textnormal{th}\) entry in
\((0,\ldots,\lambda-1)^d\). We illustrate this convention through a
well-known example. 
\begin{example}[Chair tiling]
The standard chair tiling can be decomposed into squares, yielding the
MLD-equivalent cubical substitution rule in \(d=2\) 
\[
\begin{tikzpicture}
\draw(0,0)--(.5,0)--(.5,.5)--(0,.5)--(0,0);
\draw[->](.125,.125)--(.375,.375);
\draw[|->](.75,.25)--(1.25,.25);
\draw(1.5,-.25)--(2.5,-.25)--(2.5,.75)--(1.5,.75)--(1.5,-.25);
\draw(2,-.25)--(2,.75);
\draw(1.5,.25)--(2.5,.25);
\draw[->](1.625,-.125)--(1.875,.125);
\draw[->](1.625,.625)--(1.875,.375);
\draw[->](2.125,.375)--(2.375,.625);
\draw[->](2.375,-.125)--(2.125,.125);

\draw(3.5,.5)--(4.0,.5)--(4.0,0)--(3.5,0)--(3.5,.5);
\draw[->](3.625,.375)--(3.875,.125);
\draw[|->](4.25,.25)--(4.75,.25);
\draw(5.0,.75)--(6.0,.75)--(6.0,-.25)--(5.0,-.25)--(5.0,.75);
\draw(5.5,.75)--(5.5,-.25);
\draw(5.0,.25)--(6.0,.25);
\draw[->](5.125,.625)--(5.375,.375);
\draw[->](5.125,-.125)--(5.375,.125);
\draw[->](5.625,.125)--(5.875,-.125);
\draw[->](5.875,.625)--(5.625,.375);

\draw(.5,-2)--(0,-2)--(0,-1.5)--(.5,-1.5)--(.5,-2);
\draw[->](.375,-1.875)--(.125,-1.625);
\draw[|->](.75,-1.75)--(1.25,-1.75);
\draw(2.5,-2.25)--(1.5,-2.25)--(1.5,-1.25)--(2.5,-1.25)--(2.5,-2.25);
\draw(2,-2.25)--(2,-1.25);
\draw(2.5,-1.75)--(1.5,-1.75);
\draw[->](2.375,-2.125)--(2.125,-1.875);
\draw[->](2.375,-1.375)--(2.125,-1.625);
\draw[->](1.875,-1.625)--(1.625,-1.375);
\draw[->](1.625,-2.125)--(1.875,-1.875);

\draw(4.0,-1.5)--(3.5,-1.5)--(3.5,-2.0)--(4.0,-2.0)--(4.0,-1.5);
\draw[->](3.875,-1.625)--(3.625,-1.875);
\draw[|->](4.25,-1.75)--(4.75,-1.75);
\draw(6.0,-1.25)--(5.0,-1.25)--(5.0,-2.25)--(6.0,-2.25)--(6.0,-1.25);
\draw(5.5,-1.25)--(5.5,-2.25);
\draw(6.0,-1.75)--(5.0,-1.75);
\draw[->](5.875,-1.375)--(5.625,-1.625);
\draw[->](5.875,-2.125)--(5.625,-1.875);
\draw[->](5.375,-1.875)--(5.125,-2.125);
\draw[->](5.125,-1.375)--(5.375,-1.625);
\end{tikzpicture}.
\]
Labeling the prototiles \(\{0,1,2,3\}=\{\tikz{\draw(0,0)--(.3,0)--(.3,.3)--(0,.3)--(0,0);\draw[->](.075,.075)--(.225,.225);},
\tikz{\draw(0,0)--(.3,0)--(.3,.3)--(0,.3)--(0,0);\draw[->](.075,.225)--(.225,.075);},
\tikz{\draw(0,0)--(.3,0)--(.3,.3)--(0,.3)--(0,0);\draw[->](.225,.075)--(.075,.225);},
\tikz{\draw(0,0)--(.3,0)--(.3,.3)--(0,.3)--(0,0);\draw[->](.225,.225)--(.075,.075);}\}\), we abbreviate this substitution rule as
\begin{align*}
\varsigma(0)&=(0,1,2,0)\\
\varsigma(1)&=(0,1,1,3)\\
\varsigma(2)&=(0,2,2,3)\\
\varsigma(3)&=(3,1,2,3).
\end{align*}
or further abbreviated as a list in Sage as the following.
\begin{minted}[linenos,fontsize=\fontsize{5}{6},numbersep=5pt,frame=lines,framesep=2mm,mathescape=true,escapeinside=||]{python}
[[0,1,2,0],[0,1,1,3],[0,2,2,3],[3,1,2,3]]
\end{minted}
\indent For arbitrary dimension $d>1$, the $d$-dimensional chair substitution rule is the following, in Sage.
\begin{minted}[linenos,fontsize=\fontsize{5}{6},numbersep=5pt,frame=lines,framesep=2mm,mathescape=true,escapeinside=||]{python}
[[(2**d-1)-j if i==(2**d-1)-j else j for j in list(range(2**d))] for i in list(range(2**d))]
\end{minted}
\end{example}

One immediately encounters an obstacle that is very difficult to
overcome when trying to construct the \(AP\)-inverse limit and compute
its cohomology. The top-dimensional cells in the complex are collared
prototiles, or, for cubical substitutions, patches of size \(3^d\),
and the codimension-\(1\) cells are obtained from adjacency
information of pairs of collared tiles. The latter requires us to find
all possible patches of size \(4^d\) in \(\mathcal{T}\). In higher
dimensions, e.g. \(d=4\), this becomes unwieldy, since they are
proportional to \(m^{4d}\). In addition, it may require multiple
substitutions before the supertiles of that level observe all such
patches of size \(4^d\). In higher dimensions, these severely restrict
the substitution rules we can check. Thus, rather than the standard
collaring procedure yielding the \(AP\)-complex and inverse limit, let
us consider an alternative, the dual complex, introduced in
\cite{gahler:talk0113}, that is a variant of the Barge--Diamond
(\(BD\)) complex from \cite{bargediamond:bd08} for \(d=1\) and
\cite{bargediamondhuntonsadun:bd10} for \(d>1\).\\ 
\indent For motivation, let us first consider an alternate point of
view of how the \(AP\)-complex is formed. Suppose that we are given
the set of all collared prototiles, viewed as patches in
\(\mathcal{T}\). They form the set of top-dimensional cells. The set
of codimension-\(1\) cells consists of unions of pairs of neighboring
patches, viewed as intersections of their associated cylinder
sets. Similarly, the set of codimension-\(k\) cells consists of unions
of patches that share a common intersection. That is, the
\(AP\)-complex is the \v Cech complex associated associated to
collared prototiles. 

For cubical complexes, the top-dimensional cells are patches (up to
translation) of size \(3^d\), and the codimension-\(k\) cells are
unions of \(2^k\) patches that result in patches of \(\mathcal{T}\)
and share exactly \(2^k3^{d-k}\) tiles (so that the center tiles
neighbor along a codimension-\(k\) cell; analogously, so that the
differences are patches of size \(1\) along \(k\)-directions,
\(3^{d-k}\) along the remaining \(d-k\)-directions). Repeating this
process for each level of supertiles yields the \(AP\)-inverse
limit. From the \v Cech perspective, the fact that this inverse limit
is homeomorphic to \(\Omega_\mathcal{T}\) is because the basic open
sets of each (patches of collared supertiles for the AP point of view
compared to arbitrary patches) are mutual refinements. 

Rather than working with patches of size \(3^d\), one can use patches
of size \(2^d\), then carry through the \v Cech construction. The
resulting complex, due to \cite{gahler:talk0113}, is the \textbf{dual
  complex}. More precisely, at the initial level of the inverse limit,
the codimension-\(k\) cells are unions of \(2^k\) patches that share
exactly \(1^k2^{d-k}\) tiles. To see that this is collaring of sorts,
if a \(d\)-cell has the puncture of its first tile (in lexicographic
order) situated at the origin, one pretends that the ``center'' tile
(the actual \(d\)-cell) is the unit cube with its puncture at
\([1/2]^d\). 

\begin{figure}[t]
\centering
\begin{tikzpicture}
\draw(0,0)--(1,0)--(1,1)--(0,1)--(0,0);
\draw(.5,0)--(.5,1);
\draw(0,.5)--(1,.5);
\fill[draw opacity=0.5,fill opacity=0,pattern=north east lines](.25,.25)--(.75,.25)--(.75,.75)--(.25,.75)--(.25,.25);
\draw(.25,.25)--(.75,.25)--(.75,.75)--(.25,.75)--(.25,.25);
\draw[|->](1.25,.5)--(1.75,.5)node[midway,above]{\(\varsigma\)};
\draw[draw opacity=0.25](2,-.5)--(4,-.5)--(4,1.5)--(2,1.5)--(2,-.5);
\draw[draw opacity=0.25](2.5,-.5)--(2.5,1.5);
\draw[draw opacity=0.25](3,-.5)--(3,1.5);
\draw[draw opacity=0.25](3.5,-.5)--(3.5,1.5);
\draw[draw opacity=0.25](2,0)--(4,0);
\draw[draw opacity=0.25](2,.5)--(4,.5);
\draw[draw opacity=0.25](2,1)--(4,1);
\draw(2.25,-.25)--(3.75,-.25)--(3.75,1.25)--(2.25,1.25)--(2.25,-.25);
\draw(2.75,-.25)--(2.75,1.25);
\draw(3.25,-.25)--(3.25,1.25);
\draw(2.25,.25)--(3.75,.25);
\draw(2.25,.75)--(3.75,.75);
\fill[draw opacity=0.5,fill opacity=0,pattern=north east lines](2.50,0)--(3.00,0)--(3.00,.50)--(2.50,.50)--(2.50,0);
\fill[draw opacity=0.5,fill opacity=0,pattern=north east lines](3.00,0)--(3.50,0)--(3.50,.50)--(3.00,.50)--(3.00,0);
\fill[draw opacity=0.5,fill opacity=0,pattern=north east lines](2.50,.50)--(3.00,.50)--(3.00,1.00)--(2.50,1.00)--(2.50,.50);
\fill[draw opacity=0.5,fill opacity=0,pattern=north east lines](3.00,.50)--(3.50,.50)--(3.50,1.00)--(3.00,1.00)--(3.00,.50);
\draw(2.50,0)--(3.50,0)--(3.50,1)--(2.50,1)--(2.50,0);
\draw(2.50,.50)--(3.50,.50);
\draw(3.00,0)--(3.00,1.00);
\draw[|->](4.25,.5)--(4.75,.5)node[midway,above]{\(h\)}node[midway,below]{\(\simeq\)};
\draw[draw opacity=0.25](5,-.5)--(7,-.5)--(7,1.5)--(5,1.5)--(5,-.5);
\draw[draw opacity=0.25](5.5,-.5)--(5.5,1.5);
\draw[draw opacity=0.25](6,-.5)--(6,1.5);
\draw[draw opacity=0.25](6.5,-.5)--(6.5,1.5);
\draw[draw opacity=0.25](5,0)--(7,0);
\draw[draw opacity=0.25](5,.5)--(7,.5);
\draw[draw opacity=0.25](5,1)--(7,1);
\draw(5.00,-.50)--(6.50,-.50)--(6.50,1.00)--(5.00,1.00)--(5.00,-.50);
\draw(5.50,-.50)--(5.50,1.00);
\draw(6.00,-.50)--(6.00,1.00);
\draw(5.00,0)--(6.50,0);
\draw(5.00,.50)--(6.50,.50);
\fill[draw opacity=0.5,fill opacity=0,pattern=north east lines](5.25,-.25)--(5.75,-.25)--(5.75,.25)--(5.25,.25)--(5.25,-.25);
\fill[draw opacity=0.5,fill opacity=0,pattern=north east lines](5.75,-.25)--(6.25,-.25)--(6.25,.25)--(5.75,.25)--(5.75,-.25);
\fill[draw opacity=0.5,fill opacity=0,pattern=north east lines](5.25,.25)--(5.75,.25)--(5.75,.75)--(5.25,.75)--(5.25,.25);
\fill[draw opacity=0.5,fill opacity=0,pattern=north east lines](5.75,.25)--(6.25,.25)--(6.25,.75)--(5.75,.75)--(5.75,.25);
\draw(5.25,-.25)--(6.25,-.25)--(6.25,.75)--(5.25,.75)--(5.25,-.25);
\draw(5.25,.25)--(6.25,.25);
\draw(5.75,-.25)--(5.75,.75);
\end{tikzpicture}
\caption{A substitution rule with even expansion applied to a patch of
  size \(2^d\) representing a top-dimensional cell (center tile
  shaded) returns a patch (light gray lines) whose center tiles
  (shaded) has collars (black lines) that are not actual patches,
  therefore is not cellular. One needs to compose the substitution
  with a homotopy \(h\) to produce a cellular map.}
\label{figure:dual-even-expansion}
\end{figure}%

To construct the complex at the supertile level, we first note that if
\(\lambda\) is even, the substitution rule is not a cellular map
(Figure \ref{figure:dual-even-expansion}). Assuming that each patch of
size \(2^d\) has the puncture of its first tile located at the origin,
we define a supertile to be the collection of all subpatches of size
\(2^d\) (they will overlap) formed from tiles of the actual supertile
with punctures located at \((0,\ldots,\lambda)^d\). The substitution
rule then takes a supertile of the dual complex and sends it to all
subpatches of size \(2^d\). Higher-level supertiles are analogous. 

We observe that if we focus on the ``center'' tiles (and treat the
surrounding patches of size \(2^d\) as decorations), this map is the
ordinary substitution rule followed by a fixed homotopy, therefore the
maps of the
resulting inverse limit are stationary up to homotopy. 
Furthermore, the supertiles of
this construction and the supertiles of the Anderson--Putnam
construction are cofinal refinements of each other. This inverse limit
is thus shape equivalent to \(\Omega_\mathcal{T}\), and has the same \v Cech
cohomology groups as those of \(\Omega_\mathcal{T}\)
\cite[Theorem 16]{MR298634}.\footnote{We are indebted to the
referee for this observation.}

For all of the computations that are performed directly using Sage, we
will use the dual complex and the resulting inverse limit. 

\subsection{Distinguishing tiling spaces through cup product}
\label{sec:distinguishing}
\indent As an application of the ring structure of cohomology on cubical
substitutions, let us describe pairs of tiling spaces in \(d=2\) whose
cohomology groups coincide, but the ring structure differs. We consider
\(d=2\), since the required computations are much easier to visualize and
interpret. These examples are inspired by
\cite[Example 6.7]{solomyaktrevino:st23} due to the induced substitution
matrix on \(\check{H}^1\) containing two eigenvalues that do not come
from the expansion.\\
\indent We will consider two classes of examples. In all of the complexes that we construct, we take as positive left-to-right and down-to-up orientations of \(1\)-cells, and right-handed orientation of \(2\)-cells.
\begin{example}[Expansion \(5\) Product]
\label{ex:1}
  
Let us consider the following two one-dimensional substitutions on two and three prototiles, respectively.
\begin{minted}[linenos,fontsize=\fontsize{5}{6},numbersep=5pt,frame=lines,framesep=2mm,mathescape=true,escapeinside=||]{python}
|\(\varsigma_1\)|=[[0,0,0,0,1],[1,0,0,0,0]]
|\(\varsigma_2\)|=[[0,2,1,2,0],[0,1,1,1,0],[0,2,2,2,0]]
\end{minted}
One easily checks that both substitutions are primitive and recognizable.

To compute the first cohomology group, we use the Barge--Diamond
filtration from \cite{bargediamond:bd08}, where one successively blows
up lower-dimensional skeleta to become \(d\)-dimensional (here
\(d=1\)) complexes. We will do this in detail for \(\varsigma_1\),
then just state the main parts of the computation for
\(\varsigma_2\).\\ 
\indent Let \(\epsilon>0\). We replace the prototiles themselves,
assumed to be of length \(1\), by tiles of length \(1-2\epsilon\),
denoted using the original symbols, and replace the vertices that join
two neighboring tiles by new tiles of length \(2\epsilon\) (called
\textbf{vertex flaps} in \cite{bargediamondhuntonsadun:bd10}; more
generally, \(k\)-flaps for new tiles of size \(1-2\epsilon\) in
\(k\)-directions, \(2\epsilon\) in the remaining \(d-k\)-directions,
which are blow ups of \(k\)-cells in the tiling). These new prototiles
of length \(2\epsilon\) arise from the set of patches of length \(2\),
denoted by \(0.0\), \(0.1\), \(1.0\), and \(1.1\) for this
example. This new tiling is no longer self-similar, thus let us
construct a homotopy that maintains the same self-similar structure as
the original tiling.\\ 
\indent To maintain the same combinatorial structure as the original
substitution rule, we will require that the \(0\)-flaps be
nonexpanding, and the \(1\)-flaps (the slightly-shrunken prototiles)
to be fully expanding (more generally, for cubical substitutions,
\(k\)-flaps are expanding in \(k\)-directions, nonexpanding in the
remaining \(d-k\)-directions). More precisely, if \(\lambda\) is the
expansion of the original substitution rule, in dimension \(1\), the
\(0\)-flaps are substituted, via expansion by \(\lambda\) then a
homotopy, to prototiles of size \(2\epsilon\) (again \(0\)-flaps), and
the \(1\)-flaps are sent, via expansion by \(\lambda\) then a
homotopy, to patches of size \(\lambda-2\epsilon\) (a combination of
\(0\)- and \(1\)-flaps).\\ 
\indent For \(\varsigma_1\), we have that the induced substitution rule, still denoted \(\varsigma_1\), is
\[
\begin{tikzpicture}
\draw(-.750,0)--(.750,0);
\draw(0,0)node[below]{\(0\)};
\draw[|->](1.0,0)--(1.5,0);
\draw(1.750,0)--(7.250,0);
\draw(3.250,-.1)--(3.250,.1);
\draw(3.750,-.1)--(3.750,.1);
\draw(3.750,-.1)--(3.750,.1);
\draw(5.250,-.1)--(5.250,.1);
\draw(5.750,-.1)--(5.750,.1);
\draw(2.50,0)node[below]{\(0\)};
\draw(3.50,0)node[below]{\(0.0\)};
\draw(4.50,0)node[below]{\(0\)};
\draw(5.50,0)node[below]{\(0.1\)};
\draw(6.50,0)node[below]{\(1\)};

\draw(-.750,-1)--(.750,-1);
\draw(0,-1)node[below]{\(1\)};
\draw[|->](1.0,-1)--(1.5,-1);
\draw(1.750,-1)--(7.250,-1);
\draw(3.250,-1.1)--(3.250,-.9);
\draw(3.750,-1.1)--(3.750,-.9);
\draw(3.750,-1.1)--(3.750,-.9);
\draw(5.250,-1.1)--(5.250,-.9);
\draw(5.750,-1.1)--(5.750,-.9);
\draw(2.50,-1)node[below]{\(1\)};
\draw(3.50,-1)node[below]{\(1.0\)};
\draw(4.50,-1)node[below]{\(0\)};
\draw(5.50,-1)node[below]{\(0.0\)};
\draw(6.50,-1)node[below]{\(0\)};

\draw(.250,-2)--(0.750,-2);
\draw(0.50,-2)node[below]{\(0.0\)};
\draw[|->](1.0,-2)--(1.5,-2);
\draw(1.750,-2)--(2.250,-2);
\draw(2.0,-2)node[below]{\(1.0\)};

\draw(.250,-3)--(.750,-3);
\draw(.50,-3)node[below]{\(0.1\)};
\draw[|->](1.0,-3)--(1.5,-3);
\draw(1.750,-3)--(2.250,-3);
\draw(2.0,-3)node[below]{\(1.1\)};

\draw(.250,-4)--(.750,-4);
\draw(.50,-4)node[below]{\(1.0\)};
\draw[|->](1.0,-4)--(1.5,-4);
\draw(1.750,-4)--(2.250,-4);
\draw(2.0,-4)node[below]{\(0.0\)};

\draw(.250,-5)--(.750,-5);
\draw(.50,-5)node[below]{\(1.1\)};
\draw[|->](1.0,-5)--(1.5,-5);
\draw(1.750,-5)--(2.250,-5);
\draw(2.0,-5)node[below]{\(0.1\)};
\end{tikzpicture}.
\]
\indent The Barge--Diamond filtration is a filtration of complexes
\(S_0\subseteq S_1\subseteq\cdots\subseteq S_d\). \(S_0\) is formed by
taking the \(0\)-flaps and performing identifications similar to the
Anderson--Putnam construction, \(S_1\) is formed by taking \(S_0\),
then adding in the \(1\)-flaps (then possibly performing
identifications if there are \(2\)-flaps), etc. For \(\varsigma_1\),
\(S_0\), drawn in solid lines, is a loop 
\[
\begin{tikzpicture}[on top/.style={preaction={draw=white,-,line width=#1}},on top/.default=4pt]
\draw(0.5,0)--(0,-1);
\draw[on top](0.5,-1)--(0,0);
\draw(0,0)--(0.5,0);
\draw(0,-1)--(0.5,-1);
\draw[dashed](-1.5,0)--(0,0);
\draw[dashed](-1.5,-1)--(0,-1);
\draw[dashed](2,0)--(0.5,0);
\draw[dashed](2,-1)--(0.5,-1);
\draw(0.25,0)node[above]{\(0.0\)};
\draw(0.375,-0.75)node[above,rotate=-63.435]{\(0.1\)};
\draw(0.125,-0.75)node[above,rotate=63.435]{\(1.0\)};
\draw(0.25,-1)node[below]{\(1.1\)};
\draw(-0.75,0)node[above]{\(0\)};
\draw(1.25,0)node[above]{\(0\)};
\draw(-0.75,-1)node[below]{\(1\)};
\draw(1.25,-1)node[below]{\(1\)};
\end{tikzpicture},
\]
where the dashed lines indicate the \(1\)-flaps that yield, for example, identification of the left endpoints of \(0.0\) and \(0.1\). The complex \(S_1\) is then the same picture, together with the dashed lines turned solid. We denote \(\Xi_i=\varprojlim(S_i,\varsigma_1)\). Then \(\Xi_1=\Omega_{\varsigma_1}\).\\
\indent One can check that \(\varsigma_1:S_0\rightarrow S_0\) is an orientation-reversing self-homeomorphism, thus \(\check{H}^1(\Xi_0)=\mathbb{Z}\). For \(\Xi_1\), we use the long exact sequence in relative cohomology of a pair \((\Xi_1,\Xi_0)\)
\[
\begin{tikzcd}
0\arrow[r]&\check{H}^0(\Xi_1,\Xi_0)\arrow[r]&\check{H}^0(\Xi_1)\arrow[r]\arrow[d,phantom,""{coordinate,name=X}]&\check{H}^0(\Xi_0)\arrow[lld,rounded corners,to path={--([xshift=2ex]\tikztostart.east)|-(X.center)\tikztonodes-|([xshift=-2ex]\tikztotarget.west)--(\tikztotarget)},swap,"\delta" at end]&\\
&\check{H}^1(\Xi_1,\Xi_0)\arrow[r]&\check{H}^1(\Xi_1)\arrow[r]&\check{H}^1(\Xi_0)\arrow[r]&0
\end{tikzcd}.
\]
Since \(S_0\) has a single component, over reduced cohomology, the bottom row becomes a short exact sequence.\\
\indent \(S_1/S_0\) is a wedge of two circles, with its inverse limit computed exactly from the original substitution matrix. Thus
\[
\check{H}^1(\Xi_1,\Xi_0)=\varinjlim\left(\mathbb{Z}^2,\left(\begin{array}{cc}
4&1\\
4&1
\end{array}\right)\right)=\mathbb{Z}[1/5].
\]
Since \(\check{H}^1(\Xi_0)=\mathbb{Z}\) is free, the bottom short exact sequence splits, yielding
\begin{align*}
\check{H}^1(\Omega_{\varsigma_1})&=\mathbb{Z}[1/5]\oplus\mathbb{Z}\\
\check{H}^0(\Omega_{\varsigma_1})&=\mathbb{Z}.
\end{align*}
\indent The calculation for \(\varsigma_2\) is much easier. It is a proper substitution, thus forces the border (see Remark \ref{rem:force}), and
\[
\check{H}^1(\Omega_{\varsigma_2})=\varinjlim\left(\mathbb{Z}^3,\left(\begin{array}{ccc}2&1&2\\2&3&0\\2&0&3\end{array}\right)\right).
\]
The matrix has eigenvectors
$\mathbf{e}_5=(1,1,1)$ and
$\mathbf{e}_3=(0,2,-1)$
with the subscripts their respective eigenvalues. Direct computation then shows that the direct limit splits as a direct sum, giving us
$$\check{H}^1(\Omega_{\varsigma_2})=\mathbb{Z}[1/5]\oplus\mathbb{Z}[1/3] \hspace{.7in}\mbox{ and }\hspace{.7in} \check{H}^0(\Omega_{\varsigma_2})=\mathbb{Z}.$$
\indent Let us also compute the induced substitution matrices on \(\check{H}^1(\Omega_{\varsigma_1})\) using the dual complex. It is computed via direct limit of the matrix
$$ \sigma_1^1=\left(\begin{array}{ccc}
1&1&0\\
3&1&1\\
1&1&0
\end{array}\right)$$
with eigenvalues of \(3,-1\). The eigenvalue of \(-1\) matches the orientation reversal of \(\varsigma_1\) on its respective \(S_0\).\\
\indent We finally consider the product tiling space \(\Omega_{\varsigma_1}\times\Omega_{\varsigma_2}\). A direct application of K\"unneth theorem yields that
\begin{align*}
\check{H}^2(\Omega_{\varsigma_1}\times\Omega_{\varsigma_2})&=\mathbb{Z}[1/25]\oplus\mathbb{Z}[1/15]\oplus\mathbb{Z}[1/5]\oplus\mathbb{Z}[1/3]\\
\check{H}^1(\Omega_{\varsigma_1}\times\Omega_{\varsigma_2})&=\mathbb{Z}[1/5]^2\oplus\mathbb{Z}[1/3]\oplus\mathbb{Z}\\
\check{H}^0(\Omega_{\varsigma_1}\times\Omega_{\varsigma_2})&=\mathbb{Z}
\end{align*}
with \(\check{H}^2(\Omega_{\varsigma_1}\times\Omega_{\varsigma_2})\) having eigenvalues \(25,15,-5,-3\) and \(\check{H}^1(\Omega_{\varsigma_1}\times\Omega_{\varsigma_2})\) having eigenvalues \(5_2,5_1,3_2,-1_1\), with the subscripts indicating the \(1\)-dimensional substitution rule that gives rise to the corresponding eigenvalue.\\
\indent Lifting each of the cohomology classes in \(\check{H}^1(\Omega_{\varsigma_1}\times\Omega_{\varsigma_2})\) to \(1\)-cochains, applying the cubical cup product formula given in \cite{KM:ring}, then projecting to cohomology classes in \(\check{H}^2(\Omega_{\varsigma_1}\times\Omega_{\varsigma_2})\) gives that, in the order given in the list of eigenvalues, the cup product is the resulting bilinear form
$$\check{H}^1\times \check{H}^1\rightarrow \check{H}^2:(a,b)\mapsto B(a,b) := \sum_{*\in\{25,15,-5,-3\}}B_*(a,b)\cdot e_*,$$ where $e_n\in \check{H}^2$ is a generator with eigenvalue $n$, and $B_n$ is the bilinear form given by the matrices
\begin{eqnarray*}
B_{25}&:=\left(\begin{array}{cccc}
0&-1&0&0\\
1&0&0&0\\
0&0&0&0\\
0&0&0&0
\end{array}\right),\hspace{1in}
B_{15}&:=\left(\begin{array}{cccc}
0&0&0&0\\
0&0&1&0\\
0&-1&0&0\\
0&0&0&0
\end{array}\right),\\
B_{-5}&:=\left(\begin{array}{cccc}
0&0&0&-1\\
0&0&0&0\\
0&0&0&0\\
1&0&0&0
\end{array}\right),\hspace{1in}
B_{-3}&:=\left(\begin{array}{cccc}
0&0&0&0\\
0&0&0&0\\
0&0&0&-1\\
0&0&1&0
\end{array}\right).
\end{eqnarray*}
The important observation is that, as expected, there are always
cohomology classes in
\(\check{H}^1(\Omega_{\varsigma_1}\times\Omega_{\varsigma_2};\mathbb{Z})\)
cupping to any of the cohomology classes in
\(\check{H}^2(\Omega_{\varsigma_1}\times\Omega_{\varsigma_2};\mathbb{Z})\).\\ 
\indent Lastly, let us compute the frequency module. The substitution
matrix on the \(2\)-cells of the dual complex is the matrix
$\sigma^2=\sigma_{A1}$ found in Appendix \ref{app:1}, which has
Perron--Frobenius eigenvector 
\[
\resizebox{0.9\hsize}{!}{
\((95,19,25,5,25,5,5,1,76,20,20,4,38,10,10,2,38,10,10,2,76,20,20,4,19,5,5,1,38,76,10,20,10,20,2,4)\)
}
\]
Its sum is \(750\), giving us that the frequency module is
\(\frac{1}{6}\mathbb{Z}[1/5]\). For illustration, we also compute the
frequency modules for each of the one-dimensional substitutions. The
substitution matrices on the \(1\)-cells of the corresponding dual
complexes are 
$$
\sigma_1^1=\left(\begin{array}{*{4}c}
3&3&4&3\\
1&1&0&1\\
1&0&1&1\\
0&1&0&0
\end{array}\right)\hspace{.7in}\mbox{ and }\hspace{.7in}
\sigma_2^1=\left(\begin{array}{*{9}c}
1&1&1&1&1&1&1&1&1\\
0&0&0&1&0&0&1&1&0\\
1&1&1&0&1&1&0&0&1\\
0&0&0&2&0&0&2&2&0\\
1&1&1&0&0&0&0&0&0\\
0&0&0&0&2&2&0&0&2\\
0&0&0&1&0&0&1&1&0\\
1&1&1&0&0&0&0&0&0\\
1&1&1&0&1&1&0&0&1
\end{array}\right)
$$
with Perron--Frobenius eigenvectors that are duals of lifts of the generators of \(\check{H}^1(\Omega_{\varsigma_i})\) are $(19,5,5,1)$ and $(5,1,4,2,2,4,1,2,4)$
 with sums \(30\) and \(25\), respectively, giving frequency modules \(\frac{1}{6}\mathbb{Z}[1/5]\) and \(\mathbb{Z}[1/5]\), respectively. As expected, since the Ruelle--Sullivan map is a ring homomorphism, \(\frac{1}{6}\mathbb{Z}[1/5]\cdot\mathbb{Z}[1/5]=\frac{1}{6}\mathbb{Z}[1/5]\).
\end{example}
\begin{example}[Expansion \(4\) Product]
\label{ex:2}
  
Let us consider the following two one-dimensional substitutions on two prototiles.
\begin{minted}[linenos,fontsize=\fontsize{5}{6},numbersep=5pt,frame=lines,framesep=2mm,mathescape=true,escapeinside=||]{python}
|\(\varsigma_1\)|=[[0,0,0,1],[0,1,1,1]]
|\(\varsigma_2\)|=[[0,1,1,0],[1,0,0,1]]
\end{minted}
One, again, easily checks that both substitutions are primitive and recognizable. In particular, \(\varsigma_2\) is the (square of the) Thue--Morse substitution, whose cohomology groups were computed in \cite{bargediamond:bd08}, being
$$ \check{H}^1(\Omega_{\varsigma_2})=\mathbb{Z}[1/4]\oplus\mathbb{Z}\hspace{.7in}\mbox{ and }\hspace{.7in} \check{H}^0(\Omega_{\varsigma_2})=\mathbb{Z}.$$
\indent \(\varsigma_1\) is a proper substitution, thus forces the border (see Remark \ref{rem:force}), and
\[
\check{H}^1(\Omega_{\varsigma_1})=\varinjlim\left(\mathbb{Z}^2,\left(\begin{array}{cc}3&1\\1&3\end{array}\right)\right).
\]
The matrix has eigenvectors
$\mathbf{e}_4=\left(\begin{array}{cc}1\\1\end{array}\right)$ and
$\mathbf{e}_2=\left(\begin{array}{cc}1\\-1\end{array}\right)$
with the subscripts their respective eigenvalues. Direct computation then shows that the direct limit splits as a direct sum, giving us
$$\check{H}^1(\Omega_{\varsigma_1})=\mathbb{Z}[1/4]\oplus\mathbb{Z}[1/2] \hspace{.7in}\mbox{ and }\hspace{.7in} \check{H}^0(\Omega_{\varsigma_1})=\mathbb{Z}.$$
The K\"unneth theorem yields that, for the product \(\Omega_{\varsigma_1}\times\Omega_{\varsigma_2}\),
\begin{align*}
\check{H}^2(\Omega_{\varsigma_1}\times\Omega_{\varsigma_2})&=\mathbb{Z}[1/16]\oplus\mathbb{Z}[1/8]\oplus\mathbb{Z}[1/4]\oplus\mathbb{Z}[1/2]\\
\check{H}^1(\Omega_{\varsigma_1}\times\Omega_{\varsigma_2})&=\mathbb{Z}[1/4]^2\oplus\mathbb{Z}[1/2]\oplus\mathbb{Z}\\
\check{H}^0(\Omega_{\varsigma_1}\times\Omega_{\varsigma_2})&=\mathbb{Z}
\end{align*}
with \(\check{H}^2(\Omega_{\varsigma_1}\times\Omega_{\varsigma_2})\)
having eigenvalues \(16,8,4,2\) and
\(\check{H}^1(\Omega_{\varsigma_1}\times\Omega_{\varsigma_2})\) having
eigenvalues \(4_2,4_1,2_1,1_2\), with the subscripts indicating the
\(1\)-dimensional substitution rule that gives rise to the
corresponding eigenvalue.

The same procedure as the previous example gives that, in the order given in
the list of eigenvalues, the cup product is the resulting bilinear form
\[
\check{H}^1\times \check{H}^1\rightarrow \check{H}^2\co
(a,b)\mapsto B(a,b) := \sum_{*\in \{16,8,4,2\}}B_{*}(a,b)\cdot e_{*},
\]
where $e_n\in \check{H}^2$ is a generator with eigenvalue $n$, and $B_n$ is the bilinear form given by the matrices
\begin{eqnarray*}
B_{16}&:= \left(\begin{array}{cccc}
0&-1&0&0\\
1&0&0&0\\
0&0&0&0\\
0&0&0&0
\end{array}\right),\hspace{1in}
B_8 &:= \left(\begin{array}{cccc}
0&0&-1&0\\
0&0&0&0\\
1&0&0&0\\
0&0&0&0
\end{array}\right),\\
B_4 &:= \left(\begin{array}{cccc}
0&0&0&0\\
0&0&0&1\\
0&0&0&0\\
0&-1&0&0
\end{array}\right),\hspace{1in}
B_2&:=\left(\begin{array}{cccc}
0&0&0&0\\
0&0&0&0\\
0&0&0&1\\
0&0&-1&0
\end{array}\right).
\end{eqnarray*}
Again, we observe that, as expected, there are always cohomology classes in
\(\check{H}^1(\Omega_{\varsigma_1}\times\Omega_{\varsigma_2})\) cupping to any of
the cohomology classes in
\(\check{H}^2(\Omega_{\varsigma_1}\times\Omega_{\varsigma_2})\).\\
\indent Lastly, let us compute the frequency module. The substitution matrix on the \(2\)-cells of the dual complex is $\sigma^2 = \sigma_{A2}$ found in Appendix \ref{app:1} which has Perron--Frobenius eigenvector
\[
(1,2,1,2,1,2,1,2,2,2,2,2,1,1,1,1).
\]
Its sum is \(24\), giving us that the frequency module is \(\frac{1}{3}\mathbb{Z}[1/2]\). For illustration, we also compute the frequency modules for each of the one-dimensional substitutions. The substitution matrices on the \(1\)-cells of the corresponding dual complexes are
$$
\sigma_1^1=\left(\begin{array}{*{4}c}
2&2&0&0\\
1&1&1&1\\
1&1&1&1\\
0&0&2&2
\end{array}\right)\hspace{.6in}\mbox{ and } \hspace{.6in}
\sigma_2^1=\left(\begin{array}{*{4}c}
1&0&1&1\\
1&2&1&1\\
1&1&2&1\\
1&1&0&1
\end{array}\right)
$$
with Perron--Frobenius eigenvectors that are duals of lifts of the
generators of \(\check{H}^1(\Omega_{\varsigma_i})\) are $(1,1,1,1)$
and $(1,2,2,1)$ with sums \(4\) and \(6\), respectively, giving
frequency modules \(\mathbb{Z}[1/2]\) and
\(\frac{1}{3}\mathbb{Z}[1/2]\), respectively. As expected, since the
Ruelle--Sullivan map is a ring homomorphism,
\(\mathbb{Z}[1/2]\cdot\frac{1}{3}\mathbb{Z}[1/2]=\frac{1}{3}\mathbb{Z}[1/2]\). 
\end{example}
\begin{figure}[t]
\centering
\begin{tikzpicture}
\draw(0,0)--(.5,0)--(.5,.5)--(0,.5)--(0,0);
\draw[|->](.75,.25)--(1.25,.25);
\draw(.25,.25)node{\(0\)};
\draw(1.5,-1)--(4,-1)--(4,1.5)--(1.5,1.5)--(1.5,-1);
\draw(1.5,-.5)--(4,-.5);
\draw(1.5,0)--(4,0);
\draw(1.5,.5)--(4,.5);
\draw(1.5,1)--(4,1);
\draw(2,-1)--(2,1.5);
\draw(2.5,-1)--(2.5,1.5);
\draw(3,-1)--(3,1.5);
\draw(3.5,-1)--(3.5,1.5);
\draw(1.75,-.75)node{\(2\)};
\draw(1.75,-.25)node{\(2\)};
\draw(1.75,.25)node{\(2\)};
\draw(1.75,.75)node{\(2\)};
\draw(1.75,1.25)node{\(2\)};
\draw(2.25,-.75)node{\(0\)};
\draw(2.25,-.25)node{\(0\)};
\draw(2.25,.25)node{\(0\)};
\draw(2.25,.75)node{\(1\)};
\draw(2.25,1.25)node{\(1\)};
\draw(2.75,-.75)node{\(1\)};
\draw(2.75,-.25)node{\(0\)};
\draw(2.75,.25)node{\(0\)};
\draw(2.75,.75)node{\(1\)};
\draw(2.75,1.25)node{\(1\)};
\draw(3.25,-.75)node{\(0\)};
\draw(3.25,-.25)node{\(1\)};
\draw(3.25,.25)node{\(0\)};
\draw(3.25,.75)node{\(1\)};
\draw(3.25,1.25)node{\(1\)};
\draw(3.75,-.75)node{\(1\)};
\draw(3.75,-.25)node{\(0\)};
\draw(3.75,.25)node{\(0\)};
\draw(3.75,.75)node{\(0\)};
\draw(3.75,1.25)node{\(0\)};

\draw(5,0)--(5.5,0)--(5.5,.5)--(5,.5)--(5,0);
\draw[|->](5.75,.25)--(6.25,.25);
\draw(5.25,.25)node{\(1\)};
\draw(6.5,-1)--(9,-1)--(9,1.5)--(6.5,1.5)--(6.5,-1);
\draw(6.5,-.5)--(9,-.5);
\draw(6.5,0)--(9,0);
\draw(6.5,.5)--(9,.5);
\draw(6.5,1)--(9,1);
\draw(7,-1)--(7,1.5);
\draw(7.5,-1)--(7.5,1.5);
\draw(8,-1)--(8,1.5);
\draw(8.5,-1)--(8.5,1.5);
\draw(6.75,-.75)node{\(2\)};
\draw(6.75,-.25)node{\(2\)};
\draw(6.75,.25)node{\(2\)};
\draw(6.75,.75)node{\(2\)};
\draw(6.75,1.25)node{\(2\)};
\draw(7.25,-.75)node{\(0\)};
\draw(7.25,-.25)node{\(0\)};
\draw(7.25,.25)node{\(1\)};
\draw(7.25,.75)node{\(0\)};
\draw(7.25,1.25)node{\(0\)};
\draw(7.75,-.75)node{\(0\)};
\draw(7.75,-.25)node{\(0\)};
\draw(7.75,.25)node{\(1\)};
\draw(7.75,.75)node{\(1\)};
\draw(7.75,1.25)node{\(1\)};
\draw(8.25,-.75)node{\(0\)};
\draw(8.25,-.25)node{\(0\)};
\draw(8.25,.25)node{\(0\)};
\draw(8.25,.75)node{\(1\)};
\draw(8.25,1.25)node{\(0\)};
\draw(8.75,-.75)node{\(0\)};
\draw(8.75,-.25)node{\(1\)};
\draw(8.75,.25)node{\(0\)};
\draw(8.75,.75)node{\(0\)};
\draw(8.75,1.25)node{\(0\)};

\draw(10,0)--(10.5,0)--(10.5,.5)--(10,.5)--(10,0);
\draw[|->](10.75,.25)--(11.25,.25);
\draw(10.25,.25)node{\(2\)};
\draw(11.5,-1)--(14,-1)--(14,1.5)--(11.5,1.5)--(11.5,-1);
\draw(11.5,-.5)--(14,-.5);
\draw(11.5,0)--(14,0);
\draw(11.5,.5)--(14,.5);
\draw(11.5,1)--(14,1);
\draw(12,-1)--(12,1.5);
\draw(12.5,-1)--(12.5,1.5);
\draw(13,-1)--(13,1.5);
\draw(13.5,-1)--(13.5,1.5);
\draw(11.75,-.75)node{\(0\)};
\draw(11.75,-.25)node{\(0\)};
\draw(11.75,.25)node{\(0\)};
\draw(11.75,.75)node{\(1\)};
\draw(11.75,1.25)node{\(0\)};
\draw(12.25,-.75)node{\(2\)};
\draw(12.25,-.25)node{\(2\)};
\draw(12.25,.25)node{\(2\)};
\draw(12.25,.75)node{\(2\)};
\draw(12.25,1.25)node{\(2\)};
\draw(12.75,-.75)node{\(2\)};
\draw(12.75,-.25)node{\(2\)};
\draw(12.75,.25)node{\(2\)};
\draw(12.75,.75)node{\(2\)};
\draw(12.75,1.25)node{\(2\)};
\draw(13.25,-.75)node{\(2\)};
\draw(13.25,-.25)node{\(2\)};
\draw(13.25,.25)node{\(2\)};
\draw(13.25,.75)node{\(2\)};
\draw(13.25,1.25)node{\(2\)};
\draw(13.75,-.75)node{\(2\)};
\draw(13.75,-.25)node{\(2\)};
\draw(13.75,.25)node{\(2\)};
\draw(13.75,.75)node{\(2\)};
\draw(13.75,1.25)node{\(2\)};
\end{tikzpicture}
\label{figure:non-product5}
\caption{Substitution rule for an expansion \(5\) two-dimensional cubical substitution that has the cohomology groups of a product, but does not have the ring structure of a product.}
\end{figure}%

\begin{example}[Expansion \(5\) Non-product]
  \label{ex:3}
We consider the following two-dimensional substitution on three prototiles.
\begin{minted}[linenos,fontsize=\fontsize{5}{6},numbersep=5pt,frame=lines,framesep=2mm,mathescape=true,escapeinside=||]{python}
|\(\varsigma\)|=[[2,2,2,2,2,0,0,0,1,1,1,0,0,1,1,0,1,0,1,1,1,0,0,0,0],[2,2,2,2,2,0,0,1,0,0,0,0,1,1,1,0,0,0,1,0,0,1,0,0,0],[0,0,0,1,0,2,2,2,2,2,2,2,2,2,2,2,2,2,2,2,2,2,2,2,2]]
\end{minted}
This does not force the border (e.g. the horizontal \(1\)-cell containing the prototiles \(0,1\) on its bottom and prototiles \(0,1\) on its top stays branched regardless of the power of the substitution), but is primitive and recognizable.

Using the dual complex, \(\check{H}^2(\Omega_\varsigma)\) and \(\check{H}^1(\Omega_\varsigma)\) are computed via direct limits under the matrices
\begin{align*}
\sigma^2&=\left(\begin{array}{*{15}c}
1&0&3&1&2&1&1&2&0&1&1&0&2&1&2\\
1&0&3&1&2&1&1&2&0&1&1&0&2&1&2\\
2&1&3&1&3&2&2&2&0&1&1&20&2&1&2\\
5&4&6&2&5&5&5&1&2&2&1&40&3&1&2\\
1&0&3&1&2&1&1&2&0&1&1&0&2&1&2\\
6&3&9&3&5&5&5&1&3&2&2&0&4&1&2\\
-3&-1&-3&-1&-2&-2&-2&0&-1&0&-1&0&-1&0&0\\
-3&-1&-3&-1&-2&-2&-2&0&-1&0&-1&0&-1&0&0\\
1&0&3&1&2&1&1&2&0&1&1&0&2&1&2\\
3&1&3&1&2&2&2&0&1&0&1&0&1&0&0\\
6&2&6&2&4&4&4&0&2&0&2&0&2&0&0\\
0&0&3&1&0&0&0&0&0&0&0&15&0&0&0\\
3&2&6&2&3&3&3&1&2&2&1&0&3&1&2\\
6&3&9&3&5&5&5&1&3&2&2&0&4&1&2\\
3&1&3&1&2&2&2&0&1&0&1&0&1&0&0
\end{array}\right)\\
\sigma^1&=\left(\begin{array}{cccc}
5&0&0&0\\
0&3&1&0\\
0&4&1&4\\
0&0&1&3
\end{array}\right)
\end{align*}
with eigenvectors
\begin{align*}
\mathbf{e}_{5_1}^1&=(1,0,0,0)   &\mathbf{e}_{25}^2&=(1,1,2,4,1,3,-1,-1,1,1,2,1,2,3,1)\\
\mathbf{e}_{5_2}^1&=(0,1,2,1) &\mathbf{e}_{15}^2&=(1,1,0,0,1,3,-1,-1,1,1,2,-1,2,3,1)\\
\mathbf{e}_3^1&=(0,1,0,-1) &\mathbf{e}_{-5}^2&=(1,1,-4,-8,1,3,-1,-1,1,1,2,1,2,3,1)\\
\mathbf{e}_{-1}^1&=(0,1,-4,1) &\mathbf{e}_{-3}^2&=\left(\begin{array}{c}101,101,90,-72,101,-201,151,151,\\101,-151,-302,-11,-50,-201,-151\end{array}\right)
\end{align*}
where the subscripts indicate the corresponding eigenvalues. The cup product is the bilinear form
$$\check{H}^1\times \check{H}^1\rightarrow \check{H}^2:(a,b)\mapsto B(a,b) := \sum_{*\in \{25,15,-5,-3\}}B_{*}(a,b)\cdot \mathbf{e}_{*}^2,$$ where $\mathbf{e}_{n}^2\in \check{H}^2$ is a generator with eigenvalue $n$, and $B_n$ is the bilinear form given by the matrices
\begin{align*}
B_{25}&:=\left(\begin{array}{cccc}
0&-1&0&0\\
1&0&0&0\\
0&0&0&0\\
0&0&0&0
\end{array}\right),
& B_{15}&:=\left(\begin{array}{cccc}
0&0&-1&0\\
0&0&0&0\\
1&0&0&0\\
0&0&0&0
\end{array}\right), \\
B_{-5} &:= \left(\begin{array}{cccc}
0&0&0&-1\\
0&0&0&0\\
0&0&0&0\\
1&0&0&0
\end{array}\right),
& B_{-3}&:=\left(\begin{array}{cccc}
0&0&0&0\\
0&0&0&0\\
0&0&0&0\\
0&0&0&0
\end{array}\right).
\end{align*}
The last component of the bilinear form shows that there are no cohomology classes in \(\check{H}^1(\Omega_\varsigma)\) cupping to \(\mathbf{e}_{-3}^2\) in \(\check{H}^2(\Omega_\varsigma)\), contrasting the example coming from the product. This coincides with the fact that, up to a coboundary, the lifts of the eigenvectors \(\mathbf{e}_1^1\) and \(\mathbf{e}_{-1}^1\) to \(1\)-cochains are
\begin{align*}
\widetilde{\mathbf{e}_3^1}&=(0,0,0,0,0,1,1,0,0,1,0,0,-1,1)\\
\widetilde{\mathbf{e}_{-1}^1}&=(0,0,0,0,0,1,1,-4,-4,1,0,0,1,1)
\end{align*}
with the first five entries vertical \(1\)-cells, the remaining horizontal \(1\)-cells, therefore cannot cup nontrivially (to \(\mathbf{e}_{-3}^2\)).\\
\indent At this point, we have yet to calculate the actual cohomology groups of \(\Omega_\varsigma\) over \(\mathbb{Z}\). We follow the two-dimensional analogue of the Barge--Diamond filtration presented in \cite{bargediamondhuntonsadun:bd10}. For the sake of brevity, and since the calculations will be done in detail in Example \ref{ex:4}, here we skip most of the details and just state the conclusion.\\
\indent The long exact sequence in relative cohomology of a pair \((\Xi_1,\Xi_0)\) yields the short exact sequence
\[
\begin{tikzcd}
0\arrow[r]&\check{H}^1(\Xi_1,\Xi_0)\arrow[r]&\check{H}^1(\Xi_1)\arrow[r]&\check{H}^1(\Xi_0)\arrow[r]&0
\end{tikzcd}
\]
where
\[
\check{H}^1(\Xi_1,\Xi_0)=\varinjlim\left(\mathbb{Z}^3,\left(\begin{array}{*{3}c}
5&0&0\\
0&4&1\\
0&1&4
\end{array}\right)\right),
\]
and of \((\Xi_2,\Xi_1)\) yields
\[
\begin{tikzcd}
0\arrow[r]&\check{H}^1(\Xi_2,\Xi_1)\arrow[r]&\check{H}^1(\Xi_2)\arrow[r]\arrow[d,phantom,""{coordinate,name=Y}]&\check{H}^1(\Xi_1)\arrow[lld,rounded corners,to path={--([xshift=2ex]\tikztostart.east)|-(Y.center)\tikztonodes-|([xshift=-2ex]\tikztotarget.west)--(\tikztotarget)},swap,"\delta_1^1" at end]&\\
&\check{H}^2(\Xi_2,\Xi_1)\arrow[r]&\check{H}^2(\Xi_2)\arrow[r]&\check{H}^2(\Xi_1)\arrow[r]&0
\end{tikzcd}
\]
where
\[
\check{H}^2(\Xi_2,\Xi_1)=\varinjlim\left(\mathbb{Z}^3,\left(\begin{array}{*{3}c}
11&9&5\\
14&6&5\\
4&1&20
\end{array}\right)\right).
\]
The eigenvalues of these matrices are \(5,3\) and \(25,15,-3\), respectively. A short calculation gives that \(\check{H}^1(\Xi_0)=\mathbb{Z}\). Thus \(\check{H}^1(\Xi_1,\Xi_0)\) is a direct summand of \(\check{H}^1(\Xi_1)\). By Schur's lemma, \(\delta_1^1=0\). Furthermore, since \(\Xi_2\) is a wedge of \(2\)-spheres, \(\check{H}^1(\Xi_2,\Xi_1)=0\). Therefore \(\check{H}^1(\Omega_\varsigma)=\check{H}^1(\Xi_2)=\check{H}^1(\Xi_1)\), which contains \(\check{H}^1(\Xi_1,\Xi_0)\) as a direct summand.\\
\indent However, the direct limit of the matrix
\[
\sigma=\left(\begin{array}{cc}4&1\\1&4\end{array}\right)
\]
\emph{does not split as a direct sum \(\mathbb{Z}[1/5]\oplus\mathbb{Z}[1/3]\)\footnote{If one has the matrix \(\sigma=\left(\begin{array}{cc}2&1\\1&2\end{array}\right)\) instead, by \cite[Exercise 7.5.2 (a)]{lindmarcus:lm21}, the direct limit does not split as a direct sum of the eigenspaces. However, it is still a direct sum \(\mathbb{Z}[1/3]\oplus\mathbb{Z}\) with \(\mathbb{Z}[1/3]\) generated by the eigenvector with eigenvalue \(3\), since \((\varinjlim(\mathbb{Z}^2,\sigma))/\mathbb{Z}[1/3]=\mathbb{Z}\), and the short exact sequence
\[
\begin{tikzcd}[ampersand replacement=\&]
0\arrow[r]\&\mathbb{Z}[1/3]\arrow[r]\&\varinjlim(\mathbb{Z}^2,\sigma)\arrow[r]\&\mathbb{Z}\arrow[r]\&0
\end{tikzcd}
\]
still splits. Note that \(\mathbb{Z}\) \emph{is not generated by the eigenvector with eigenvalue \(1\)!}}!} To see this, consider any
\[
\iota:\mathbb{Z}[1/5]\rightarrow\varinjlim(\mathbb{Z}^2,\sigma).
\]
Since \(1\in\mathbb{Z}[1/5]\) is infinitely divisible by \(5\), its image must be as well, thus this map must be of the form \(1\mapsto(c,c)\), an eigenvector with eigenvalue \(5\), with \(c\in\mathbb{Z}\). By writing
\[
\begin{tikzcd}
\mathbb{Z}\arrow[r,"\iota_0"]\arrow[d,"\cdot 5"]&\mathbb{Z}^2\arrow[d,"\sigma"]\\
\mathbb{Z}\arrow[r,"\iota_1"]\arrow[d,"\cdot 5"]&\mathbb{Z}^2\arrow[d,"\sigma"]\\
\vdots\arrow[d]&\vdots\arrow[d]\\
\mathbb{Z}[1/5]\arrow[r,"\iota"]&\varinjlim(\mathbb{Z}^2,\sigma)
\end{tikzcd}
\]
where each \(\iota_n:\mathbb{Z}\rightarrow\mathbb{Z}^2\) is \(1\mapsto(c,c)\),
\[
\coker\iota=\varinjlim\coker\iota_n.
\]
\indent Since \(\mathbb{Z}\leq\coker\iota_n\), without loss of generality, let us assume \(c=1\) and \(\coker\iota_n=\mathbb{Z}\), with \(\mathbb{Z}^2\rightarrow\coker\iota_n\) given by \((a,b)\mapsto a-b\). Choosing a right splitting \(\coker\iota_n\rightarrow\mathbb{Z}^2\) to be \(1\mapsto(1,0)\) gives that the induced map
\[
\sigma_\ast:\coker\iota_n\rightarrow\coker\iota_{n+1}
\]
is multiplication by \(3\), and
\[
\coker\iota=\mathbb{Z}[1/3].
\]
\indent To see that this right splitting does not induce a splitting in the limit, take \(1/3\in\mathbb{Z}[1/3]\). This maps to \((1/3,0)\). By definition, if it belongs to the direct limit, there must be some sufficiently large \(k\in\mathbb{N}\) so that \(\sigma^k(1/3,0)\in\mathbb{Z}^2\), but
\begin{align*}
\sigma^k\left(\begin{array}{c}\frac{1}{3}\\0\end{array}\right)&=\frac{1}{2}\sigma^k\left(\left(\begin{array}{c}\frac{1}{3}\\\frac{1}{3}\end{array}\right)+\left(\begin{array}{c}\frac{1}{3}\\-\frac{1}{3}\end{array}\right)\right)\\
&=\frac{1}{2\cdot 3}\left(5^k\left(\begin{array}{c}1\\1\end{array}\right)+3^k\left(\begin{array}{c}1\\-1\end{array}\right)\right)\notin\mathbb{Z}^2.
\end{align*}
\indent Given any right splitting, since the image of \(1/3^n\) for sufficiently large \(n\in\mathbb{N}\) under the splitting can be decomposed like this to not land in \(\mathbb{Z}^2\) under any \(\sigma^k\), any short exact sequence of the form
\[
\begin{tikzcd}
0\arrow[r]&\mathbb{Z}[1/5]\arrow[r,"\iota"]&\varinjlim(\mathbb{Z}^2,\sigma)\arrow[r]&\mathbb{Z}[1/3]\arrow[r]&0
\end{tikzcd}
\]
cannot be split.\\
\indent Finally, applying the same argument to the map
\[
\mathbb{Z}[1/3]\rightarrow\varinjlim(\mathbb{Z}^2,\sigma)
\]
shows that the direct limit does not split as \(\mathbb{Z}[1/5]\oplus\mathbb{Z}[1/3]\).\\
\indent Therefore, this is an example where the cohomology groups over \(\mathbb{Q}\) are indistinguishable from those of a product, and cup product can be used to differentiate the spaces there, but the cohomology groups over \(\mathbb{Z}\) \emph{also already differ from that of a product!}\\
\indent Lastly, let us compute the frequency module. The substitution matrix on the \(2\)-cells of the dual complex is $\sigma = \sigma_{A3}$ found in Appendix \ref{app:1} which has Perron--Frobenius eigenvector
\[
\resizebox{0.9\hsize}{!}{
\((1913,962,5250,1750,2213,1775,1775,1750,1113,612,912,813,4575,1750,1750,675,17500,1725,675,813,375,1350,474).\)
}
\]
Its sum is \(52500\), giving us that the frequency module is \(\frac{1}{84}\mathbb{Z}[1/5]\), which is different from the frequency module of the product.\\
\indent For illustration, we also compute the frequency modules for substitution on the vertical and horizontal \(1\)-cells of the dual complex. The substitution matrices are
$$\sigma_v^1=\left(\begin{array}{*{5}c}
0&0&3&0&0\\
0&0&1&0&0\\
5&5&0&5&5\\
0&0&1&0&0\\
0&0&0&0&0
\end{array}\right)
\hspace{.4in}\mbox{ and }\hspace{.4in}\sigma_h^1=\left(\begin{array}{*{9}c}
0&0&0&4&3&0&0&0&3\\
2&2&2&0&0&0&0&0&0\\
0&0&0&0&1&1&1&1&1\\
1&1&0&0&0&0&0&0&0\\
1&1&2&0&0&0&0&0&0\\
1&1&1&1&1&0&0&1&1\\
0&0&0&0&0&0&0&0&0\\
0&0&0&0&0&4&4&3&0\\
0&0&0&0&0&0&0&0&0
\end{array}\right).$$
The Perron--Frobenius eigenvectors that are duals of lifts of the generators of \(\check{H}^1(\Omega_\varsigma)\) of eigenvalue \(5\) are $(3,1,5,1,0)$ and $(4,6,5,2,4,7,0,14,0)$. They have sums \(10\) and \(42\), respectively, giving frequency modules \(\frac{1}{2}\mathbb{Z}[1/5]\) and \(\frac{1}{42}\mathbb{Z}[1/5]\), respectively.
\end{example}

\begin{figure}[t]
\centering
\begin{tikzpicture}
\draw(0,0)--(.5,0)--(.5,.5)--(0,.5)--(0,0);
\draw[|->](.75,.25)--(1.25,.25);
\draw(.25,.25)node{\(0\)};
\draw(1.5,-.75)--(3.5,-.75)--(3.5,1.25)--(1.5,1.25)--(1.5,-.75);
\draw(1.5,-.25)--(3.5,-.25);
\draw(1.5,.25)--(3.5,.25);
\draw(1.5,.75)--(3.5,.75);
\draw(2,-.75)--(2,1.25);
\draw(2.5,-.75)--(2.5,1.25);
\draw(3,-.75)--(3,1.25);
\draw(1.75,-.5)node{\(0\)};
\draw(1.75,0)node{\(0\)};
\draw(1.75,.5)node{\(1\)};
\draw(1.75,1)node{\(0\)};
\draw(2.25,-.5)node{\(2\)};
\draw(2.25,0)node{\(2\)};
\draw(2.25,.5)node{\(2\)};
\draw(2.25,1)node{\(2\)};
\draw(2.75,-.5)node{\(1\)};
\draw(2.75,0)node{\(0\)};
\draw(2.75,.5)node{\(1\)};
\draw(2.75,1)node{\(0\)};
\draw(3.25,-.5)node{\(0\)};
\draw(3.25,0)node{\(1\)};
\draw(3.25,.5)node{\(0\)};
\draw(3.25,1)node{\(1\)};

\draw(4.5,0)--(5.0,0)--(5.0,.5)--(4.5,.5)--(4.5,0);
\draw[|->](5.25,.25)--(5.75,.25);
\draw(4.75,.25)node{\(1\)};
\draw(6.0,-.75)--(8.0,-.75)--(8.0,1.25)--(6.0,1.25)--(6.0,-.75);
\draw(6.0,-.25)--(8.0,-.25);
\draw(6.0,.25)--(8.0,.25);
\draw(6.0,.75)--(8.0,.75);
\draw(6.5,-.75)--(6.5,1.25);
\draw(7.0,-.75)--(7.0,1.25);
\draw(7.5,-.75)--(7.5,1.25);
\draw(6.25,-.5)node{\(0\)};
\draw(6.25,0)node{\(0\)};
\draw(6.25,.5)node{\(0\)};
\draw(6.25,1)node{\(1\)};
\draw(6.75,-.5)node{\(2\)};
\draw(6.75,0)node{\(2\)};
\draw(6.75,.5)node{\(2\)};
\draw(6.75,1)node{\(2\)};
\draw(7.25,-.5)node{\(1\)};
\draw(7.25,0)node{\(1\)};
\draw(7.25,.5)node{\(1\)};
\draw(7.25,1)node{\(1\)};
\draw(7.75,-.5)node{\(1\)};
\draw(7.75,0)node{\(0\)};
\draw(7.75,.5)node{\(0\)};
\draw(7.75,1)node{\(1\)};

\draw(9,0)--(9.5,0)--(9.5,.5)--(9,.5)--(9,0);
\draw[|->](9.75,.25)--(10.25,.25);
\draw(9.25,.25)node{\(2\)};
\draw(10.5,-.75)--(12.5,-.75)--(12.5,1.25)--(10.5,1.25)--(10.5,-.75);
\draw(10.5,-.25)--(12.5,-.25);
\draw(10.5,.25)--(12.5,.25);
\draw(10.5,.75)--(12.5,.75);
\draw(11,-.75)--(11,1.25);
\draw(11.5,-.75)--(11.5,1.25);
\draw(12,-.75)--(12,1.25);
\draw(10.75,-.5)node{\(2\)};
\draw(10.75,0)node{\(2\)};
\draw(10.75,.5)node{\(2\)};
\draw(10.75,1)node{\(2\)};
\draw(11.25,-.5)node{\(1\)};
\draw(11.25,0)node{\(0\)};
\draw(11.25,.5)node{\(0\)};
\draw(11.25,1)node{\(1\)};
\draw(11.75,-.5)node{\(2\)};
\draw(11.75,0)node{\(2\)};
\draw(11.75,.5)node{\(2\)};
\draw(11.75,1)node{\(2\)};
\draw(12.25,-.5)node{\(2\)};
\draw(12.25,0)node{\(2\)};
\draw(12.25,.5)node{\(2\)};
\draw(12.25,1)node{\(2\)};
\end{tikzpicture}
\label{figure:non-product4}
\caption{Substitution rule for an expansion \(4\) two-dimensional cubical substitution that has the cohomology groups of a product, but does not have the ring structure of a product.}
\end{figure}%

\begin{example}[Expansion \(4\) Non-product]
  \label{ex:4}

In the final example of this section, we consider the following two-dimensional substitution on three prototiles.
\begin{minted}[linenos,fontsize=\fontsize{5}{6},numbersep=5pt,frame=lines,framesep=2mm,mathescape=true,escapeinside=||]{python}
|\(\varsigma\)|=[[0,0,1,0,2,2,2,2,1,0,1,0,0,1,0,1],[0,0,0,1,2,2,2,2,1,1,1,1,1,0,0,1],[2,2,2,2,1,0,0,1,2,2,2,2,2,2,2,2]]
\end{minted}
This, again, does not force the (see Remark \ref{rem:force}) border, but is primitive and recognizable.\\
\indent Using the dual complex, \(\check{H}^2(\Omega_\varsigma)\) and \(\check{H}^1(\Omega_\varsigma)\) are computed via direct limits under the matrices
\begin{align*}
\sigma^2&=\left(\begin{array}{*{13}c}
1&2&0&0&0&2&1&3&0&0&0&2&0\\
4&4&2&-1&-1&4&-1&-1&6&8&0&0&0\\
4&2&5&0&0&3&0&0&3&4&0&0&1\\
2&2&0&0&0&1&0&3&0&0&1&3&0\\
4&2&4&1&1&2&0&2&0&0&1&3&0\\
3&2&4&1&1&3&1&2&0&0&0&2&0\\
8&4&8&1&2&4&1&5&0&0&2&5&1\\
-2&0&-4&-1&-1&-1&0&0&0&0&0&0&-1\\
2&2&1&0&0&2&-1&-1&3&4&0&0&0\\
1&1&1&0&0&1&0&0&1&8&0&0&0\\
4&2&4&1&1&2&0&2&0&0&1&3&0\\
6&2&8&1&2&3&1&2&0&0&1&2&1\\
2&0&4&1&1&1&0&0&0&0&0&0&1
\end{array}\right)\\
\sigma^1&=\left(\begin{array}{*{4}c}
4&0&0&0\\
1&3&-1&-1\\
1&-1&2&0\\
1&-1&0&2
\end{array}\right),
\end{align*}
with eigenvectors
\begin{align*}
\mathbf{e}_{4_1}^1&=(2,0,1,1) &   \mathbf{e}_{16}^2&=(1,2,2,1,2,2,4,-1,1,1,2,3,1)\\
\mathbf{e}_{4_2}^1&=(0,2,-1,-1) & \mathbf{e}_8^2&=(1,-2,0,1,2,2,4,-1,-1,-1,2,3,1)\\
\mathbf{e}_2^1&=(0,0,1,-1) &  \mathbf{e}_4^2&=(1,-4,-1,1,2,2,4,-1,-2,1,2,3,1)\\
\mathbf{e}_1^1&=(0,1,1,1) & \mathbf{e}_2^2&=(1,0,-1,1,0,0,0,1,0,0,0,-1,-1),
\end{align*}
where the subscripts indicate the corresponding eigenvalues. The cup product is the bilinear form
$$\check{H}^1\times \check{H}^1\rightarrow \check{H}^2:(a,b)\mapsto B(a,b) := \sum_{*\in\{16,8,4,2\}}B_*(a,b)\cdot \mathbf{e}^2_*,$$ where $\mathbf{e}^2_n\in \check{H}^2$ is a generator with eigenvalue $n$, and $B_n$ is the bilinear form given by the matrices
\begin{align*}
B_{16}&:=\left(\begin{array}{*{4}c}
0&2&0&0\\
-2&0&0&0\\
0&0&0&0\\
0&0&0&0
\end{array}\right),
&B_8&:=\left(\begin{array}{*{4}c}
0&0&-2&0\\
0&0&0&0\\
2&0&0&0\\
0&0&0&0
\end{array}\right), \\
B_4&:=\left(\begin{array}{*{4}c}
0&0&0&-2\\
0&0&0&0\\
0&0&0&0\\
2&0&0&0
\end{array}\right),
&B_2&:=\left(\begin{array}{*{4}c}
0&0&0&0\\
0&0&0&0\\
0&0&0&0\\
0&0&0&0
\end{array}\right).
\end{align*}
Like the previous example, the last component of the bilinear form shows that there are no cohomology classes in \(\check{H}^1(\Omega_\varsigma)\) cupping to \(\mathbf{e}_2^2\) in \(\check{H}^2(\Omega_\varsigma)\), contrasting the example coming from the product. This coincides with the fact that, up to a coboundary, the lifts of the eigenvectors \(\mathbf{e}_2^1\) and \(\mathbf{e}_1^1\) to \(1\)-cochains are
\begin{align*}
\widetilde{\mathbf{e}_2^1}&=(0,0,0,0,0,0,0,0,1,1,0,1,1,-1)\\
\widetilde{\mathbf{e}_1^1}&=(0,0,0,0,0,-1,0,0,1,1,-1,1,1,1)
\end{align*}
with the first five entries vertical \(1\)-cells, the remaining horizontal \(1\)-cells, therefore cannot cup nontrivially (to \(\mathbf{e}_2^2\)).\\
\indent Let us now calculate the actual cohomology groups of \(\Omega_\varsigma\) over \(\mathbb{Z}\) following the two-dimensional analogue of the Barge--Diamond filtration. Let us begin by drawing \(S_0\) restricted to the eventual range
\[
\begin{tikzpicture}
\draw(0,0)--(1,0)--(1,1)--(0,1)--(0,0);
\draw(.25,.25)node{\(1\)};
\draw(.25,.75)node{\(0\)};
\draw(.75,.25)node{\(0\)};
\draw(.75,.75)node{\(0\)};
\draw(.5,0)node[below=-0.5ex]{\tiny\(0\)};
\draw(1,.5)node[right=-0.5ex]{\tiny\(0\)};
\draw(.5,1)node[above=-0.5ex]{\tiny\(0\)};
\draw(0,.5)node[left=-0.5ex]{\tiny\(0\)};
\draw(0,0)node[below left=-0.5ex]{\tiny\(0\)};
\draw(1,0)node[below right=-0.5ex]{\tiny\(0\)};
\draw(1,1)node[above right=-0.5ex]{\tiny\(0\)};
\draw(0,1)node[above left=-0.5ex]{\tiny\(0\)};

\draw(2,0)--(3,0)--(3,1)--(2,1)--(2,0);
\draw(2.25,.25)node{\(1\)};
\draw(2.25,.75)node{\(0\)};
\draw(2.75,.25)node{\(2\)};
\draw(2.75,.75)node{\(2\)};
\draw(2.5,0)node[below=-0.5ex]{\tiny\(2\)};
\draw(3,.5)node[right=-0.5ex]{\tiny\(2\)};
\draw(2.5,1)node[above=-0.5ex]{\tiny\(1\)};
\draw(2,.5)node[left=-0.5ex]{\tiny\(0\)};
\draw(2,0)node[below left=-0.5ex]{\tiny\(0\)};
\draw(3,0)node[below right=-0.5ex]{\tiny\(2\)};
\draw(3,1)node[above right=-0.5ex]{\tiny\(1\)};
\draw(2,1)node[above left=-0.5ex]{\tiny\(0\)};

\draw(4,0)--(5,0)--(5,1)--(4,1)--(4,0);
\draw(4.25,.25)node{\(1\)};
\draw(4.25,.75)node{\(1\)};
\draw(4.75,.25)node{\(0\)};
\draw(4.75,.75)node{\(0\)};
\draw(4.5,0)node[below=-0.5ex]{\tiny\(0\)};
\draw(5,.5)node[right=-0.5ex]{\tiny\(0\)};
\draw(4.5,1)node[above=-0.5ex]{\tiny\(2\)};
\draw(4,.5)node[left=-0.5ex]{\tiny\(1\)};
\draw(4,0)node[below left=-0.5ex]{\tiny\(0\)};
\draw(5,0)node[below right=-0.5ex]{\tiny\(0\)};
\draw(5,1)node[above right=-0.5ex]{\tiny\(0\)};
\draw(4,1)node[above left=-0.5ex]{\tiny\(1\)};

\draw(6,0)--(7,0)--(7,1)--(6,1)--(6,0);
\draw(6.25,.25)node{\(1\)};
\draw(6.25,.75)node{\(1\)};
\draw(6.75,.25)node{\(1\)};
\draw(6.75,.75)node{\(0\)};
\draw(6.5,0)node[below=-0.5ex]{\tiny\(1\)};
\draw(7,.5)node[right=-0.5ex]{\tiny\(1\)};
\draw(6.5,1)node[above=-0.5ex]{\tiny\(2\)};
\draw(6,.5)node[left=-0.5ex]{\tiny\(1\)};
\draw(6,0)node[below left=-0.5ex]{\tiny\(0\)};
\draw(7,0)node[below right=-0.5ex]{\tiny\(1\)};
\draw(7,1)node[above right=-0.5ex]{\tiny\(0\)};
\draw(6,1)node[above left=-0.5ex]{\tiny\(1\)};

\draw(8,0)--(9,0)--(9,1)--(8,1)--(8,0);
\draw(8.25,.25)node{\(1\)};
\draw(8.25,.75)node{\(1\)};
\draw(8.75,.25)node{\(2\)};
\draw(8.75,.75)node{\(2\)};
\draw(8.5,0)node[below=-0.5ex]{\tiny\(2\)};
\draw(9,.5)node[right=-0.5ex]{\tiny\(2\)};
\draw(8.5,1)node[above=-0.5ex]{\tiny\(3\)};
\draw(8,.5)node[left=-0.5ex]{\tiny\(1\)};
\draw(8,0)node[below left=-0.5ex]{\tiny\(0\)};
\draw(9,0)node[below right=-0.5ex]{\tiny\(2\)};
\draw(9,1)node[above right=-0.5ex]{\tiny\(1\)};
\draw(8,1)node[above left=-0.5ex]{\tiny\(1\)};

\draw(10,0)--(11,0)--(11,1)--(10,1)--(10,0);
\draw(10.25,.25)node{\(2\)};
\draw(10.25,.75)node{\(2\)};
\draw(10.75,.25)node{\(0\)};
\draw(10.75,.75)node{\(0\)};
\draw(10.5,0)node[below=-0.5ex]{\tiny\(3\)};
\draw(11,.5)node[right=-0.5ex]{\tiny\(0\)};
\draw(10.5,1)node[above=-0.5ex]{\tiny\(4\)};
\draw(10,.5)node[left=-0.5ex]{\tiny\(2\)};
\draw(10,0)node[below left=-0.5ex]{\tiny\(1\)};
\draw(11,0)node[below right=-0.5ex]{\tiny\(0\)};
\draw(11,1)node[above right=-0.5ex]{\tiny\(0\)};
\draw(10,1)node[above left=-0.5ex]{\tiny\(2\)};

\draw(12,0)--(13,0)--(13,1)--(12,1)--(12,0);
\draw(12.25,.25)node{\(2\)};
\draw(12.25,.75)node{\(2\)};
\draw(12.75,.25)node{\(1\)};
\draw(12.75,.75)node{\(0\)};
\draw(12.5,0)node[below=-0.5ex]{\tiny\(4\)};
\draw(13,.5)node[right=-0.5ex]{\tiny\(1\)};
\draw(12.5,1)node[above=-0.5ex]{\tiny\(4\)};
\draw(12,.5)node[left=-0.5ex]{\tiny\(2\)};
\draw(12,0)node[below left=-0.5ex]{\tiny\(1\)};
\draw(13,0)node[below right=-0.5ex]{\tiny\(1\)};
\draw(13,1)node[above right=-0.5ex]{\tiny\(0\)};
\draw(12,1)node[above left=-0.5ex]{\tiny\(2\)};

\draw(14,0)--(15,0)--(15,1)--(14,1)--(14,0);
\draw(14.25,.25)node{\(2\)};
\draw(14.25,.75)node{\(2\)};
\draw(14.75,.25)node{\(2\)};
\draw(14.75,.75)node{\(2\)};
\draw(14.5,0)node[below=-0.5ex]{\tiny\(5\)};
\draw(15,.5)node[right=-0.5ex]{\tiny\(2\)};
\draw(14.5,1)node[above=-0.5ex]{\tiny\(5\)};
\draw(14,.5)node[left=-0.5ex]{\tiny\(2\)};
\draw(14,0)node[below left=-0.5ex]{\tiny\(1\)};
\draw(15,0)node[below right=-0.5ex]{\tiny\(2\)};
\draw(15,1)node[above right=-0.5ex]{\tiny\(1\)};
\draw(14,1)node[above left=-0.5ex]{\tiny\(2\)};
\end{tikzpicture}
\]
where the boundary identifications are \emph{only along the respective positions}, e.g. none of the bottom horizontal \(1\)-cells are identified to any of the top horizontal \(1\)-cells, despite the numbering, and none of the bottom right \(0\)-cells are identified to any of the top right \(0\)-cells, despite the numbering.\\
\indent We arrange the \(1\)-cells and the \(0\)-cells on \(S_0\) in the order
\[
\begin{tikzpicture}
\draw(0,0)--(1,0)--(1,1)--(0,1)--(0,0);
\draw(.5,0)node[below=-0.5ex]{\tiny\(0\)};
\draw(1,.5)node[right=-0.5ex]{\tiny\(1\)};
\draw(.5,1)node[above=-0.5ex]{\tiny\(2\)};
\draw(0,.5)node[left=-0.5ex]{\tiny\(3\)};

\draw(2.0,0)--(3.0,0)--(3.0,1)--(2.0,1)--(2.0,0);
\draw(2,0)node[below left=-0.5ex]{\tiny\(0\)};
\draw(3,0)node[below right=-0.5ex]{\tiny\(1\)};
\draw(3,1)node[above right=-0.5ex]{\tiny\(2\)};
\draw(2,1)node[above left=-0.5ex]{\tiny\(3\)};
\end{tikzpicture}
\]
then by the numbering in \(S_0\). This gives us the coboundary matrices
\begin{small}
\begin{align*}
\partial_0^1&=\left(\begin{array}{*{18}c}
1&0&0&0&0&0&1&0&0&-1&0&0&0&0&0&-1&0&0\\
0&0&1&0&0&0&0&0&1&0&-1&0&0&0&0&-1&0&0\\
1&0&0&0&0&0&1&0&0&0&0&-1&0&0&0&0&-1&0\\
0&1&0&0&0&0&0&1&0&0&0&-1&0&0&0&0&-1&0\\
0&0&1&0&0&0&0&0&1&0&0&0&-1&0&0&0&-1&0\\
0&0&0&1&0&0&1&0&0&0&0&0&0&-1&0&0&0&-1\\
0&0&0&0&1&0&0&1&0&0&0&0&0&-1&0&0&0&-1\\
0&0&0&0&0&1&0&0&1&0&0&0&0&0&-1&0&0&-1
\end{array}\right)\\
\partial_0^0&=\left(\begin{array}{*{10}c}
-1&0&1&0&0&0&0&0&0&0\\
-1&0&0&1&0&0&0&0&0&0\\
-1&0&0&0&1&0&0&0&0&0\\
0&-1&1&0&0&0&0&0&0&0\\
0&-1&0&1&0&0&0&0&0&0\\
0&-1&0&0&1&0&0&0&0&0\\
0&0&-1&0&0&1&0&0&0&0\\
0&0&0&-1&0&1&0&0&0&0\\
0&0&0&0&-1&0&1&0&0&0\\
0&0&0&0&0&1&0&-1&0&0\\
0&0&0&0&0&0&1&-1&0&0\\
0&0&0&0&0&1&0&0&-1&0\\
0&0&0&0&0&0&1&0&-1&0\\
0&0&0&0&0&1&0&0&0&-1\\
0&0&0&0&0&0&1&0&0&-1\\
-1&0&0&0&0&0&0&1&0&0\\
-1&0&0&0&0&0&0&0&1&0\\
0&-1&0&0&0&0&0&0&0&1
\end{array}\right).
\end{align*}
\end{small}
Direct calculation shows that
\begin{align*}
\check{H}^2(\Xi_0)&=0\\
\check{H}^1(\Xi_0)&=\mathbb{Z}=\langle\overline{(1,0,0,0,-1,0,0,1,0,0,-1,0,-1,0,0,1,1,0)}\rangle\\
\check{H}^0(\Xi_0)&=\mathbb{Z}
\end{align*}
where the vector is a \(1\)-cochain in \(S_0\).\\
\indent The extra cells that we add on to \(S_0\) to form \(S_1\) are the same as the previous example, being the vertical \(1\)-flaps
\[
\begin{tikzpicture}
\draw(0,0)--(1,0)--(1,2)--(0,2)--(0,0);
\draw(.25,1)node{\(0\)};
\draw(.75,1)node{\(0\)};
\draw(0,1)node[left=-0.5ex]{\tiny\(0\)};
\draw(1,1)node[right=-0.5ex]{\tiny\(0\)};

\draw(2,0)--(3,0)--(3,2)--(2,2)--(2,0);
\draw(2.25,1)node{\(0\)};
\draw(2.75,1)node{\(1\)};
\draw(2,1)node[left=-0.5ex]{\tiny\(0\)};
\draw(3,1)node[right=-0.5ex]{\tiny\(1\)};

\draw(4,0)--(5,0)--(5,2)--(4,2)--(4,0);
\draw(4.25,1)node{\(0\)};
\draw(4.75,1)node{\(2\)};
\draw(4,1)node[left=-0.5ex]{\tiny\(0\)};
\draw(5,1)node[right=-0.5ex]{\tiny\(2\)};

\draw(6,0)--(7,0)--(7,2)--(6,2)--(6,0);
\draw(6.25,1)node{\(1\)};
\draw(6.75,1)node{\(0\)};
\draw(6,1)node[left=-0.5ex]{\tiny\(1\)};
\draw(7,1)node[right=-0.5ex]{\tiny\(0\)};

\draw(8,0)--(9,0)--(9,2)--(8,2)--(8,0);
\draw(8.25,1)node{\(1\)};
\draw(8.75,1)node{\(1\)};
\draw(8,1)node[left=-0.5ex]{\tiny\(1\)};
\draw(9,1)node[right=-0.5ex]{\tiny\(1\)};

\draw(10,0)--(11,0)--(11,2)--(10,2)--(10,0);
\draw(10.25,1)node{\(1\)};
\draw(10.75,1)node{\(2\)};
\draw(10,1)node[left=-0.5ex]{\tiny\(1\)};
\draw(11,1)node[right=-0.5ex]{\tiny\(2\)};

\draw(0,-3)--(1,-3)--(1,-1)--(0,-1)--(0,-3);
\draw(.25,-2)node{\(2\)};
\draw(.75,-2)node{\(0\)};
\draw(0,-2)node[left=-0.5ex]{\tiny\(2\)};
\draw(1,-2)node[right=-0.5ex]{\tiny\(0\)};

\draw(2,-3)--(3,-3)--(3,-1)--(2,-1)--(2,-3);
\draw(2.25,-2)node{\(2\)};
\draw(2.75,-2)node{\(1\)};
\draw(2,-2)node[left=-0.5ex]{\tiny\(2\)};
\draw(3,-2)node[right=-0.5ex]{\tiny\(1\)};

\draw(4,-3)--(5,-3)--(5,-1)--(4,-1)--(4,-3);
\draw(4.25,-2)node{\(2\)};
\draw(4.75,-2)node{\(2\)};
\draw(4,-2)node[left=-0.5ex]{\tiny\(2\)};
\draw(5,-2)node[right=-0.5ex]{\tiny\(2\)};
\end{tikzpicture}
\]
where the top and the bottom horizontal \(1\)-cells are identified to their appropriate bottom and top horizontal \(1\)-cells in \(S_0\), and the horizontal \(1\)-flaps
\[
\begin{tikzpicture}
\draw(0,0)--(2,0)--(2,1)--(0,1)--(0,0);
\draw(1,.25)node{\(0\)};
\draw(1,.75)node{\(0\)};
\draw(1,0)node[below=-0.5ex]{\tiny\(0\)};
\draw(1,1)node[above=-0.5ex]{\tiny\(0\)};

\draw(3,0)--(5,0)--(5,1)--(3,1)--(3,0);
\draw(4,.25)node{\(0\)};
\draw(4,.75)node{\(1\)};
\draw(4,0)node[below=-0.5ex]{\tiny\(0\)};
\draw(4,1)node[above=-0.5ex]{\tiny\(1\)};

\draw(6,0)--(8,0)--(8,1)--(6,1)--(6,0);
\draw(7,.25)node{\(1\)};
\draw(7,.75)node{\(0\)};
\draw(7,0)node[below=-0.5ex]{\tiny\(1\)};
\draw(7,1)node[above=-0.5ex]{\tiny\(0\)};

\draw(9,0)--(11,0)--(11,1)--(9,1)--(9,0);
\draw(10,.25)node{\(1\)};
\draw(10,.75)node{\(1\)};
\draw(10,0)node[below=-0.5ex]{\tiny\(1\)};
\draw(10,1)node[above=-0.5ex]{\tiny\(1\)};

\draw(12,0)--(14,0)--(14,1)--(12,1)--(12,0);
\draw(13,.25)node{\(2\)};
\draw(13,.75)node{\(2\)};
\draw(13,0)node[below=-0.5ex]{\tiny\(2\)};
\draw(13,1)node[above=-0.5ex]{\tiny\(2\)};
\end{tikzpicture}
\]
where the left and the right vertical \(1\)-cells are identified to their appropriate right and left vertical \(1\)-cells in \(S_0\). \(S_1/S_0\) has three generating loops up to homotopy, one from the vertical \(1\)-flaps, one from the first four horizontal \(1\)-flaps, and one from the last horizontal \(1\)-flap. We then compute that
\[
\check{H}^1(\Xi_1,\Xi_0)=\varinjlim\left(\mathbb{Z}^3,\left(\begin{array}{*{3}c}
4&0&0\\
0&3&1\\
0&1&3
\end{array}\right)\right)
\]
with the order of the entries of the matrix given as described. The matrix has eigenvectors
\begin{align*}
\mathbf{e}_{4_1}&=(1,0,0)\\
\mathbf{e}_{4_2}&=(0,1,1)\\
\mathbf{e}_2&=(0,1,-1)
\end{align*}
with the subscripts their respective eigenvalues. Direct computation then shows that the direct limit splits as a direct sum, giving us
\[
\check{H}^1(\Xi_1,\Xi_0)=\mathbb{Z}[1/4]^2\oplus\mathbb{Z}[1/2].
\]
\indent \(\check{H}^2(\Xi_1,\Xi_0)\) requires more effort. We observe that the generating \(2\)-spheres of \(S_1/S_0\), up to homotopy, occur when \(1\)-flaps form a ``tube'', e.g. the first four of the horizontal \(1\)-flaps. Let us represent all such tubes by graphs that represents the horizontal and vertical cross sections of vertical and horizontal \(1\)-flaps, respectively, with the vertices the \(1\)-cells in \(S_1\) not attached to \(S_0\), and the edges the \(1\)-flaps themselves. In this example, they are, drawn in the orientation matching the cross sections of the respective \(1\)-flaps,
\[
\begin{tikzpicture}
\draw(.5,0)--(0,-.5);
\draw(.5,-.5)--(0,-1);
\draw(.5,0)--(0,-1);
\draw(.5,-.5)--(0,0);
\draw(0,-.5)--(.5,-1);
\draw(0,0)--(.5,-1);
\draw(0,0)--(.5,0);
\draw(0,-.5)--(.5,-.5);
\draw(0,-1)--(.5,-1);
\draw(0,0)node[left=-0.5ex]{\tiny\(0\)};
\draw(0,-.5)node[left=-0.5ex]{\tiny\(1\)};
\draw(0,-1)node[left=-0.5ex]{\tiny\(2\)};
\draw(.5,0)node[right=-0.5ex]{\tiny\(0\)};
\draw(.5,-.5)node[right=-0.5ex]{\tiny\(1\)};
\draw(.5,-1)node[right=-0.5ex]{\tiny\(2\)};

\draw(2.0,-.25)--(1.5,-.75);
\draw(2.0,-.75)--(1.5,-.25);
\draw(1.5,-.25)--(1.5,-.75);
\draw(2.0,-.25)--(2.0,-.75);
\draw(1.5,-.25)node[above=-0.5ex]{\tiny\(0\)};
\draw(1.5,-.75)node[below=-0.5ex]{\tiny\(0\)};
\draw(2.0,-.25)node[above=-0.5ex]{\tiny\(1\)};
\draw(2.0,-.75)node[below=-0.5ex]{\tiny\(1\)};

\draw(3.0,-.25)--(3.0,-.75);
\draw(3.0,-.25)node[above=-0.5ex]{\tiny\(2\)};
\draw(3.0,-.75)node[below=-0.5ex]{\tiny\(2\)};
\end{tikzpicture}
\]
with labels of the vertices corresponding to the labels of the \(1\)-cells marked in the set of \(1\)-flaps.\\
\indent The first graph corresponding to vertical \(1\)-flaps has a cycle basis consisting of four elements, therefore four associated generating \(2\)-spheres in \(S_1/S_0\), which we pick to be
\[
\begin{tikzpicture}[on top/.style={preaction={draw=white,-,line width=#1}},on top/.default=4pt]
\draw(.5,0)--(0,-.5);
\draw[on top](.5,-.5)--(0,0);
\draw(0,0)--(.5,0);
\draw(0,-.5)--(.5,-.5);
\draw(0,0)node[left=-0.5ex]{\tiny\(0\)};
\draw(0,-.5)node[left=-0.5ex]{\tiny\(1\)};
\draw(.5,0)node[right=-0.5ex]{\tiny\(0\)};
\draw(.5,-.5)node[right=-0.5ex]{\tiny\(1\)};

\draw(2.0,-.5)--(1.5,-1);
\draw[on top](1.5,-.5)--(2.0,-1);
\draw(1.5,-.5)--(2.0,-.5);
\draw(1.5,-1)--(2.0,-1);
\draw(1.5,-.5)node[left=-0.5ex]{\tiny\(1\)};
\draw(1.5,-1)node[left=-0.5ex]{\tiny\(2\)};
\draw(2.0,-.5)node[right=-0.5ex]{\tiny\(1\)};
\draw(2.0,-1)node[right=-0.5ex]{\tiny\(2\)};

\draw(3.5,0)--(3,-1);
\draw[on top](3,0)--(3.5,-1);
\draw(3,0)--(3.5,0);
\draw(3,-1)--(3.5,-1);
\draw(3,0)node[left=-0.5ex]{\tiny\(0\)};
\draw(3,-1)node[left=-0.5ex]{\tiny\(2\)};
\draw(3.5,0)node[right=-0.5ex]{\tiny\(0\)};
\draw(3.5,-1)node[right=-0.5ex]{\tiny\(2\)};

\draw(4.5,-.5)--(5.0,-.5);
\draw[on top](4.5,0)--(5.0,-1);
\draw(5.0,-.5)--(4.5,0);
\draw(4.5,-.5)--(5.0,-1);
\draw(4.5,0)node[left=-0.5ex]{\tiny\(0\)};
\draw(4.5,-.5)node[left=-0.5ex]{\tiny\(1\)};
\draw(5.0,-.5)node[right=-0.5ex]{\tiny\(1\)};
\draw(5.0,-1)node[right=-0.5ex]{\tiny\(2\)};
\end{tikzpicture}
\]
oriented so that the top edge is left-to-right. In fact, let us orient for all other cycles the same way. The substitution rule then dictates that
\[
\begin{tikzpicture}[on top/.style={preaction={draw=white,-,line width=#1}},on top/.default=4pt]
\draw(.5,0)--(0,-.5);
\draw[on top](.5,-.5)--(0,0);
\draw(0,0)--(.5,0);
\draw(0,-.5)--(.5,-.5);
\draw(0,0)node[left=-0.5ex]{\tiny\(0\)};
\draw(0,-.5)node[left=-0.5ex]{\tiny\(1\)};
\draw(.5,0)node[right=-0.5ex]{\tiny\(0\)};
\draw(.5,-.5)node[right=-0.5ex]{\tiny\(1\)};
\draw[|->](1,-.5)--(1.5,-.5);
\draw(2,-.5)--(2.5,0);
\draw(2,-.5)--(2.5,-.5);
\draw(2,-.5)node[left=-0.5ex]{\tiny\(1\)};
\draw(2.5,0)node[right=-0.5ex]{\tiny\(0\)};
\draw(2.5,-.5)node[right=-0.5ex]{\tiny\(1\)};
\draw(3,-.5)node{\(+\)};
\draw(3.5,0)--(4,0);
\draw(3.5,0)--(4,-.5);
\draw(3.5,0)node[left=-0.5ex]{\tiny\(0\)};
\draw(4,0)node[right=-0.5ex]{\tiny\(0\)};
\draw(4,-.5)node[right=-0.5ex]{\tiny\(1\)};
\draw(4.5,-.5)node{\(+\)};
\draw(5,0)--(5.5,0);
\draw(5,-.5)--(5.5,0);
\draw(5,0)node[left=-0.5ex]{\tiny\(0\)};
\draw(5,-.5)node[left=-0.5ex]{\tiny\(1\)};
\draw(5.5,0)node[right=-0.5ex]{\tiny\(0\)};
\draw(6,-.5)node{\(+\)};
\draw(6.5,0)--(7,0);
\draw(6.5,-.5)--(7,0);
\draw(6.5,0)node[left=-0.5ex]{\tiny\(0\)};
\draw(6.5,-.5)node[left=-0.5ex]{\tiny\(1\)};
\draw(7,0)node[right=-0.5ex]{\tiny\(0\)};

\draw(9.0,-.5)--(8.5,-1);
\draw[on top](9.0,-1)--(8.5,-.5);
\draw(8.5,-.5)--(9.0,-.5);
\draw(8.5,-1)--(9.0,-1);
\draw(8.5,-.5)node[left=-0.5ex]{\tiny\(1\)};
\draw(8.5,-1)node[left=-0.5ex]{\tiny\(2\)};
\draw(9.0,-.5)node[right=-0.5ex]{\tiny\(1\)};
\draw(9.0,-1)node[right=-0.5ex]{\tiny\(2\)};
\draw[|->](9.5,-.5)--(10.0,-.5);
\draw(11.0,-.5)--(10.5,-1);
\draw[on top](11.0,-1)--(10.5,-.5);
\draw(10.5,-.5)--(11.0,-.5);
\draw(10.5,-1)--(11.0,-1);
\draw(10.5,-.5)node[left=-0.5ex]{\tiny\(1\)};
\draw(10.5,-1)node[left=-0.5ex]{\tiny\(2\)};
\draw(11.0,-.5)node[right=-0.5ex]{\tiny\(1\)};
\draw(11.0,-1)node[right=-0.5ex]{\tiny\(2\)};
\draw(11.5,-.5)node{\(+\)};
\draw(12.0,-1)--(12.5,0);
\draw[on top](12.0,0)--(12.5,-1);
\draw(12.0,0)--(12.5,0);
\draw(12.0,-1)--(12.5,-1);
\draw(12.0,0)node[left=-0.5ex]{\tiny\(0\)};
\draw(12.0,-1)node[left=-0.5ex]{\tiny\(2\)};
\draw(12.5,0)node[right=-0.5ex]{\tiny\(0\)};
\draw(12.5,-1)node[right=-0.5ex]{\tiny\(2\)};
\draw(13.0,-.5)node{\(+\)};
\draw(13.5,-1)--(14.0,0);
\draw[on top](13.5,0)--(14.0,-1);
\draw(13.5,0)--(14.0,0);
\draw(13.5,-1)--(14.0,-1);
\draw(13.5,0)node[left=-0.5ex]{\tiny\(0\)};
\draw(13.5,-1)node[left=-0.5ex]{\tiny\(2\)};
\draw(14.0,0)node[right=-0.5ex]{\tiny\(0\)};
\draw(14.0,-1)node[right=-0.5ex]{\tiny\(2\)};
\draw(14.5,-.5)node{\(+\)};
\draw(15.0,-1)--(15.5,0);
\draw[on top](15.0,-.5)--(15.5,-1);
\draw(15.0,-.5)--(15.5,0);
\draw(15.0,-1)--(15.5,-1);
\draw(15.0,-.5)node[left=-0.5ex]{\tiny\(1\)};
\draw(15.0,-1)node[left=-0.5ex]{\tiny\(2\)};
\draw(15.5,0)node[right=-0.5ex]{\tiny\(0\)};
\draw(15.5,-1)node[right=-0.5ex]{\tiny\(2\)};

\draw(.5,-2)--(0,-3);
\draw[on top](.5,-3)--(0,-2);
\draw(0,-2)--(.5,-2);
\draw(0,-3)--(.5,-3);
\draw(0,-2)node[left=-0.5ex]{\tiny\(0\)};
\draw(0,-3)node[left=-0.5ex]{\tiny\(2\)};
\draw(.5,-2)node[right=-0.5ex]{\tiny\(0\)};
\draw(.5,-3)node[right=-0.5ex]{\tiny\(2\)};
\draw[|->](1.0,-2.5)--(1.5,-2.5);
\draw(2.0,-3)--(2.5,-2);
\draw[on top](2.0,-2.5)--(2.5,-3);
\draw(2.0,-2.5)--(2.5,-2);
\draw(2.0,-3)--(2.5,-3);
\draw(2.0,-2.5)node[left=-0.5ex]{\tiny\(1\)};
\draw(2.0,-3)node[left=-0.5ex]{\tiny\(2\)};
\draw(2.5,-2)node[right=-0.5ex]{\tiny\(0\)};
\draw(2.5,-3)node[right=-0.5ex]{\tiny\(2\)};
\draw(3.0,-2.5)node{\(+\)};
\draw(3.5,-3)--(4.0,-2.5);
\draw[on top](3.5,-2)--(4.0,-3);
\draw(3.5,-2)--(4.0,-2.5);
\draw(3.5,-3)--(4.0,-3);
\draw(3.5,-2)node[left=-0.5ex]{\tiny\(0\)};
\draw(3.5,-3)node[left=-0.5ex]{\tiny\(2\)};
\draw(4.0,-2.5)node[right=-0.5ex]{\tiny\(1\)};
\draw(4.0,-3)node[right=-0.5ex]{\tiny\(2\)};
\draw(4.5,-2.5)node{\(+\)};
\draw(5.0,-3)--(5.5,-2);
\draw[on top](5.0,-2.5)--(5.5,-3);
\draw(5.0,-2.5)--(5.5,-2);
\draw(5.0,-3)--(5.5,-3);
\draw(5.0,-2.5)node[left=-0.5ex]{\tiny\(1\)};
\draw(5.0,-3)node[left=-0.5ex]{\tiny\(2\)};
\draw(5.5,-2)node[right=-0.5ex]{\tiny\(0\)};
\draw(5.5,-3)node[right=-0.5ex]{\tiny\(2\)};
\draw(6.0,-2.5)node{\(+\)};
\draw(6.5,-3)--(7.0,-2);
\draw[on top](6.5,-2)--(7.0,-3);
\draw(6.5,-2)--(7.0,-2);
\draw(6.5,-3)--(7.0,-3);
\draw(6.5,-2)node[left=-0.5ex]{\tiny\(0\)};
\draw(6.5,-3)node[left=-0.5ex]{\tiny\(2\)};
\draw(7.0,-2)node[right=-0.5ex]{\tiny\(0\)};
\draw(7.0,-3)node[right=-0.5ex]{\tiny\(2\)};

\draw(9.0,-2.5)--(8.5,-2.5);
\draw[on top](9.0,-3)--(8.5,-2);
\draw(8.5,-2)--(9.0,-2.5);
\draw(8.5,-2.5)--(9.0,-3);
\draw(8.5,-2)node[left=-0.5ex]{\tiny\(0\)};
\draw(8.5,-2.5)node[left=-0.5ex]{\tiny\(1\)};
\draw(9.0,-2.5)node[right=-0.5ex]{\tiny\(1\)};
\draw(9.0,-3)node[right=-0.5ex]{\tiny\(2\)};
\draw[|->](9.5,-2.5)--(10.0,-2.5);
\draw(11.0,-3)--(10.5,-2.5);
\draw(10.5,-2.5)--(11.0,-2.5);
\draw(10.5,-2.5)node[left=-0.5ex]{\tiny\(1\)};
\draw(11.0,-2.5)node[right=-0.5ex]{\tiny\(1\)};
\draw(11.0,-3)node[right=-0.5ex]{\tiny\(2\)};
\draw(11.5,-2.5)node{\(+\)};
\draw(12.0,-2)--(12.5,-2);
\draw(12.0,-2)--(12.5,-3);
\draw(12.0,-2)node[left=-0.5ex]{\tiny\(0\)};
\draw(12.5,-2)node[right=-0.5ex]{\tiny\(0\)};
\draw(12.5,-3)node[right=-0.5ex]{\tiny\(2\)};
\draw(13.0,-2.5)node{\(-\)};
\draw(13.5,-2.5)--(14.0,-2);
\draw[on top](13.5,-2)--(14.0,-3);
\draw(13.5,-2)--(14.0,-2);
\draw(13.5,-2.5)--(14.0,-3);
\draw(13.5,-2)node[left=-0.5ex]{\tiny\(0\)};
\draw(13.5,-2.5)node[left=-0.5ex]{\tiny\(1\)};
\draw(14.0,-2)node[right=-0.5ex]{\tiny\(0\)};
\draw(14.0,-3)node[right=-0.5ex]{\tiny\(2\)};
\draw(14.5,-2.5)node{\(+\)};
\draw(15.0,-2.5)--(15.5,-2);
\draw[on top](15.0,-2)--(15.5,-3);
\draw(15.0,-2)--(15.5,-2);
\draw(15.0,-2.5)--(15.5,-3);
\draw(15.0,-2)node[left=-0.5ex]{\tiny\(0\)};
\draw(15.0,-2.5)node[left=-0.5ex]{\tiny\(1\)};
\draw(15.5,-2)node[right=-0.5ex]{\tiny\(0\)};
\draw(15.5,-3)node[right=-0.5ex]{\tiny\(2\)};
\end{tikzpicture}.
\]
The cycles not in our basis decompose as
\[
\begin{tikzpicture}[on top/.style={preaction={draw=white,-,line width=#1}},on top/.default=4pt]
\draw(0,-1.0)--(.5,0);
\draw[on top](0,-.5)--(.5,-1.0);
\draw(0,-.5)--(.5,0);
\draw(0,-1.0)--(.5,-1.0);
\draw(0,-.5)node[left=-0.5ex]{\tiny\(1\)};
\draw(0,-1.0)node[left=-0.5ex]{\tiny\(2\)};
\draw(.5,0)node[right=-0.5ex]{\tiny\(0\)};
\draw(.5,-1.0)node[right=-0.5ex]{\tiny\(2\)};
\draw(1,-.5)node{\(=\)};
\draw(2.0,0)--(1.5,-1);
\draw[on top](1.5,0)--(2.0,-1);
\draw(1.5,0)--(2.0,0);
\draw(1.5,-1)--(2.0,-1);
\draw(1.5,0)node[left=-0.5ex]{\tiny\(0\)};
\draw(1.5,-1)node[left=-0.5ex]{\tiny\(2\)};
\draw(2.0,0)node[right=-0.5ex]{\tiny\(0\)};
\draw(2.0,-1)node[right=-0.5ex]{\tiny\(2\)};
\draw(2.5,-.5)node{\(-\)};
\draw(3.5,0)--(3,-.5);
\draw[on top](3.5,-.5)--(3,0);
\draw(3,0)--(3.5,0);
\draw(3,-.5)--(3.5,-.5);
\draw(3,0)node[left=-0.5ex]{\tiny\(0\)};
\draw(3,-.5)node[left=-0.5ex]{\tiny\(1\)};
\draw(3.5,0)node[right=-0.5ex]{\tiny\(0\)};
\draw(3.5,-.5)node[right=-0.5ex]{\tiny\(1\)};
\draw(4,-.5)node{\(-\)};
\draw(4.5,-.5)--(5.0,-.5);
\draw[on top](4.5,0)--(5.0,-1);
\draw(5.0,-.5)--(4.5,0);
\draw(4.5,-.5)--(5.0,-1);
\draw(4.5,0)node[left=-0.5ex]{\tiny\(0\)};
\draw(4.5,-.5)node[left=-0.5ex]{\tiny\(1\)};
\draw(5.0,-.5)node[right=-0.5ex]{\tiny\(1\)};
\draw(5.0,-1)node[right=-0.5ex]{\tiny\(2\)};

\draw(0,-3.0)--(.5,-2.5);
\draw[on top](0,-2.0)--(.5,-3.0);
\draw(0,-2.0)--(.5,-2.5);
\draw(0,-3.0)--(.5,-3.0);
\draw(0,-2.0)node[left=-0.5ex]{\tiny\(0\)};
\draw(0,-3.0)node[left=-0.5ex]{\tiny\(2\)};
\draw(.5,-2.5)node[right=-0.5ex]{\tiny\(1\)};
\draw(.5,-3.0)node[right=-0.5ex]{\tiny\(2\)};
\draw(1,-2.5)node{\(=\)};
\draw(2.0,-2.5)--(1.5,-3);
\draw[on top](1.5,-2.5)--(2.0,-3);
\draw(1.5,-2.5)--(2.0,-2.5);
\draw(1.5,-3)--(2.0,-3);
\draw(1.5,-2.5)node[left=-0.5ex]{\tiny\(1\)};
\draw(1.5,-3)node[left=-0.5ex]{\tiny\(2\)};
\draw(2.0,-2.5)node[right=-0.5ex]{\tiny\(1\)};
\draw(2.0,-3)node[right=-0.5ex]{\tiny\(2\)};
\draw(2.5,-2.5)node{\(+\)};
\draw(3,-2.5)--(3.5,-2.5);
\draw[on top](3,-2)--(3.5,-3);
\draw(3.5,-2.5)--(3,-2);
\draw(3,-2.5)--(3.5,-3);
\draw(3,-2)node[left=-0.5ex]{\tiny\(0\)};
\draw(3,-2.5)node[left=-0.5ex]{\tiny\(1\)};
\draw(3.5,-2.5)node[right=-0.5ex]{\tiny\(1\)};
\draw(3.5,-3)node[right=-0.5ex]{\tiny\(2\)};
\end{tikzpicture}.
\]
A similar calculation for the second and the third graphs gives us
\[
\begin{tikzpicture}[on top/.style={preaction={draw=white,-,line width=#1}},on top/.default=4pt]
\draw(.5,-.75)--(0,-.25);
\draw[on top](0,-.75)--(.5,-.25);
\draw(0,-.75)--(0,-.25);
\draw(.5,-.75)--(.5,-.25);
\draw(0,-.25)node[above=-0.5ex]{\tiny\(0\)};
\draw(0,-.75)node[below=-0.5ex]{\tiny\(0\)};
\draw(.5,-.25)node[above=-0.5ex]{\tiny\(1\)};
\draw(.5,-.75)node[below=-0.5ex]{\tiny\(1\)};
\draw[|->](1.0,-.5)--(1.5,-.5);
\draw(2.5,-.75)--(2.0,-.25);
\draw(2.0,-.75)--(2.0,-.25);
\draw(2.0,-.25)node[above=-0.5ex]{\tiny\(0\)};
\draw(2.0,-.75)node[below=-0.5ex]{\tiny\(0\)};
\draw(2.5,-.75)node[below=-0.5ex]{\tiny\(1\)};
\draw(3.0,-.5)node{\(+\)};
\draw(3.5,-.75)--(3.5,-.25);
\draw(3.5,-.25)node[above=-0.5ex]{\tiny\(2\)};
\draw(3.5,-.75)node[below=-0.5ex]{\tiny\(2\)};
\draw(4.0,-.5)node{\(+\)};
\draw(4.5,-.75)--(5.0,-.25);
\draw(5.0,-.75)--(5.0,-.25);
\draw(4.5,-.75)node[below=-0.5ex]{\tiny\(0\)};
\draw(5.0,-.25)node[above=-0.5ex]{\tiny\(1\)};
\draw(5.0,-.75)node[below=-0.5ex]{\tiny\(1\)};
\draw(5.5,-.5)node{\(+\)};
\draw(6.5,-.75)--(6.0,-.25);
\draw(6.5,-.75)--(6.5,-.25);
\draw(6.0,-.25)node[above=-0.5ex]{\tiny\(0\)};
\draw(6.5,-.25)node[above=-0.5ex]{\tiny\(1\)};
\draw(6.5,-.75)node[below=-0.5ex]{\tiny\(1\)};

\draw(.5,-2.25)--(.5,-1.75);
\draw(.5,-1.75)node[above=-0.5ex]{\tiny\(2\)};
\draw(.5,-2.25)node[below=-0.5ex]{\tiny\(2\)};
\draw[|->](1.0,-2.0)--(1.5,-2.0);
\draw(2.0,-2.25)--(2.0,-1.75);
\draw(2.0,-1.75)node[above=-0.5ex]{\tiny\(2\)};
\draw(2.0,-2.25)node[below=-0.5ex]{\tiny\(2\)};
\draw(2.5,-2.0)node{\(+\)};
\draw(3.0,-2.25)--(3.0,-1.75);
\draw(3.0,-1.75)node[above=-0.5ex]{\tiny\(1\)};
\draw(3.0,-2.25)node[below=-0.5ex]{\tiny\(1\)};
\draw(3.5,-2.0)node{\(+\)};
\draw(4.0,-2.25)--(4.0,-1.75);
\draw(4.0,-1.75)node[above=-0.5ex]{\tiny\(2\)};
\draw(4.0,-2.25)node[below=-0.5ex]{\tiny\(2\)};
\draw(4.5,-2.0)node{\(+\)};
\draw(5.0,-2.25)--(5.0,-1.75);
\draw(5.0,-1.75)node[above=-0.5ex]{\tiny\(2\)};
\draw(5.0,-2.25)node[below=-0.5ex]{\tiny\(2\)};
\end{tikzpicture}.
\]
\indent Observing that the substitution rule is trivial on the second and the third graphs indicates that we can calculate \(\check{H}^2(\Xi_1,\Xi_0)\) using only the vertical \(1\)-flaps, which is
\[
\check{H}^2(\Xi_1,\Xi_0)=\varinjlim\left(\mathbb{Z}^4,\left(\begin{array}{*{4}c}
0&0&0&0\\
-1&1&3&-1\\
-2&1&3&-1\\
0&0&0&0
\end{array}\right)\right)
\]
where the entries are ordered in the cycle basis we selected. Finally,
\[
\check{H}^2(\Xi_1,\Xi_0)=\mathbb{Z}[1/4]=\langle(0,1,1,0)\rangle.
\]
Thus the long exact sequence in relative cohomology of a pair \((\Xi_1,\Xi_0)\) reads
\[
\begin{tikzcd}
0\arrow[r]&\check{H}^0(\Xi_1,\Xi_0)\arrow[r]&\check{H}^0(\Xi_1)\arrow[r]\arrow[d,phantom,""{coordinate,name=X}]&\mathbb{Z}\arrow[lld,rounded corners,to path={--([xshift=2ex]\tikztostart.east)|-(X.center)\tikztonodes-|([xshift=-2ex]\tikztotarget.west)--(\tikztotarget)},swap,"\delta_0^0" at end]&\\
&\mathbb{Z}[1/4]^2\oplus\mathbb{Z}[1/2]\arrow[r]&\check{H}^1(\Xi_1)\arrow[r]\arrow[d,phantom,""{coordinate,name=Y}]&\mathbb{Z}\arrow[lld,rounded corners,to path={--([xshift=2ex]\tikztostart.east)|-(Y.center)\tikztonodes-|([xshift=-2ex]\tikztotarget.west)--(\tikztotarget)},swap,"\delta_0^1" at end]&\\
&\mathbb{Z}[1/4]\arrow[r]&\check{H}^2(\Xi_1)\arrow[r]&0\arrow[r]&0
\end{tikzcd}.
\]
\indent Recalling that
\[
\check{H}^1(\Xi_0)=\mathbb{Z}=\langle\overline{(1,0,0,0,-1,0,0,1,0,0,-1,0,-1,0,0,1,1,0)}\rangle
\]
we check the image of the generator under \(\delta_0^1\), which is
\[
\begin{tikzpicture}
\draw(0,0)--(1,0)--(1,1)--(0,1)--(0,0);
\draw(.25,.25)node{\(1\)};
\draw(.25,.75)node{\(\)};
\draw(.75,.25)node{\(0\)};
\draw(.75,.75)node{\(\)};
\draw(.5,0)node[below=-0.5ex]{\tiny\(0\)};
\draw(1,.5)node[right=-0.5ex]{\tiny\(\)};
\draw(.5,1)node[above=-0.5ex]{\tiny\(\)};
\draw(0,.5)node[left=-0.5ex]{\tiny\(\)};
\draw(0,0)node[below left=-0.5ex]{\tiny\(\)};
\draw(1,0)node[below right=-0.5ex]{\tiny\(\)};
\draw(1,1)node[above right=-0.5ex]{\tiny\(\)};
\draw(0,1)node[above left=-0.5ex]{\tiny\(\)};

\draw(1.5,.5)node{\(-\)};

\draw(2,0)--(3,0)--(3,1)--(2,1)--(2,0);
\draw(2.25,.25)node{\(2\)};
\draw(2.25,.75)node{\(\)};
\draw(2.75,.25)node{\(1\)};
\draw(2.75,.75)node{\(\)};
\draw(2.5,0)node[below=-0.5ex]{\tiny\(4\)};
\draw(3,.5)node[right=-0.5ex]{\tiny\(\)};
\draw(2.5,1)node[above=-0.5ex]{\tiny\(\)};
\draw(2,.5)node[left=-0.5ex]{\tiny\(\)};
\draw(2,0)node[below left=-0.5ex]{\tiny\(\)};
\draw(3,0)node[below right=-0.5ex]{\tiny\(\)};
\draw(3,1)node[above right=-0.5ex]{\tiny\(\)};
\draw(2,1)node[above left=-0.5ex]{\tiny\(\)};

\draw(3.5,.5)node{\(+\)};

\draw(4,0)--(5,0)--(5,1)--(4,1)--(4,0);
\draw(4.25,.25)node{\(\)};
\draw(4.25,.75)node{\(\)};
\draw(4.75,.25)node{\(1\)};
\draw(4.75,.75)node{\(0\)};
\draw(4.5,0)node[below=-0.5ex]{\tiny\(\)};
\draw(5,.5)node[right=-0.5ex]{\tiny\(1\)};
\draw(4.5,1)node[above=-0.5ex]{\tiny\(\)};
\draw(4,.5)node[left=-0.5ex]{\tiny\(\)};
\draw(4,0)node[below left=-0.5ex]{\tiny\(\)};
\draw(5,0)node[below right=-0.5ex]{\tiny\(\)};
\draw(5,1)node[above right=-0.5ex]{\tiny\(\)};
\draw(4,1)node[above left=-0.5ex]{\tiny\(\)};

\draw(5.5,.5)node{\(-\)};

\draw(6,0)--(7,0)--(7,1)--(6,1)--(6,0);
\draw(6.25,.25)node{\(\)};
\draw(6.25,.75)node{\(0\)};
\draw(6.75,.25)node{\(\)};
\draw(6.75,.75)node{\(2\)};
\draw(6.5,0)node[below=-0.5ex]{\tiny\(\)};
\draw(7,.5)node[right=-0.5ex]{\tiny\(\)};
\draw(6.5,1)node[above=-0.5ex]{\tiny\(1\)};
\draw(6,.5)node[left=-0.5ex]{\tiny\(\)};
\draw(6,0)node[below left=-0.5ex]{\tiny\(\)};
\draw(7,0)node[below right=-0.5ex]{\tiny\(\)};
\draw(7,1)node[above right=-0.5ex]{\tiny\(\)};
\draw(6,1)node[above left=-0.5ex]{\tiny\(\)};

\draw(7.5,.5)node{\(-\)};

\draw(8,0)--(9,0)--(9,1)--(8,1)--(8,0);
\draw(8.25,.25)node{\(\)};
\draw(8.25,.75)node{\(1\)};
\draw(8.75,.25)node{\(\)};
\draw(8.75,.75)node{\(2\)};
\draw(8.5,0)node[below=-0.5ex]{\tiny\(\)};
\draw(9,.5)node[right=-0.5ex]{\tiny\(\)};
\draw(8.5,1)node[above=-0.5ex]{\tiny\(3\)};
\draw(8,.5)node[left=-0.5ex]{\tiny\(\)};
\draw(8,0)node[below left=-0.5ex]{\tiny\(\)};
\draw(9,0)node[below right=-0.5ex]{\tiny\(\)};
\draw(9,1)node[above right=-0.5ex]{\tiny\(\)};
\draw(8,1)node[above left=-0.5ex]{\tiny\(\)};

\draw(9.5,.5)node{\(+\)};

\draw(10,0)--(11,0)--(11,1)--(10,1)--(10,0);
\draw(10.25,.25)node{\(1\)};
\draw(10.25,.75)node{\(0\)};
\draw(10.75,.25)node{\(\)};
\draw(10.75,.75)node{\(\)};
\draw(10.5,0)node[below=-0.5ex]{\tiny\(\)};
\draw(11,.5)node[right=-0.5ex]{\tiny\(\)};
\draw(10.5,1)node[above=-0.5ex]{\tiny\(\)};
\draw(10,.5)node[left=-0.5ex]{\tiny\(0\)};
\draw(10,0)node[below left=-0.5ex]{\tiny\(\)};
\draw(11,0)node[below right=-0.5ex]{\tiny\(\)};
\draw(11,1)node[above right=-0.5ex]{\tiny\(\)};
\draw(10,1)node[above left=-0.5ex]{\tiny\(\)};

\draw(11.5,.5)node{\(+\)};

\draw(12,0)--(13,0)--(13,1)--(12,1)--(12,0);
\draw(12.25,.25)node{\(1\)};
\draw(12.25,.75)node{\(1\)};
\draw(12.75,.25)node{\(\)};
\draw(12.75,.75)node{\(\)};
\draw(12.5,0)node[below=-0.5ex]{\tiny\(\)};
\draw(13,.5)node[right=-0.5ex]{\tiny\(\)};
\draw(12.5,1)node[above=-0.5ex]{\tiny\(\)};
\draw(12,.5)node[left=-0.5ex]{\tiny\(1\)};
\draw(12,0)node[below left=-0.5ex]{\tiny\(\)};
\draw(13,0)node[below right=-0.5ex]{\tiny\(\)};
\draw(13,1)node[above right=-0.5ex]{\tiny\(\)};
\draw(12,1)node[above left=-0.5ex]{\tiny\(\)};

\draw[|->](-1.5,-2)--(-1.0,-2);

\draw(-.5,-2)node{\(-\)};

\draw(0,-3)--(1,-3)--(1,-1)--(0,-1)--(0,-3);
\draw(.25,-2)node{\(1\)};
\draw(.75,-2)node{\(0\)};
\draw(0,-2)node[left=-0.5ex]{\tiny\(1\)};
\draw(1,-2)node[right=-0.5ex]{\tiny\(0\)};

\draw(1.5,-2)node{\(+\)};

\draw(2,-3)--(3,-3)--(3,-1)--(2,-1)--(2,-3);
\draw(2.25,-2)node{\(2\)};
\draw(2.75,-2)node{\(1\)};
\draw(2,-2)node[left=-0.5ex]{\tiny\(2\)};
\draw(3,-2)node[right=-0.5ex]{\tiny\(1\)};

\draw(3.5,-2)node{\(-\)};

\draw(4,-2.5)--(6,-2.5)--(6,-1.5)--(4,-1.5)--(4,-2.5);
\draw(5,-2.25)node{\(1\)};
\draw(5,-1.75)node{\(0\)};
\draw(5,-2.5)node[below=-0.5ex]{\tiny\(1\)};
\draw(5,-1.5)node[above=-0.5ex]{\tiny\(0\)};

\draw(6.5,-2)node{\(-\)};

\draw(7,-3)--(8,-3)--(8,-1)--(7,-1)--(7,-3);
\draw(7.25,-2)node{\(0\)};
\draw(7.75,-2)node{\(2\)};
\draw(7,-2)node[left=-0.5ex]{\tiny\(0\)};
\draw(8,-2)node[right=-0.5ex]{\tiny\(2\)};

\draw(8.5,-2)node{\(-\)};

\draw(9,-3)--(10,-3)--(10,-1)--(9,-1)--(9,-3);
\draw(9.25,-2)node{\(1\)};
\draw(9.75,-2)node{\(2\)};
\draw(9,-2)node[left=-0.5ex]{\tiny\(1\)};
\draw(10,-2)node[right=-0.5ex]{\tiny\(2\)};

\draw(-.5,-4.5)node{\(+\)};

\draw(0,-5.0)--(2,-5.0)--(2,-4.0)--(0,-4.0)--(0,-5.0);
\draw(1,-4.75)node{\(1\)};
\draw(1,-4.25)node{\(0\)};
\draw(1,-5.0)node[below=-0.5ex]{\tiny\(1\)};
\draw(1,-4.0)node[above=-0.5ex]{\tiny\(0\)};

\draw(2.5,-4.5)node{\(+\)};

\draw(3,-5.0)--(5,-5.0)--(5,-4.0)--(3,-4.0)--(3,-5.0);
\draw(4,-4.75)node{\(1\)};
\draw(4,-4.25)node{\(1\)};
\draw(4,-5.0)node[below=-0.5ex]{\tiny\(1\)};
\draw(4,-4.0)node[above=-0.5ex]{\tiny\(1\)};
\end{tikzpicture}
\]
and is contractible, so \(\delta_0^1=0\). This is the same conclusion reached by Schur's lemma. Thus, working over reduced cohomology,
\begin{align*}
\check{H}^2(\Xi_1)&=\mathbb{Z}[1/4]\\
\check{H}^1(\Xi_1)&=\mathbb{Z}[1/4]^2\oplus\mathbb{Z}[1/2]\oplus\mathbb{Z}.
\end{align*}
\indent Towards the relative cohomology of a pair \((\Xi_2,\Xi_1)\), we have that
\[
\check{H}^2(\Xi_2,\Xi_1)=\varinjlim\left(\mathbb{Z}^3,\left(\begin{array}{*{3}c}
7&5&4\\
5&7&4\\
2&2&12
\end{array}\right)\right)
\]
with eigenvectors
\begin{align*}
\mathbf{e}_{16}&=(1,1,1)\\
\mathbf{e}_8&=(1,1,-1)\\
\mathbf{e}_2&=(1,-1,0).
\end{align*}
A straightforward calculation and the observation that \(S_2/S_1\) is a wedge of \(2\)-spheres give us
\begin{align*}
\check{H}^2(\Xi_2,\Xi_1)&=\mathbb{Z}[1/16]\oplus\mathbb{Z}[1/8]\oplus\mathbb{Z}[1/2]\\
\check{H}^1(\Xi_2,\Xi_1)&=0.
\end{align*}
\indent The long exact sequence in relative cohomology of a pair \((\Xi_2,\Xi_1)\) then reads
\[
\begin{tikzcd}
0\arrow[r]&\check{H}^0(\Xi_2,\Xi_1)\arrow[r]&\check{H}^0(\Xi_2)\arrow[r]\arrow[d,phantom,""{coordinate,name=X}]&\mathbb{Z}\arrow[lld,rounded corners,to path={--([xshift=2ex]\tikztostart.east)|-(X.center)\tikztonodes-|([xshift=-2ex]\tikztotarget.west)--(\tikztotarget)},swap,"\delta_1^0" at end]&\\
&0\arrow[r]&\check{H}^1(\Xi_2)\arrow[r]\arrow[d,phantom,""{coordinate,name=Y}]&\mathbb{Z}[1/4]^2\oplus\mathbb{Z}[1/2]\oplus\mathbb{Z}\arrow[lld,rounded corners,to path={--([xshift=2ex]\tikztostart.east)|-(Y.center)\tikztonodes-|([xshift=-2ex]\tikztotarget.west)--(\tikztotarget)},swap,"\delta_1^1" at end]&\\
&\mathbb{Z}[1/16]\oplus\mathbb{Z}[1/8]\oplus\mathbb{Z}[1/2]\arrow[r]&\check{H}^2(\Xi_2)\arrow[r]&\mathbb{Z}[1/4]\arrow[r]&0
\end{tikzcd}.
\]
By Schur's lemma, to show that \(\delta_1^1=0\), it suffices to check that \(\delta_1^1(0,1,-1)=0\), where \((0,1,-1)\) is the generator of \(\mathbb{Z}[1/2]\leq\check{H}^1(\Xi_1)\), which is straightforward. In fact, \(\delta_1^1\) applied to each of the generating loops in \(S_1/S_0\) (then lifted to a \(1\)-cochain in \(S_1\), modulo homotopy) yields as many \(2\)-flaps of the same type on one side of a loop as the other.\\
\indent In the very last step, we observe that \(\check{H}^2(\Xi_2,\Xi_1)\) and \(\check{H}^2(\Xi_1)\) have every element \(2\)-divisible. We first show that \(\check{H}^2(\Xi_2)\) is \(2\)-divisible as well. Consider any \(c\in\check{H}^2(\Xi_2)\). Let \(d\in\check{H}^2(\Xi_1)\) so that \(c\mapsto d\). By assumption, \(d/2\in\check{H}^2(\Xi_1)\). Let \(c'\in\check{H}^2(\Xi_2)\) so that \(c'\mapsto d/2\). Then \(c-2c'\mapsto 0\), and \(c-2c'\in\check{H}^2(\Xi_2,\Xi_1)\). Again, by assumption, \((c-2c')/2\in\check{H}^2(\Xi_2,\Xi_1)\leq\check{H}^2(\Xi_2)\). Then \(2((c-2c')/2+c')=c\), and \(\check{H}^2(\Xi_2)\) is \(2\)-divisible. To construct a splitting map \(\check{H}^2(\Xi_1)\rightarrow\check{H}^2(\Xi_2)\), take any \(c\in\check{H}^2(\Xi_2)\) so that \(c\mapsto 1\). The map is then defined by \(a/2^n\mapsto ac/2^n\), where \(a/2^n\in\check{H}^2(\Xi_1)\), and \(c/2^n\) exists since \(\check{H}^2(\Xi_2)\) is \(2\)-divisible.\\
\indent Thus the bottom short exact sequence splits, giving us that, over \(\mathbb{Z}\),
\begin{align*}
\check{H}^2(\Omega_\varsigma)&=\mathbb{Z}[1/16]\oplus\mathbb{Z}[1/8]\oplus\mathbb{Z}[1/4]\oplus\mathbb{Z}[1/2]\\
\check{H}^1(\Omega_\varsigma)&=\mathbb{Z}[1/4]^2\oplus\mathbb{Z}[1/2]\oplus\mathbb{Z}\\
\check{H}^0(\Omega_\varsigma)&=\mathbb{Z},
\end{align*}
which are exactly the cohomology groups of a product. Therefore, this is an example over \(\mathbb{Z}\) that the cup product structure discerns from a product space.\\
\indent Lastly, let us compute the frequency module. The substitution matrix on the \(2\)-cells of the dual complex is $\sigma = \sigma_{A4}$ which is found in Appendix \ref{app:1} and has Perron--Frobenius eigenvector
\[
\resizebox{0.9\hsize}{!}{
\((2512,2048,2048,2304,55,112,2352,121,696,2352,1792,4096,80,313,791,712,448,216,976,384,168).\)
}
\]
Its sum is \(24576\), giving us that the frequency module is \(\frac{1}{3}\mathbb{Z}[1/2]\), making it indistinguishable from the frequency module of the product.

For illustration, we also compute the frequency modules for substitution on the vertical and horizontal \(1\)-cells of the dual complex. There are two components to the vertical \(1\)-cells. The substitution matrices are then
$$
\sigma_{v_1}^1=\left(\begin{array}{*{1}c}
4
\end{array}\right),\;
\sigma_{v_2}^1=\left(\begin{array}{*{4}c}
1&1&1&1\\
1&1&1&1\\
0&0&0&0\\
2&2&2&2
\end{array}\right),\;\mbox{ and }\;
\sigma_h^1=\left(\begin{array}{*{9}c}
2&0&0&1&1&1&1&1&0\\
1&1&1&1&1&1&1&1&1\\
0&1&1&0&0&0&0&0&0\\
0&0&0&1&0&0&1&0&0\\
0&0&0&0&1&1&0&1&0\\
0&1&1&0&0&1&0&0&1\\
0&0&0&0&0&0&0&0&0\\
1&0&0&1&1&0&1&1&0\\
0&1&1&0&0&0&0&0&2
\end{array}\right).
$$
By replacing \(\mathbf{e}_{4_1}^1\) with \(\mathbf{e}_{4_1}^1+\mathbf{e}_{4_2}^1=(2,2,0,0)\), the Perron--Frobenius eigenvectors that are duals of lifts of the generators of \(\check{H}^1(\Omega_\varsigma)\) of eigenvalue \(4\) are
$(1), (1,1,0,2),$ and $(2,3,1,0,1,2,0,1,2)$
with the sum of the first two the lift of the new \(\mathbf{e}_{4_1}^1\). They have sums \(1\), \(4\), and \(12\), respectively, giving frequency modules \(\mathbb{Z}[1/2]\), \(\mathbb{Z}[1/2]\), and \(\frac{1}{3}\mathbb{Z}[1/2]\), respectively.
\end{example}

\begin{figure}[t]
  \includegraphics[width=3in]{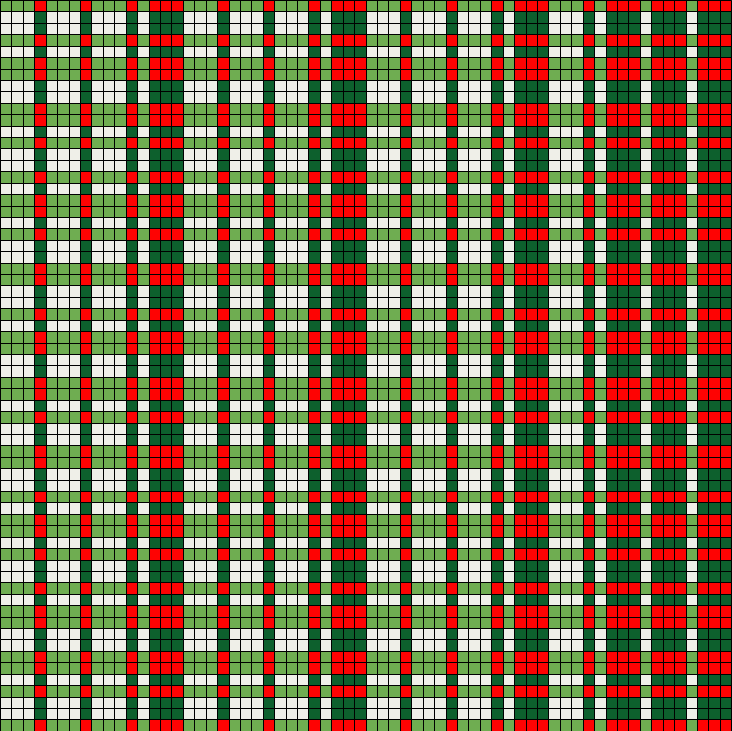} \hspace{.3in}
  \includegraphics[width=3in]{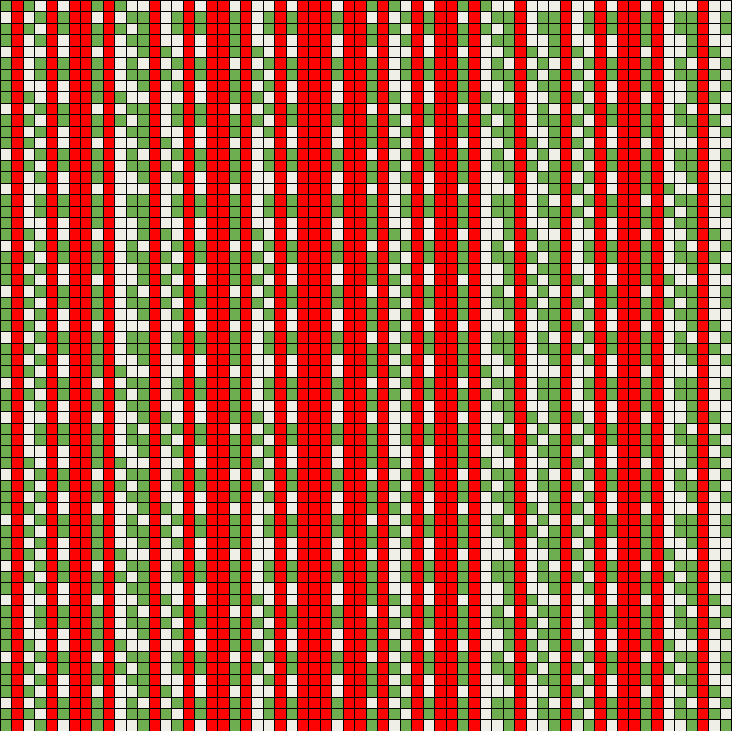}
    \caption{Patches from Examples \ref{ex:2} and \ref{ex:4}, respectively. That their tiling spaces are not homeomorphic can only be detected through the cup product.}
\end{figure}

\section{Breakdown of the Chern character in dimension four}
\label{sec:breakdown}
This section contains what we consider the main result of this paper:
a cubical substitution in dimension four where the Chern character
cannot factor as an integral isomorphism. The case of dimension four
is very special, since the complex $K$-theory and Chern character can
be completely determined from the cohomology ring \cite{K4}, even
though the $K$-theory ring and  cohomology ring may not be isomorphic.
An immediate consequence of \cite{K4} is that if $\Omega$ is a
$4$-dimensional complex and the Chern character on $K^0(\Omega)$ does
\emph{not} factor through $H^4(\Omega;\bZ)$, then there must be an
element $[c]$ of $H^2(\Omega;\bZ)$ whose cup square
$[c]^2=[c]\smile[c]$ is not divisible by $2$ in $H^4(\Omega;\bZ)$.
\subsection{The example}

\begin{figure}[t]
\centering
\begin{tikzpicture}[decoration=snake]
\draw(0,0)--(.5,0)--(.5,.5)--(0,.5)--(0,0);
\draw[|->](.75,.25)--(1.25,.25);
\draw(.25,.25)node{\(0\)};
\draw(1.5,-.5)--(2,-.5)--(2,1)--(1.5,1)--(1.5,-.5);
\draw(1.5,0)--(2,0);
\draw(1.5,.5)--(2,.5);
\draw(1.75,-.25)node{\(0\)};
\draw(1.75,.25)node{\(1\)};
\draw(1.75,.75)node{\(0\)};
\draw[|->](2.25,.25)--(2.75,.25);
\draw(3.0,-.5)--(4.5,-.5)--(4.5,1)--(3.0,1)--(3.0,-.5);
\draw(3.0,0)--(4.5,0);
\draw(3.0,.5)--(4.5,.5);
\draw(3.5,-.5)--(3.5,1);
\draw(4.0,-.5)--(4.0,1);
\draw(3.25,-.25)node{\(0\)};
\draw(3.25,.25)node{\(1\)};
\draw(3.25,.75)node{\(0\)};
\draw(3.75,-.25)node{\(2\)};
\draw(3.75,.25)node{\(3\)};
\draw(3.75,.75)node{\(2\)};
\draw(4.25,-.25)node{\(0\)};
\draw(4.25,.25)node{\(1\)};
\draw(4.25,.75)node{\(0\)};
\end{tikzpicture}
\caption{A two-dimensional analogue of the main example in section \ref{subsec:ex}. We begin with a one-dimensional checkerboard substitution on \(8\) prototiles (drawn vertically in the first map; only the first supertile is drawn, the rest are obtained by \(+1\mod 8\)), then extending it along the remaining coordinate (horizontally) by another (different) checkerboard (second map). The resulting tiling has rotational symmetry about axis the pattern is extended along (horizontal here).}
\label{fig:example}
\end{figure}%

\label{subsec:ex}
Let us consider the following four-dimensional expansion \(3\) substitution on eight prototiles.
\begin{minted}[linenos,fontsize=\fontsize{5}{6},numbersep=5pt,frame=lines,framesep=2mm,mathescape=true,escapeinside=||]{python}
|\(\varsigma\)|=[
[0,1,0,2,3,2,0,1,0,2,3,2,0,1,0,2,3,2,0,1,0,2,3,2,0,1,0,2,3,2,0,1,0,2,3,2,0,1,0,0,1,0,0,1,0,2,3,2,0,1,0,2,3,2,0,1,0,2,3,2,0,1,0,2,3,2,0,1,0,2,3,2,0,1,0,2,3,2,0,1,0],
[1,2,1,3,4,3,1,2,1,3,4,3,1,2,1,3,4,3,1,2,1,3,4,3,1,2,1,3,4,3,1,2,1,3,4,3,1,2,1,1,2,1,1,2,1,3,4,3,1,2,1,3,4,3,1,2,1,3,4,3,1,2,1,3,4,3,1,2,1,3,4,3,1,2,1,3,4,3,1,2,1],
[2,3,2,4,5,4,2,3,2,4,5,4,2,3,2,4,5,4,2,3,2,4,5,4,2,3,2,4,5,4,2,3,2,4,5,4,2,3,2,2,3,2,2,3,2,4,5,4,2,3,2,4,5,4,2,3,2,4,5,4,2,3,2,4,5,4,2,3,2,4,5,4,2,3,2,4,5,4,2,3,2],
[3,4,3,5,6,5,3,4,3,5,6,5,3,4,3,5,6,5,3,4,3,5,6,5,3,4,3,5,6,5,3,4,3,5,6,5,3,4,3,3,4,3,3,4,3,5,6,5,3,4,3,5,6,5,3,4,3,5,6,5,3,4,3,5,6,5,3,4,3,5,6,5,3,4,3,5,6,5,3,4,3],
[4,5,4,6,7,6,4,5,4,6,7,6,4,5,4,6,7,6,4,5,4,6,7,6,4,5,4,6,7,6,4,5,4,6,7,6,4,5,4,4,5,4,4,5,4,6,7,6,4,5,4,6,7,6,4,5,4,6,7,6,4,5,4,6,7,6,4,5,4,6,7,6,4,5,4,6,7,6,4,5,4],
[5,6,5,7,0,7,5,6,5,7,0,7,5,6,5,7,0,7,5,6,5,7,0,7,5,6,5,7,0,7,5,6,5,7,0,7,5,6,5,5,6,5,5,6,5,7,0,7,5,6,5,7,0,7,5,6,5,7,0,7,5,6,5,7,0,7,5,6,5,7,0,7,5,6,5,7,0,7,5,6,5],
[6,7,6,0,1,0,6,7,6,0,1,0,6,7,6,0,1,0,6,7,6,0,1,0,6,7,6,0,1,0,6,7,6,0,1,0,6,7,6,6,7,6,6,7,6,0,1,0,6,7,6,0,1,0,6,7,6,0,1,0,6,7,6,0,1,0,6,7,6,0,1,0,6,7,6,0,1,0,6,7,6],
[7,0,7,1,2,1,7,0,7,1,2,1,7,0,7,1,2,1,7,0,7,1,2,1,7,0,7,1,2,1,7,0,7,1,2,1,7,0,7,7,0,7,7,0,7,1,2,1,7,0,7,1,2,1,7,0,7,1,2,1,7,0,7,1,2,1,7,0,7,1,2,1,7,0,7,1,2,1,7,0,7]
]
\end{minted}
    The underlying idea is that in dimension four, we look for a substitution with a torsion term in \(\check{H}^4(\Omega_\varsigma;\mathbb{Z})\). Furthermore, since the Chern character in dimension four involves a division by \(2\), for it to not be an isomorphism over \(\mathbb{Z}\), we want \(\check{H}^4(\Omega_\varsigma;\mathbb{Z})\) to contain even torsion, since odd torsion has multiplication by \(2\) invertible.\\
\indent This example is built from checkerboard patterns (Figure \ref{fig:example}). It appears that although checkerboard patterns do not necessarily force the border (see Remark \ref{rem:force}), the level of supertiles required to exhaust all \(1\)-collared\footnote{Although the top-dimensional cells are ``half-collars'', the construction of the complex requires knowledge of how they are identified in the tiling, thus another ``half-collar'' of the half-collared prototiles, i.e. \(1\)-collared prototiles.} prototiles to construct the dual complex is not too large. Since we want even torsion, we create a checkerboard pattern on \([0,3]^3\), then extend in the last dimension with yet another checkerboard pattern that depends on the prototile in the pattern on \([0,3]^3\), taking the result to be the supertiles, so as to introduce four-fold rotational symmetry with odd expansion. For this particular example, the dual complex has \(1120\) \(4\)-cells, \(1232\) \(3\)-cells, \(480\) \(2\)-cells, \(88\) \(1\)-cells, and \(8\) \(0\)-cells. One can get by with fewer prototiles or simpler patterns (see the next subsection), but this is the first example we discovered where everything we desire is evident.\\
\indent We only focus on \(\check{H}^2(\Omega_\varsigma;\mathbb{Z})\) and \(\check{H}^4(\Omega_\varsigma;\mathbb{Z})\), since we want to identify a class \([c]\in\check{H}^2(\Omega_\varsigma;\mathbb{Z})\) so that \([c]^2=[c]\smile[c]\in\check{H}^4(\Omega_\varsigma;\mathbb{Z})\) is not (uniquely) divisible by \(2\) (see (\ref{eqn:chern})). \(\check{H}^2(\Omega_\varsigma;\mathbb{Z})\) is the direct limit computed via the matrix $\sigma^2= \sigma^2_{A5}$, which is found in Appendix \ref{app:1}, and \(\check{H}^4(\Omega_\varsigma;\mathbb{Z})\) is computed via a matrix too large to be listed here (\(354\times 354\)), but it suffices to know that for the dual complex, \(\check{H}^4\cong \mathbb{Z}_4^3\oplus F\) for some torsion-free group $F$. Since we only care about the even torsion terms, there being $\mathbb{Z}_4^3$ in \(\check{H}^4\) of the dual complex, we write the substitution action on the dual complex as
\begin{equation}
  \label{eqn:sigma4}
\sigma^4 = \left(\begin{array}{c|c}\sigma^4_{\textnormal{torsion}}&\ast\\\hline 0&\sigma^4_{\textnormal{torsion-free}}\end{array}\right) = \left(
\begin{array}{@{}c|c@{}}
\begin{array}{ccc}
1&0&0\\
0&2&3\\
0&3&2
\end{array}
&\ast\\\hline 0&\sigma^4_{\textnormal{torsion-free}}
\end{array}
\right),
\end{equation}
where $\sigma^4_{\textnormal{torsion}}$ is the action restricted to the torsion part of $\check{H}^4$.

\indent Let us consider the eigenvector of \(\sigma^2\) with eigenvalue \(3\)
\[
v=(0,0,0,0,1,1,0,0,1,0,0,0,0,0,-1,1,1,0,0,0,0,0,0,0).
\]
Lifting it to a \(2\)-cochain, cupping it with itself via the cubical cup product formula given in \cite{KM:ring}, then projecting to \(\check{H}^4\) of the dual complex gives the class
\[
w=v\smile v = (2,0,2,\cdots)\in\mathbb{Z}_4^3\oplus F
\]
which is an eigenvector with eigenvalue \(9\). By (\ref{eqn:sigma4}), the
first copy of $\mathbb{Z}_4$ is preserved in the direct limit. Since the
entry of \(w\) in this copy of $\mathbb{Z}_4$ is \(2\), \(w=v\smile v\) is
not (uniquely) divisible by \(2\), and the Chern character cannot be an isomorphism
over \(\mathbb{Z}\).  Incidentally, this provides a negative answer to
the question posed on the last line of \cite{MR4030282}.

Following the same procedure as before, the Perron--Frobenius
eigenvector of the substitution matrix on the \(4\)-cells of the dual
complex sums to \(29952\), giving us that the frequency module for
this cubical substitution tiling is $\frac{1}{3328}\mathbb{Z}[1/3]$. 

\subsubsection{Additional examples}
We have discovered numerous other examples with the same features as the example
above. Here we mention two simple ones.

The first example consists on a cubical substitution on three
prototiles and expansion 3. Although this substitution rule is easier
to define, the cohomological computations are more
labor-intensive. The substitution rule is listed below and is the
analogue of the two-dimensional substitution rule found in Figure
\ref{fig:ex2}. We leave the details to the reader. 

\begin{figure}[t]
\centering
\begin{subfigure}{\textwidth}
\centering
\begin{tikzpicture}
\draw(0,0)--(.5,0)--(.5,.5)--(0,.5)--(0,0);
\draw[|->](.75,.25)--(1.25,.25);
\draw(.25,.25)node{\(0\)};
\draw(1.5,-.5)--(3,-.5)--(3,1)--(1.5,1)--(1.5,-.5);
\draw(1.5,0)--(3,0);
\draw(1.5,.5)--(3,.5);
\draw(2,-.5)--(2,1);
\draw(2.5,-.5)--(2.5,1);
\draw(1.75,-.25)node{\(1\)};
\draw(1.75,.25)node{\(0\)};
\draw(1.75,.75)node{\(1\)};
\draw(2.25,-.25)node{\(0\)};
\draw(2.25,.25)node{\(0\)};
\draw(2.25,.75)node{\(0\)};
\draw(2.75,-.25)node{\(1\)};
\draw(2.75,.25)node{\(0\)};
\draw(2.75,.75)node{\(1\)};

\draw(4.0,0)--(4.5,0)--(4.5,.5)--(4.0,.5)--(4.0,0);
\draw[|->](4.75,.25)--(5.25,.25);
\draw(4.25,.25)node{\(1\)};
\draw(5.5,-.5)--(7.0,-.5)--(7.0,1)--(5.5,1)--(5.5,-.5);
\draw(5.5,0)--(7.0,0);
\draw(5.5,.5)--(7.0,.5);
\draw(6.0,-.5)--(6.0,1);
\draw(6.5,-.5)--(6.5,1);
\draw(5.75,-.25)node{\(0\)};
\draw(5.75,.25)node{\(1\)};
\draw(5.75,.75)node{\(0\)};
\draw(6.25,-.25)node{\(1\)};
\draw(6.25,.25)node{\(1\)};
\draw(6.25,.75)node{\(1\)};
\draw(6.75,-.25)node{\(0\)};
\draw(6.75,.25)node{\(1\)};
\draw(6.75,.75)node{\(0\)};
\end{tikzpicture}
\caption{}
\label{fig:squiral}
\end{subfigure}

\hfill

\begin{subfigure}{\textwidth}
\centering
\begin{tikzpicture}
\draw(0,0)--(.5,0)--(.5,.5)--(0,.5)--(0,0);
\draw[|->](.75,.25)--(1.25,.25);
\draw(.25,.25)node{\(0\)};
\draw(1.5,-.5)--(3,-.5)--(3,1)--(1.5,1)--(1.5,-.5);
\draw(1.5,0)--(3,0);
\draw(1.5,.5)--(3,.5);
\draw(2,-.5)--(2,1);
\draw(2.5,-.5)--(2.5,1);
\draw(1.75,-.25)node{\(0\)};
\draw(1.75,.25)node{\(1\)};
\draw(1.75,.75)node{\(0\)};
\draw(2.25,-.25)node{\(1\)};
\draw(2.25,.25)node{\(1\)};
\draw(2.25,.75)node{\(1\)};
\draw(2.75,-.25)node{\(0\)};
\draw(2.75,.25)node{\(1\)};
\draw(2.75,.75)node{\(0\)};

\draw(4.0,0)--(4.5,0)--(4.5,.5)--(4.0,.5)--(4.0,0);
\draw[|->](4.75,.25)--(5.25,.25);
\draw(4.25,.25)node{\(1\)};
\draw(5.5,-.5)--(7.0,-.5)--(7.0,1)--(5.5,1)--(5.5,-.5);
\draw(5.5,0)--(7.0,0);
\draw(5.5,.5)--(7.0,.5);
\draw(6.0,-.5)--(6.0,1);
\draw(6.5,-.5)--(6.5,1);
\draw(5.75,-.25)node{\(1\)};
\draw(5.75,.25)node{\(2\)};
\draw(5.75,.75)node{\(1\)};
\draw(6.25,-.25)node{\(2\)};
\draw(6.25,.25)node{\(2\)};
\draw(6.25,.75)node{\(2\)};
\draw(6.75,-.25)node{\(1\)};
\draw(6.75,.25)node{\(2\)};
\draw(6.75,.75)node{\(1\)};

\draw(8,0)--(8.5,0)--(8.5,.5)--(8,.5)--(8,0);
\draw[|->](8.75,.25)--(9.25,.25);
\draw(8.25,.25)node{\(2\)};
\draw(9.5,-.5)--(11,-.5)--(11,1)--(9.5,1)--(9.5,-.5);
\draw(9.5,0)--(11,0);
\draw(9.5,.5)--(11,.5);
\draw(10,-.5)--(10,1);
\draw(10.5,-.5)--(10.5,1);
\draw(9.75,-.25)node{\(2\)};
\draw(9.75,.25)node{\(0\)};
\draw(9.75,.75)node{\(2\)};
\draw(10.25,-.25)node{\(0\)};
\draw(10.25,.25)node{\(0\)};
\draw(10.25,.75)node{\(0\)};
\draw(10.75,-.25)node{\(2\)};
\draw(10.75,.25)node{\(0\)};
\draw(10.75,.75)node{\(2\)};
\end{tikzpicture}
\caption{}
\label{fig:ex2}
\end{subfigure}
\caption{\subref{fig:squiral} The substitution rule that produces the squiral tiling. \subref{fig:ex2} The four-dimensional version of a three-prototile variant of the squiral substitution has the same features as the example in (\ref{eqn:sigma4}).}
\end{figure}%


\begin{minted}[linenos,fontsize=\fontsize{5}{6},numbersep=5pt,frame=lines,framesep=2mm,mathescape=true,escapeinside=||]{python}
|\(\varsigma\)|=[
[0,1,0,1,0,1,0,1,0,1,0,1,0,1,0,1,0,1,0,1,0,1,0,1,0,1,0,1,0,1,0,1,0,1,0,1,0,1,0,1,1,1,0,1,0,1,0,1,0,1,0,1,0,1,0,1,0,1,0,1,0,1,0,1,0,1,0,1,0,1,0,1,0,1,0,1,0,1,0,1,0],
[1,2,1,2,1,2,1,2,1,2,1,2,1,2,1,2,1,2,1,2,1,2,1,2,1,2,1,2,1,2,1,2,1,2,1,2,1,2,1,2,2,2,1,2,1,2,1,2,1,2,1,2,1,2,1,2,1,2,1,2,1,2,1,2,1,2,1,2,1,2,1,2,1,2,1,2,1,2,1,2,1],
[2,0,2,0,2,0,2,0,2,0,2,0,2,0,2,0,2,0,2,0,2,0,2,0,2,0,2,0,2,0,2,0,2,0,2,0,2,0,2,0,0,0,2,0,2,0,2,0,2,0,2,0,2,0,2,0,2,0,2,0,2,0,2,0,2,0,2,0,2,0,2,0,2,0,2,0,2,0,2,0,2]
]
\end{minted}

The second example is a cubical substitution on two prototiles and expansion 3. Its construction is much more involved, which we elect to forgo, and does not have two-dimensional analogues. We again leave the computational details to the reader.

\begin{minted}[linenos,fontsize=\fontsize{5}{6},numbersep=5pt,frame=lines,framesep=2mm,mathescape=true,escapeinside=||]{python}
|\(\varsigma\)|=[
[0,0,0,0,0,1,0,1,0,0,0,1,0,0,0,1,0,0,0,1,0,1,0,0,0,0,0,0,0,1,0,0,0,1,0,0,0,0,0,0,1,1,0,1,0,1,0,0,0,1,0,0,0,1,0,1,0,1,0,0,0,0,0,1,0,0,0,1,0,0,0,1,0,0,0,0,0,1,0,1,0],
[1,1,1,1,1,0,1,0,1,1,1,0,1,1,1,0,1,1,1,0,1,0,1,1,1,1,1,1,1,0,1,1,1,0,1,1,1,1,1,1,0,0,1,0,1,0,1,1,1,0,1,1,1,0,1,0,1,0,1,1,1,1,1,0,1,1,1,0,1,1,1,0,1,1,1,1,1,0,1,0,1]
]
\end{minted}

\begin{figure}[t]
\centering
\begin{tikzpicture}[decoration=snake]
\draw(.25,.25)node{\(0\)};
\draw(.25,.75)node{\(0\)};
\draw(.75,.25)node{\(0\)};
\draw(.75,.75)node{\(1\)};
\draw(0,0)--(1,0)--(1,1)--(0,1)--(0,0);
\draw(.5,0)--(.5,1);
\draw(0,.5)--(1,.5);
\fill[draw opacity=0.5,fill opacity=0,pattern=north east lines](.25,.25)--(.75,.25)--(.75,.75)--(.25,.75)--(.25,.25);
\draw(.25,.25)--(.75,.25)--(.75,.75)--(.25,.75)--(.25,.25);
\draw[|->](1.25,.5)--(1.75,.5)node[midway,above]{\(\sim\)};
\draw(2.25,.25)node{\(1\)};
\draw(2.25,.75)node{\(1\)};
\draw(2.75,.25)node{\(1\)};
\draw(2.75,.75)node{\(0\)};
\draw(2,0)--(3,0)--(3,1)--(2,1)--(2,0);
\draw(2.5,0)--(2.5,1);
\draw(2,.5)--(3,.5);
\fill[draw opacity=0.5,fill opacity=0,pattern=north east lines](2.25,.25)--(2.75,.25)--(2.75,.75)--(2.25,.75)--(2.25,.25);
\draw(2.25,.25)--(2.75,.25)--(2.75,.75)--(2.25,.75)--(2.25,.25);
\end{tikzpicture}
\caption{Two ``half-collared'' prototiles in the dual complex where quotienting by the \(\mathbb{Z}_2\)-action that interchanges the underlying prototiles identifies them.}
\label{fig:squiral-z/2}
\end{figure}%

\indent\indent All of the previous examples had a torsion-free generator in \(\check{H}^2\) squaring to twice a torsion generator in \(\check{H}^4\). To demonstrate that \(\mathbb{RP}^4\)-like behavior can also occur, that is, a torsion generator in \(\check{H}^2\) can square to twice a torsion generator in \(\check{H}^4\), we first give some motivations.\\
\indent\indent It turns out that the first example we provided in this section is a four-dimensional analogue of the three-prototile version of the squiral substitution (encoded slightly differently than usual). In \(d=2\), the squiral tiling was shown in \cite{MR3227148} to be a substitution tiling with expansion \(3\) on squares using two prototiles (Figure \ref{fig:squiral}), and it was computed there that the tiling has singular continuous spectrum. If one writes down the induced substitution rule on collared prototiles (from the dual complex), one can impose a \(\mathbb{Z}_2\)-action that interchanges the two prototiles (Figure \ref{fig:squiral-z/2}). It turns out that quotienting the squiral tiling space by this action results in another substitution tiling space, but with \(\mathbb{Z}_2\)-torsion in \(\check{H}^2\) when there was no torsion in the original squiral tiling.\\
\indent Motivated by this, let us consider the following generalization of the squiral substitution to arbitrary \(d\) that is different from the one in the first example.

\begin{minted}[linenos,fontsize=\fontsize{5}{6},numbersep=5pt,frame=lines,framesep=2mm,mathescape=true,escapeinside=||]{python}
|\(\varsigma\)|=[[(i+1)%2 if 1 not in t else i%2 for t in list(cartesian_product([range(3)]*d))] for i in range(2)]
\end{minted}

\noindent For \(d=4\), the rule reads as the following.

\begin{minted}[linenos,fontsize=\fontsize{5}{6},numbersep=5pt,frame=lines,framesep=2mm,mathescape=true,escapeinside=||]{python}
|\(\varsigma\)|=[
[1,0,1,0,0,0,1,0,1,0,0,0,0,0,0,0,0,0,1,0,1,0,0,0,1,0,1,0,0,0,0,0,0,0,0,0,0,0,0,0,0,0,0,0,0,0,0,0,0,0,0,0,0,0,1,0,1,0,0,0,1,0,1,0,0,0,0,0,0,0,0,0,1,0,1,0,0,0,1,0,1],
[0,1,0,1,1,1,0,1,0,1,1,1,1,1,1,1,1,1,0,1,0,1,1,1,0,1,0,1,1,1,1,1,1,1,1,1,1,1,1,1,1,1,1,1,1,1,1,1,1,1,1,1,1,1,0,1,0,1,1,1,0,1,0,1,1,1,1,1,1,1,1,1,0,1,0,1,1,1,0,1,0]
]
\end{minted}

\indent We computed there to be \(478\) ``half-collared'' prototiles, and quotienting by the \(\mathbb{Z}_2\)-action resulted in \(239\) prototiles with an associated induced substitution rule, \(\overline{\varsigma}\), that is too large to be written here. Fortunately, since the induced substitution on the half-collared prototiles of the squiral tiling forces the border (see Remark \ref{rem:force}), as does the resulting induced substitution rule, thus it suffices to work with the uncollared \(AP\)-complex, which we denote \(\overline{\Gamma}\) (reserving \(\Gamma\) for the original substitution rule \(\varsigma\)). It has \(239\) \(4\)-cells, \(160\) \(3\)-cells, \(48\) \(2\)-cells, \(8\) \(1\)-cells, and \(1\) \(0\)-cell.\\
\indent We have that
\begin{align*}
\check{H}^2(\overline{\Gamma};\mathbb{Z})&=\mathbb{Z}_2\oplus\mathbb{Z}^9\\
\check{H}^4(\overline{\Gamma};\mathbb{Z})&=\mathbb{Z}_2^{14}\oplus\mathbb{Z}_4\oplus\mathbb{Z}^{126}.
\end{align*}
The induced substitution matrices are
\begin{align*}
\sigma^2&=\left(\begin{array}{c|c}\sigma^2_\textnormal{torsion}&\ast\\\hline 0&\sigma^2_\textnormal{torsion-free}\end{array}\right)=\left(\begin{array}{c|c}1&\ast\\\hline 0&\sigma^2_\textnormal{torsion-free}\end{array}\right)\\
\sigma^4&= \left(\begin{array}{c|c}\sigma^4_{\textnormal{torsion}}&\ast\\\hline 0&\sigma^4_{\textnormal{torsion-free}}\end{array}\right) = \left(
\begin{array}{c|c}\textnormal{Id}_{15}&\ast\\\hline 0&\sigma^4_{\textnormal{torsion-free}}\end{array}\right),
\end{align*}
so the torsion terms stay in the limit and belong to the respective cohomology groups of \(\Omega_{\overline{\varsigma}}\).\\
\indent Let
\[
v=(1,0,0,0,0,0,0,0,0,0)\in\check{H}^2(\overline{\Gamma};\mathbb{Z}).
\]
Lifting it to a \(2\)-cochain, cupping it with itself, then projecting to \(\check{H}^4(\overline{\Gamma};\mathbb{Z})\) gives the class
\[
w=v\smile v=(0,0,0,0,0,0,0,0,0,0,0,0,0,0,2,0,\ldots)\in\check{H}^2(\overline{\Gamma};\mathbb{Z})
\]
where the single nonzero entry corresponds exactly to twice the generator of \(\mathbb{Z}_4\leq\check{H}^4(\overline{\Gamma};\mathbb{Z})\).

\subsection{Gap labels}
What does the example above say about the gap-labeling conjecture?
What it \emph{does not} say is that it is false: the example above
relied heavily on the cohomology of the tiling space having
non-trivial torsion terms, which were the source of the breakdown of
the Chern character map. However, \textbf{the torsion is not detected
  by the Ruelle-Sullivan current $C_\mu$}, a homomorphism to
$\mathbb{R}$, meaning that the torsion is not detected in the
frequency module. In addition, it is not clear that the image of the
trace map on $K_0(\mathcal{A}_p(\Omega))$ is contained in the image of
the trace map on $K_0(\mathcal{A}_{AF}(\Omega))$, so we cannot compare
$\tau(K_0(\mathcal{A}_p(\Omega)))$ with the frequency module, as we
only know how to compute the image of the tracial state of the
AF-algebra $\mathcal{A}_{AF}(\Omega)$. 

What it \emph{does} say is that the gap in several papers relying on
the Chern character factoring through integral cohomology to give
(\ref{eqn:wishful}) is real---our example shows it. We do not believe
that a true counterexample to the gap-labeling conjecture---which
would only exist in dimensions greater than three---will be found
using cubical substitutions, and so the problem of finding a
counterexample becomes much harder, as one needs to define a
non-cubical substitution rule in four dimensions (or higher) and
compute its cohomology ring. This task seems out of reach at the
moment. 

\subsection{An aperiodic counterexample to equivariant gap-labeling in dimension 4}
\label{subsec:equivCounter}

\begin{example}
To provide an additional counterexample for Conjecture
\ref{conj:equivlabeling} that is a substitution tiling, we construct
the simplest possible one that is not a solenoid by building off of
Theorem \ref{thm:equivlabeling}. We aim for its complex to have the
same \(3\)-skeleton as \(T^4\), but with an additional \(4\)-cell
attached in the same way as the original \(4\)-cell in \(T^4\).\\ 
\indent To do so, consider the following four-dimensional substitution with expansion 3 on two prototiles.

\begin{minted}[linenos,fontsize=\fontsize{5}{6},numbersep=5pt,frame=lines,framesep=2mm,mathescape=true,escapeinside=||]{python}
|\(\varsigma\)|=[
[0,0,0,0,0,0,0,0,0,0,0,0,0,0,0,0,0,0,0,0,0,0,0,0,0,0,0,0,0,0,0,0,0,0,0,0,0,0,0,0,1,0,0,0,0,0,0,0,0,0,0,0,0,0,0,0,0,0,0,0,0,0,0,0,0,0,0,0,0,0,0,0,0,0,0,0,0,0,0,0,0],
[0,0,0,0,0,0,0,0,0,0,0,0,0,0,0,0,0,0,0,0,0,0,0,0,0,0,0,0,0,0,0,0,0,0,0,0,0,0,0,0,0,0,0,0,0,0,0,0,0,0,0,0,0,0,0,0,0,0,0,0,0,0,0,0,0,0,0,0,0,0,0,0,0,0,0,0,0,0,0,0,0]
]
\end{minted}

\noindent Its two-dimensional analogue is Figure \ref{fig:ex3}.\\
\indent This substitution is primitive and recognizable, and clearly
forces the border (see Remark \ref{rem:force}), thus we can use the \emph{uncollared}
\(AP\)-complex instead, which is the \(AP\)-construction without using
collared prototiles. By the same theorem in \cite{AP}, its inverse
limit is homeomorphic to the tiling space.\\ 
\indent The \(2^4\) patch of all \(0\)'s indicates that the complex
contains a copy of \(T^4\). The prototile \(1\) being present in all
possible positions in a \(2^4\) patch asserts that its associated
\(4\)-cell is attached to the \(3\)-skeleton in the same exact way as
the \(4\)-cell associated to the prototile \(0\).\\ 
\indent All of the coboundary maps are trivial, giving us that the cohomology groups we are interested in are
\begin{align*}
\check{H}^4(\Omega_\varsigma;\mathbb{Z})&=\varinjlim\left(\mathbb{Z}^2,\left(\begin{array}{cc} 
80&1\\
81&0
\end{array}\right)\right)\\
\check{H}^2(\Omega_\varsigma;\mathbb{Z})&=\varinjlim\left(\mathbb{Z}^6,9\cdot\textnormal{Id}\right),
\end{align*}
\noindent where \(\mathbb{Z}^2\) is simultaneously \(C^4\) and \(H^4\)
of the \(AP\)-complex, and \(\mathbb{Z}^6\) is simultaneously \(C^2\)
and \(H^2\) of the \(AP\)-complex. Due to this, the cubical cup
product on the cohomology ring coincides with the cubical cup product
on the cochains, which is easy to describe. For example, the two
generators
\(c_1=[0,1]\times[0,1]\times[0]\times[0],c_2=[0]\times[0]\times[0,1]\times[0,1]\in
C^2\) cup to \((1,1)\in C^4\), the sum of the duals of the two
prototiles. That is, the cohomology ring structure of this
\(AP\)-complex only differs from that of \(T^4\) by the cup product
always witnessing both \(4\)-cochains whenever it is nontrivial. 

\begin{figure}[t]
\centering
\begin{tikzpicture}
\draw(0,0)--(.5,0)--(.5,.5)--(0,.5)--(0,0);
\draw[|->](.75,.25)--(1.25,.25);
\draw(.25,.25)node{\(0\)};
\draw(1.5,-.5)--(3,-.5)--(3,1)--(1.5,1)--(1.5,-.5);
\draw(1.5,0)--(3,0);
\draw(1.5,.5)--(3,.5);
\draw(2,-.5)--(2,1);
\draw(2.5,-.5)--(2.5,1);
\draw(1.75,-.25)node{\(0\)};
\draw(1.75,.25)node{\(0\)};
\draw(1.75,.75)node{\(0\)};
\draw(2.25,-.25)node{\(0\)};
\draw(2.25,.25)node{\(1\)};
\draw(2.25,.75)node{\(0\)};
\draw(2.75,-.25)node{\(0\)};
\draw(2.75,.25)node{\(0\)};
\draw(2.75,.75)node{\(0\)};

\draw(4.0,0)--(4.5,0)--(4.5,.5)--(4.0,.5)--(4.0,0);
\draw[|->](4.75,.25)--(5.25,.25);
\draw(4.25,.25)node{\(1\)};
\draw(5.5,-.5)--(7.0,-.5)--(7.0,1)--(5.5,1)--(5.5,-.5);
\draw(5.5,0)--(7.0,0);
\draw(5.5,.5)--(7.0,.5);
\draw(6.0,-.5)--(6.0,1);
\draw(6.5,-.5)--(6.5,1);
\draw(5.75,-.25)node{\(0\)};
\draw(5.75,.25)node{\(0\)};
\draw(5.75,.75)node{\(0\)};
\draw(6.25,-.25)node{\(0\)};
\draw(6.25,.25)node{\(0\)};
\draw(6.25,.75)node{\(0\)};
\draw(6.75,-.25)node{\(0\)};
\draw(6.75,.25)node{\(0\)};
\draw(6.75,.75)node{\(0\)};
\end{tikzpicture}
\caption{The four-dimensional version of this substitution rule is the same \(3\)-skeleton as \(T^4\), with the only difference being that there are two \(4\)-cells (the two prototiles) attached to the \(3\)-cells.}
\label{fig:ex3}
\end{figure}%

One then follows the same exact argument provided in Theorem
\ref{thm:equivlabeling} and concludes that this substitution tiling
space, upon quotienting by the same \(\mathbb{Z}_2\)-action, yields a
counterexample to Conjecture \ref{conj:equivlabeling}. To be precise,
let $\Gamma$ be the AP complex of the substitution described above,
and let $\rho:\Gamma\rightarrow \mathbb{T}^4$ be the factor map which
collapses the 4 cells. There is an involution on $\Gamma$ which is
equivariant (under $\rho$) with the involution of $\mathbb{T}^4 =
\mathbb{T}^2\times\mathbb{T}^2$ used in Theorem
\ref{thm:equivlabeling}. This involution is defined by exchanging the
2-cells $c_1$ and $c_2$ described above. The induced map on cohomology
$\rho^*\co H^4(\mathbb{T}^4;\mathbb{Z})\rightarrow H^4(\Gamma;\mathbb{Z})$
sends the generator in
$H^4(\mathbb{T}^4;\mathbb{Z})$ to $(1,1)\in H^4(\Gamma;\mathbb{Z})$.
Let $S_4 = \varprojlim (\mathbb{T}^4,3\cdot \mathrm{Id})$
be the four dimensional solenoid constructed with maps
of expansion 3, and note that it is a factor of the tiling space
$\Omega$ corresponding to the substitution above. Then the argument
from Theorem \ref{thm:equivlabeling} carries over to $S_4$ in the
direct limit since the expansion is by 3, and it can be pulled back to
$H^*(\Omega;\mathbb{Z})$ using the map induced by the factor map
(because the expansions of both systems are by 3 and thus cannot be
divided by 2 by the Chern character). Thus this example provides an
aperiodic counterexample to Conjecture \ref{conj:equivlabeling}. 

\indent Lastly, let us remark on the mechanism that yields these
counterexamples. In the original \(AP\)-complex, the nontrivial
squares are of the form \((c_1+c_2)\smile(c_1+c_2)=c_1\smile
c_1+c_1\smile c_2+c_2\smile c_1+c_2\smile c_2=2c_1\smile c_2\), which
always return twice the sum of the generators in \(H^4\), and
therefore still yields integrality of the Chern character. In this
\(\mathbb{Z}_2\) quotient, \(c_1\) and \(c_2\) are identified, so
\(c_1\smile c_2=c_1\smile c_1\) (in the quotient), and we no longer
need to add cochains together to produce a nontrivial square in
$H^4$. This breaks the integrality.\\ 
\indent Does there exist other (\(\mathbb{Z}_2\)) actions that produce the same effect? There are a few obvious choices that do not appear to work.
\begin{itemize}
\item
If the substitution is on two prototiles, let the \(\mathbb{Z}_2\)-action be swapping the two prototiles, assuming that the supertiles are related by \(+1\mod 2\) (so that the action is well-defined). This does not appear to work, because ``horizontal'' cells remain ``horizontal'', and ``vertical'' cells remain ``vertical'', thus to obtain a nontrivial square, one still needs to add together both ``horizontal'' and ``vertical'', resulting in a factor of two on the \(4\)-cochain.
\item
If the substitution has its supertiles symmetric about, say, the axis along \((1,0,0,0)\), let the \(\mathbb{Z}_2\)-action be flipping the axis. For the same reason as above, this also does not appear to work.
\end{itemize}
\end{example}

\newpage 
\appendix
\section{Large matrices}
\label{app:1}
\[
\resizebox{0.9\hsize}{!}{
\(\sigma_{A1}=\left(\begin{array}{*{36}c}
3&3&3&3&4&4&3&3&3&3&4&3&3&3&4&3&3&3&4&3&3&3&4&3&3&3&4&3&3&3&3&3&4&4&3&3\\
0&0&0&0&0&0&0&0&0&0&0&0&3&3&4&3&0&0&0&0&0&0&0&0&3&3&4&3&3&0&3&0&4&0&3&0\\
1&1&1&1&0&0&1&1&1&1&0&1&1&1&0&1&1&1&0&1&1&1&0&1&1&1&0&1&1&1&1&1&0&0&1&1\\
0&0&0&0&0&0&0&0&0&0&0&0&1&1&0&1&0&0&0&0&0&0&0&0&1&1&0&1&1&0&1&0&0&0&1&0\\
1&1&0&0&1&1&1&1&1&0&1&1&1&0&1&1&1&0&1&1&1&0&1&1&1&0&1&1&1&1&0&0&1&1&1&1\\
0&0&0&0&0&0&0&0&0&0&0&0&1&0&1&1&0&0&0&0&0&0&0&0&1&0&1&1&1&0&0&0&1&0&1&0\\
0&0&1&1&0&0&0&0&0&1&0&0&0&1&0&0&0&1&0&0&0&1&0&0&0&1&0&0&0&0&1&1&0&0&0&0\\
0&0&0&0&0&0&0&0&0&0&0&0&0&1&0&0&0&0&0&0&0&0&0&0&0&1&0&0&0&0&1&0&0&0&0&0\\
3&3&3&3&4&4&3&3&3&3&4&3&0&0&0&0&3&3&4&3&3&3&4&3&0&0&0&0&0&3&0&3&0&4&0&3\\
1&1&1&1&0&0&1&1&1&1&0&1&0&0&0&0&1&1&0&1&1&1&0&1&0&0&0&0&0&1&0&1&0&0&0&1\\
1&1&0&0&1&1&1&1&1&0&1&1&0&0&0&0&1&0&1&1&1&0&1&1&0&0&0&0&0&1&0&0&0&1&0&1\\
0&0&1&1&0&0&0&0&0&1&0&0&0&0&0&0&0&1&0&0&0&1&0&0&0&0&0&0&0&0&0&1&0&0&0&0\\
0&0&0&0&0&0&0&0&0&0&0&0&6&6&8&6&0&0&0&0&0&0&0&0&6&6&8&6&6&0&6&0&8&0&6&0\\
0&0&0&0&0&0&0&0&0&0&0&0&2&2&0&2&0&0&0&0&0&0&0&0&2&2&0&2&2&0&2&0&0&0&2&0\\
0&0&0&0&0&0&0&0&0&0&0&0&2&0&2&2&0&0&0&0&0&0&0&0&2&0&2&2&2&0&0&0&2&0&2&0\\
0&0&0&0&0&0&0&0&0&0&0&0&0&2&0&0&0&0&0&0&0&0&0&0&0&2&0&0&0&0&2&0&0&0&0&0\\
3&3&3&3&4&4&3&3&3&3&4&3&0&0&0&0&0&0&0&0&0&0&0&0&0&0&0&0&0&0&0&0&0&0&0&0\\
1&1&1&1&0&0&1&1&1&1&0&1&0&0&0&0&0&0&0&0&0&0&0&0&0&0&0&0&0&0&0&0&0&0&0&0\\
1&1&0&0&1&1&1&1&1&0&1&1&0&0&0&0&0&0&0&0&0&0&0&0&0&0&0&0&0&0&0&0&0&0&0&0\\
0&0&1&1&0&0&0&0&0&1&0&0&0&0&0&0&0&0&0&0&0&0&0&0&0&0&0&0&0&0&0&0&0&0&0&0\\
0&0&0&0&0&0&0&0&0&0&0&0&0&0&0&0&6&6&8&6&6&6&8&6&0&0&0&0&0&6&0&6&0&8&0&6\\
0&0&0&0&0&0&0&0&0&0&0&0&0&0&0&0&2&2&0&2&2&2&0&2&0&0&0&0&0&2&0&2&0&0&0&2\\
0&0&0&0&0&0&0&0&0&0&0&0&0&0&0&0&2&0&2&2&2&0&2&2&0&0&0&0&0&2&0&0&0&2&0&2\\
0&0&0&0&0&0&0&0&0&0&0&0&0&0&0&0&0&2&0&0&0&2&0&0&0&0&0&0&0&0&0&2&0&0&0&0\\
0&0&0&0&0&0&0&0&0&0&0&0&3&3&4&3&0&0&0&0&0&0&0&0&3&3&4&3&3&0&3&0&4&0&3&0\\
0&0&0&0&0&0&0&0&0&0&0&0&1&1&0&1&0&0&0&0&0&0&0&0&1&1&0&1&1&0&1&0&0&0&1&0\\
0&0&0&0&0&0&0&0&0&0&0&0&1&0&1&1&0&0&0&0&0&0&0&0&1&0&1&1&1&0&0&0&1&0&1&0\\
0&0&0&0&0&0&0&0&0&0&0&0&0&1&0&0&0&0&0&0&0&0&0&0&0&1&0&0&0&0&1&0&0&0&0&0\\
3&3&3&3&4&4&3&3&3&3&4&3&0&0&0&0&0&0&0&0&0&0&0&0&0&0&0&0&0&0&0&0&0&0&0&0\\
3&3&3&3&4&4&3&3&3&3&4&3&0&0&0&0&3&3&4&3&3&3&4&3&0&0&0&0&0&3&0&3&0&4&0&3\\
1&1&1&1&0&0&1&1&1&1&0&1&0&0&0&0&0&0&0&0&0&0&0&0&0&0&0&0&0&0&0&0&0&0&0&0\\
1&1&1&1&0&0&1&1&1&1&0&1&0&0&0&0&1&1&0&1&1&1&0&1&0&0&0&0&0&1&0&1&0&0&0&1\\
1&1&0&0&1&1&1&1&1&0&1&1&0&0&0&0&0&0&0&0&0&0&0&0&0&0&0&0&0&0&0&0&0&0&0&0\\
1&1&0&0&1&1&1&1&1&0&1&1&0&0&0&0&1&0&1&1&1&0&1&1&0&0&0&0&0&1&0&0&0&1&0&1\\
0&0&1&1&0&0&0&0&0&1&0&0&0&0&0&0&0&0&0&0&0&0&0&0&0&0&0&0&0&0&0&0&0&0&0&0\\
0&0&1&1&0&0&0&0&0&1&0&0&0&0&0&0&0&1&0&0&0&1&0&0&0&0&0&0&0&0&0&1&0&0&0&0
\end{array}\right)\)
}
\]
\\

\[
\sigma_{A2}=\left(\begin{array}{*{16}c}
2&0&2&0&0&0&0&0&2&2&0&0&2&2&0&0\\
2&4&2&4&0&0&0&0&2&2&0&0&2&2&0&0\\
1&0&1&0&1&0&1&0&1&1&1&1&1&1&1&1\\
1&2&1&2&1&2&1&2&1&1&1&1&1&1&1&1\\
1&0&1&0&1&0&1&0&1&1&1&1&1&1&1&1\\
1&2&1&2&1&2&1&2&1&1&1&1&1&1&1&1\\
0&0&0&0&2&0&2&0&0&0&2&2&0&0&2&2\\
0&0&0&0&2&4&2&4&0&0&2&2&0&0&2&2\\
2&2&2&2&0&0&0&0&4&4&0&0&2&2&0&0\\
1&1&1&1&1&1&1&1&2&2&2&2&1&1&1&1\\
1&1&1&1&1&1&1&1&2&2&2&2&1&1&1&1\\
0&0&0&0&2&2&2&2&0&0&4&4&0&0&2&2\\
2&2&2&2&0&0&0&0&0&0&0&0&2&2&0&0\\
1&1&1&1&1&1&1&1&0&0&0&0&1&1&1&1\\
1&1&1&1&1&1&1&1&0&0&0&0&1&1&1&1\\
0&0&0&0&2&2&2&2&0&0&0&0&0&0&2&2
\end{array}\right)
\]
\\
\[
\sigma_{A3}=\left(\begin{array}{*{23}c}
1&1&2&3&1&1&2&2&2&1&3&3&0&0&0&0&0&1&1&2&3&3&2\\
0&0&1&1&0&0&2&3&2&0&1&1&0&0&0&0&0&0&0&2&1&1&2\\
3&3&0&0&3&4&2&0&2&3&3&3&3&3&3&3&3&4&4&2&3&3&2\\
1&1&0&0&1&0&2&0&2&1&1&1&1&1&1&1&1&0&0&2&1&1&2\\
2&2&3&3&2&2&1&2&1&2&2&2&0&0&0&0&0&2&2&1&2&2&1\\
1&1&2&1&1&1&2&4&2&1&2&2&0&0&0&0&0&1&1&2&2&2&2\\
1&1&2&3&1&2&1&2&1&1&2&2&0&0&0&0&0&2&2&1&2&2&1\\
1&1&0&0&1&1&1&0&1&1&1&1&1&1&1&1&1&1&1&1&1&1&1\\
2&2&2&1&2&1&0&0&0&2&0&0&0&0&0&0&0&1&1&0&0&0&0\\
0&0&0&0&0&0&2&2&2&0&1&1&0&0&0&0&0&0&0&2&1&1&2\\
1&1&1&1&1&1&1&1&1&1&0&0&0&0&0&0&0&1&1&1&0&0&1\\
1&1&1&0&1&0&1&1&1&1&1&1&0&0&0&0&0&0&0&1&1&1&1\\
2&2&2&2&2&2&3&3&3&2&3&3&0&0&0&0&3&2&2&3&3&3&3\\
1&1&1&1&1&1&1&1&1&1&1&1&0&0&0&0&1&1&1&1&1&1&1\\
1&1&1&1&1&1&1&1&1&1&1&1&0&0&0&0&1&1&1&1&1&1&1\\
1&1&1&1&1&1&0&0&0&1&0&0&0&0&0&0&0&1&1&0&0&0&0\\
0&0&0&0&0&0&0&0&0&0&0&0&20&20&20&20&15&0&0&0&0&0&0\\
2&2&2&2&2&2&1&1&1&2&1&1&0&0&0&0&0&2&2&1&1&1&1\\
1&1&1&1&1&1&0&0&0&1&0&0&0&0&0&0&0&1&1&0&0&0&0\\
1&1&1&0&1&0&1&1&1&1&1&1&0&0&0&0&0&0&0&1&1&1&1\\
0&0&0&0&0&0&1&1&1&0&1&1&0&0&0&0&0&0&0&1&1&1&1\\
2&2&2&2&2&2&0&0&0&2&0&0&0&0&0&0&0&2&2&0&0&0&0\\
0&0&0&2&0&2&0&0&0&0&0&0&0&0&0&0&0&2&2&0&0&0&0
\end{array}\right)
\]
\\

\[
\sigma_{A4}=\left(\begin{array}{*{21}c}
2&1&1&1&2&2&2&2&2&3&1&1&2&2&2&2&2&2&3&2&2\\
2&2&2&2&2&2&2&2&2&0&2&1&2&0&0&0&0&2&0&0&0\\
2&2&2&2&2&2&2&2&2&0&2&1&2&0&0&0&0&2&0&0&0\\
0&1&1&1&0&0&0&0&0&4&1&1&0&4&4&4&4&0&4&4&4\\
0&0&0&0&0&0&0&0&0&0&0&0&0&0&0&1&0&0&0&0&1\\
0&0&0&0&0&0&0&0&1&0&0&0&0&0&0&1&0&1&0&0&1\\
3&1&1&1&1&1&3&1&1&2&1&1&1&1&1&1&1&1&2&1&1\\
0&0&0&0&0&0&0&0&0&0&0&0&0&1&1&0&1&0&0&1&0\\
2&0&0&0&3&2&1&3&1&0&0&0&2&1&1&0&1&1&0&1&0\\
3&1&1&1&1&1&2&1&1&3&1&1&1&1&1&1&1&1&2&1&1\\
0&3&3&3&0&0&0&0&0&0&3&1&0&0&0&0&0&0&0&0&0\\
0&4&4&4&0&0&0&0&0&0&4&8&0&0&0&0&0&0&0&0&0\\
0&0&0&0&1&1&0&1&1&0&0&0&1&0&0&0&0&1&0&0&0\\
0&0&0&0&1&0&0&1&1&0&0&0&0&2&2&1&1&1&0&1&1\\
2&0&0&0&3&3&2&3&2&0&0&0&3&0&0&0&0&2&0&0&0\\
0&0&0&0&0&0&0&0&1&2&0&0&0&2&2&2&1&1&1&1&2\\
0&0&0&0&0&1&0&0&0&1&0&0&1&1&1&1&2&0&1&2&1\\
0&0&0&0&0&1&1&0&1&0&0&0&1&0&0&0&0&1&0&0&0\\
0&1&1&1&0&0&1&0&0&0&1&1&0&0&0&0&0&0&1&0&0\\
0&0&0&0&0&0&0&0&0&1&0&0&0&1&1&1&1&0&1&1&1\\
0&0&0&0&0&0&0&0&0&0&0&0&0&0&0&1&1&0&1&1&1
\end{array}\right)
\]
\\

\[
\resizebox{0.9\hsize}{!}{
\(\sigma^2_{A5}=\left(\begin{array}{*{24}c}
-8&-7&-4&-8&-2&4&0&-2&-11&-4&0&0&0&0&0&0&9&0&0&0&0&0&0&0\\
1&2&14&16&1&-2&0&1&1&14&0&0&0&0&0&0&0&0&0&0&0&0&0&0\\
0&0&9&6&0&2&-2&0&0&6&0&0&0&-2&2&0&0&2&-4&2&-2&-2&2&0\\
0&0&-3&0&0&-2&2&0&0&-3&0&0&0&2&-2&0&0&-2&4&-2&2&2&-2&0\\
9&7&6&8&3&-2&0&2&11&-2&0&0&0&0&0&0&-9&0&0&0&0&0&0&0\\
0&0&3&3&0&1&2&0&0&3&0&0&0&2&-2&0&0&-2&4&-2&2&2&-2&0\\
-3&-3&6&6&0&-2&5&0&-3&6&0&0&0&2&-2&0&3&-2&4&-2&2&2&-2&0\\
-2&-2&8&8&1&0&0&2&-4&16&0&0&0&0&0&0&3&0&0&0&0&0&0&0\\
3&3&0&0&0&0&0&0&6&0&0&0&0&0&0&0&-3&0&0&0&0&0&0&0\\
0&0&3&3&0&0&0&0&0&6&0&0&0&0&0&0&0&0&0&0&0&0&0&0\\
-13&-11&10&8&-1&2&0&-1&-13&10&9&6&0&0&0&6&6&0&2&0&2&2&0&0\\
0&0&0&0&0&0&0&0&0&0&-3&0&0&0&0&-3&3&0&-2&0&-2&-2&0&0\\
0&0&0&0&0&8&-8&0&0&0&7&9&5&-8&8&7&-7&0&-6&0&2&2&0&0\\
0&0&0&0&0&0&0&0&0&0&3&3&0&3&0&3&-3&0&-2&0&-2&-2&0&0\\
-3&-3&3&3&0&2&-2&0&-3&3&1&3&2&-2&5&1&2&0&-2&0&0&0&0&0\\
0&0&0&0&0&0&0&0&0&0&3&3&0&0&0&6&-3&0&0&0&0&0&0&0\\
-13&-11&10&8&-1&2&0&-1&-13&10&0&0&0&0&0&0&15&0&0&0&0&0&0&0\\
0&0&0&0&0&-8&8&0&0&0&2&0&-2&8&-8&2&-2&-1&16&-4&4&4&-2&0\\
0&0&0&0&0&0&0&0&0&0&0&0&0&0&0&0&0&0&5&0&2&2&0&0\\
0&0&0&0&0&0&0&0&0&0&2&0&-2&0&0&2&-2&2&0&5&0&0&2&0\\
0&0&0&0&0&4&-4&0&0&0&-2&0&2&-4&4&-2&2&0&-4&0&1&0&0&0\\
0&0&0&0&0&-4&4&0&0&0&2&0&-2&4&-4&2&-2&0&12&0&4&5&0&0\\
0&0&0&0&0&16&-16&0&0&0&-2&0&2&-16&16&-2&2&6&-24&12&-4&-4&7&0\\
0&0&0&0&0&-4&4&0&0&0&0&0&0&4&-4&0&0&-2&8&-4&2&2&-2&1
\end{array}\right).\)
}
\]

\bibliographystyle{amsalpha}
\bibliography{biblio}
\end{document}

%% file: equivlabeling.tex
\subsection{Equivariant gap labeling}
\label{subsec:equivgap}

In this subsection we mention an equivariant
variant of the gap-labeling conjecture
which is stronger, hence more likely to fail.  In fact, we produce
a counterexample, though we still show that the conjecture
is true in low dimensions.  The reason for mentioning
this equivariant conjecture here is that our methods might be usable
in some other cases to disprove this equivariant
conjecture, even if they do not disprove the original conjecture.
Suppose $\cT$ is a primitive substitution tiling as before, and furthermore
assume that there is a finite group $G$ of symmetries that
acts on the associated tiling space $\Omega$.  (Equivariant $K$-theory
works just as well with $G$ compact Lie, but that case isn't relevant
for substitution tilings.) This group $G$ will of course
act on $\cA_{p}(\Omega)$ and preserve both the unique invariant
measure on $\Omega$ and the unique tracial state
$\tau\co K_0(\cA_{p}(\Omega))\to \bR$.
We have a natural inclusion map on \emph{equivariant} $K$-theory:
\[
K^0_G(\Omega)= K_0^G(C(\Omega))\to K^G_0(\cA_{p}(\Omega)).
\]
The \textbf{equivariant gap-labeling conjecture} asserts that this
map induces an isomorphism on all trace invariants.  But there are far
more of these than in the non-equivariant case.  To explain this,
recall that equivariant $K$-theory is a module over $R(G)$, the
representation ring of $G$, which after tensoring with $\bC$
can be identified with the ring of linear combinations of characters
(of irreducible representations of $G$), or with the algebra $Z(G)$ of
class functions on $G$ (this is the center of the group ring $\bC G$).
If $A$ is a (unital, for  simplicity) $C^*$-algebra equipped with an
action of $G$, then $K_0^G(A)$ is the Grothendieck group of equivalence
classes of $G$-invariant
self-adjoint projections $e\in \mathrm{End}(V)\otimes A$, $(V,\pi)$
a finite-dimensional unitary representation space for $G$
\cite[\S11.3]{MR1656031}.  If $\tau$ is a $G$-invariant trace on $A$,
it defines a map $K_0^G(A)\to Z(G)\cong R(G)\otimes \bC$ via
$[e]\mapsto (g\mapsto \tau^V\bigl((\pi(g)\otimes 1_A) e)\bigr)$,
where $\tau^V$ is the trace
induced by $\tau$ on $\mathrm{End}(V)\otimes A$.  Another way of
thinking about this is to observe that $R(G)\otimes \bC$ is a semisimple
ring (a direct sum of copies of $\bC$, one for each conjugacy class
of $G$) and so $K_0^G(A)\otimes \bC$, as a module
over $R(G)\otimes \bC$, is a direct
sum of $\bC$-vector spaces indexed by the conjugacy classes of
$G$. The $\tau$ then gives a trace on each summand. For example, if $A=\bC$
with $\tau\co A\to\bC$ the identity map,
$K_0^G(A)\cong R(G)$, and for a finite-dimensional
representation $(V,\pi)$ of $G$, the class $[V]\in R(G)$ is
represented by $e=1_V\in V$, and $(\pi(g)\otimes 1_A) e = \pi(g)$, whose trace
under $\tau^V$ is just the character of $V$.  As another example,
if $X=G/H$ is a transitive $G$-space, then $K^0_G(X)=R(H)$ (with
$R(G)$-module structure via the restriction map $R(G)\to R(H)$).
If $\tau$ is counting measure on $G/H$ (of total mass $[G:H]$)
and if $(V,\pi)$ is a finite-dimensional representation of $H$,
consider the induced $G$-vector bundle $E=G\times_H V$, which is
corresponds to a finitely generated projective $G$-$C(X)$-module,
with a class in $K^0_G(X)$.  When $G$ is abelian, the trace $\tau^V$ sends this
to the function in $Z(G)$ given by the character of $\pi$ on $H$
extended to be $0$ off of $H$. When $G$ is not abelian, one
gets instead the character of the induced representation
$\text{Ind}_{H\uparrow G}(\pi)$, which is supported on the union
of the conjugates of $H$.  For example, if $G=S_3$, $H=A_3$ (a cyclic
group of order $3$), and $\pi$ is a nontrivial character of $H$,
then $\text{Ind}_{H\uparrow G}(\pi)$ is an irreducible representation
with character $2$ at the identity, $-1$ on the $3$-cycles,
$0$ on the $2$-cycles.

Now recall some ideas from \cite{MR0236951}
and \cite{MR0248277}.  First of all, by \cite[Proposition 3.7]{MR0248277},
every prime ideal $\fp$ of $R(G)$ has a \textbf{support}, which is
a cyclic subgroup $H$ of $G$ determined up to conjugation in $G$.
This is the smallest subgroup $H<G$ such that $\fp$ is the inverse image of a
prime ideal in $R(H)$ under the restriction map $R(G)\to R(H)$.
Segal's localization theorem
(\cite[Proposition 4.1]{MR0234452} and \cite[Theorem 1.1]{MR0236951})
says that for
a compact $G$-space such as $\Omega$, the localization of $K^0_G(\Omega)$
at $\fp$ has the property that $K^0_G(\Omega)_\fp\to K^0_G(\Omega^{(H)})_\fp$
is an isomorphism.  Here  $\Omega^{(H)}=G\cdot \Omega^H$ denotes the
$G$-saturation of the subset fixed by $H$. (When $G$ is abelian,
this is equal to $\Omega^H$.) In particular, if $\fp$ is
the prime ideal of characters that vanish at a point $g\in G$
(if we choose $g=1$, then this $\fp$ is just the augmentation ideal of
$R(G)$), then the
support of $\fp$ is the cyclic subgroup generated by $g$, and if
we invert characters that don't vanish at $g$,
which is harmless if we are evaluating characters at $g$, then we can
restrict to $G\cdot \Omega^g$.
The trace $\tau$ then gives a linear functional defined on
$K^0_G(G\cdot \Omega^g)_\fp\cong K^0_G(\Omega)_\fp$,
and by varying $g$, we get a different functional for each
conjugacy class of cyclic subgroups $\langle g\rangle$.
There are similar functionals defined on
$K^G_0(\cA_{p}(\Omega))$, after localizing at these prime ideals
(one for each element of the group).  These amount to the same as the
map to $Z(G)$ defined before, followed by evaluation at a group element
(or conjugacy class, in the noncommutative case).
\begin{conjecture}[Equivariant gap-labeling conjecture]
  \label{conj:equivlabeling}
  Suppose $\cT$ is a primitive substitution tiling as before, and further
  assume that there is a finite group $G$ of symmetries that
  acts on the associated tiling space $\Omega$.  This group will of course
  preserve the unique tracial state $\tau\co K_0(\cA_{p}(\Omega))\to \bR$.
  But $G$ also acts on $\cA_{p}(\Omega)$. Then the image of
  $K^0_G(\Omega)$ in $Z(G)$ under $\tau$ coincides with the image of
  $K^G_0(\cA_{p}(\Omega))$.
\end{conjecture} 

Let's examine a special case to see what this actually says.
\begin{example}
  \label{ex:Z2action}
  Let $G=\{1,-1\}$ be cyclic of order $2$.  Then $R(G)=\bZ[\chi]/(\chi^2-1)$,
  where $\chi$ is the sign character.  After inverting $2$, this splits
  as $\bZ[\frac12]\times\bZ[\frac12]$, with the two factors corresponding
  to the trivial representation $1$ and to $\chi$. (For a summary
  of the structure of $R(G)$, including the complete prime ideal structure,
  see \cite[\S2]{MR3044609}.)  If $A$ is a unital
  $G$-$C^*$-algebra with a $G$-invariant trace $\tau$, the map
  $K_0^G(A)\to Z(G)$ can be viewed as the pair consisting of the
  usual trace, $\tau\co K_0(A)\to \bR$, along with the restriction
  of $\tau$ to the fixed-point algebra, giving a map $K_0(A^G)\to \bR$.
  (This is because $R(G)$ has exactly two minimal prime ideals,
  the augmentation ideal, supported at $\{1\}$, and the kernel of
  $\chi$, supported on all of $G$.  The remaining prime ideals
  are all maximal ideals with finite residue field.)
  So the equivariant gap-labeling conjecture amounts to the usual
  gap-labeling conjecture \textbf{plus} the assertion that
  the image of the trace on the $G$-invariant subalgebra $\cA_{p}(\Omega)^G$
  coincides with the image of the trace on $K_0(C(\Omega)^G)=K^0(\Omega/G)$.
\end{example}
\begin{remark}
  \label{rem:historyequiv}
  To the best of our knowledge, the equivariant gap-labeling conjecture
  hasn't be formulated in this way before.  But versions of it appear in
  papers such as \cite{MR1898159,MR3304277,MR2718945,MR2669441}.
  Those sources don't deal with equivariant $K$-theory in full generality
  but deal either with the $K$-theory of the crossed product
  $\cA_{p}(\Omega)\rtimes G$ (which is the same as the equivariant $K$-theory
  by the Green-Julg Theorem --- see \cite[\S11.7]{MR1656031})
  or else with the $K$-theory of the fixed-point algebra
  $\cA_{p}(\Omega)^G$ (which corresponds to the part of the
  equivariant $K$-theory attached to the trivial representation).
\end{remark}
  
\subsection{The equivariant Chern character}
\label{subsec:equivChern}
We can study the equivariant gap-labeling conjecture using the equivariant
Chern character, just as the usual Chern character is used to study the
usual gap-labeling.  The equivariant Chern character (say for the
action of a finite group $G$ on a $G$-CW complex like the AP complex)
is defined in \cite{MR1851256}, and gives a ring homomorphism
(that becomes an isomorphism after tensoring the domain with $\bQ$)
\( \ch_G\co K^*_G(X) \to H^*_G (X; R(\text{---}) \otimes\bQ) \).
Here the right-hand side
is \textbf{Bredon cohomology} \cite{MR0214062} of the functor
$R(\text{---})\otimes \bQ\co \mathrm{Or}(G)^{\text{op}}\to \bQ\text{-Mod}$,
where $\mathrm{Or}(G)$ is
the orbit category with objects the transitive $G$-spaces
and morphisms the $G$-equivariant maps,
defined by sending $G/H$ to $R(H)\otimes \bQ$.  In the special case
where $X$ has a single orbit type, say $X=(G/H)\times Y$ with $G$
acting trivially on $Y$, we have $K^*_G(X)=R(H)\otimes K^*(Y)$,
and the equivariant Chern character is just the usual Chern character for
$Y$, tensored with $R(H)$.

Now to prove the equivariant gap-labeling conjecture for a particular
tiling space $\Omega$, one can try to replicate the strategy based
on the diagram \eqref{eq:forreal}.  The difference would be that we
use equivariant $K$-theory and cohomology, and the map $\tau$ to $\bR$
is now replaced by the map $\tau^g$ induced by $\tau$
on equivariant $K$-theory/cohomology, along with evaluation of (virtual)
characters at an element $g\in G$.  For simplicity let's take $G$ to
be abelian.  Fix $g\in G$ and apply the Segal localization theorem
with $\fp=\fp_g$ the prime ideal of $R(G)$
consisting of virtual representations
whose characters vanish at $g$.  The support of $\fp$ is the cyclic
subgroup $H=\langle g\rangle$, so after localizing at $\fp$
(in other words, inverting characters that do not vanish at $g$), we
can replace $\Omega$ by $\Omega^g$.  So diagram \eqref{eq:forreal}
becomes
\begin{equation}
\label{eq:equivChern}
\begin{tikzcd}
  K_0^G(\mathcal{A}_p(\Omega))_\fp
  \arrow[r, "\mathrm{ch}_d^G\circ \chi^{-1}"] \arrow[d, "\tau^g"]
& \check{H}^d_G(\Omega^g;R(\text{---})_\fp\otimes\bQ) \arrow[d, "C_\mu^g"] \\
\bC \arrow[r,  "\mathrm{Id}"]
& || \bC.
\end{tikzcd}
\end{equation}
This reduces to \eqref{eq:forreal} if $g=1$, but in general, if $g\ne 1$,
we need to use $\bC$ as the target instead of $\bR$ since
$G$ can have characters with non-real values at $g$.  (This already happens
for $G$ cyclic of order $>2$.) The restriction of $\tau^g$ to
$K_0^G(C(\Omega))=K^0_G(\Omega)$ via the inclusion of $C(\Omega)$ into
$\cA_p(\Omega)$ can be computed via a $G$-invariant
transversal $N$ to the $\bR$-orbits in $\Omega$.  Since $N$ is totally
disconnected (i.e., $0$-dimensional), there are no denominators
in the equivariant Chern character for $N$, which becomes an
\textbf{integral} isomorphism.  So just as in \cite[Remark 9.7]{ADRS:bloch},
we have $K_0^G(\cA_p(\Omega))\cong K_0^G(C(N)\rtimes \bZ^d)$ and
we want to relate this to
$K_0^G(C(N))=K^0_G(N)\cong \check H^0_G(N;R(\text{---}))$.  The restriction
of $\tau^g$ is $\bZ^d$-invariant, so it factors through
$H^0_G(N;R(\text{---}))_{\bZ^d}\cong H^d_G(\Omega;R(\text{---}))$.
So the equivariant gap-labeling conjecture reduces to the assertion
that for each $g\in G$, the image of $\tau^g$ on $K_0^G(\cA_p(\Omega))$ in
\eqref{eq:equivChern} coincides (integrally!) with the image of
$H^d_G(\Omega^g;R(\text{---}))$.  Just as in the non-equivariant case,
the conjecture follows from diagram \eqref{eq:equivChern}
in dimensions up to $3$ when there are no denominators in the Chern
character.
\begin{theorem}
  \label{thm:equivlabeling}
  The equivariant gap-labeling conjecture, Conjecture \ref{conj:equivlabeling},
  holds when the dimension $d$ is $\le 3$. However, it \textbf{fails}
  for a periodic lattice tiling with $\Omega=\bT^4$
  \textup{(}the $4$-torus\textup{)} and $G$ the cyclic
  group of order $2$ interchanging the two factors in $\bT^4=\bT^2\times \bT^2$.
\end{theorem}
\begin{proof}
  We have already explained the proof when $d\le 3$.  It remains to show
  why the conjecture fails for $\Omega=\bT^2\times \bT^2$ and $G$ the
  cyclic group of order $2$ interchanging the two factors.   By Example
  \ref{ex:Z2action}, it suffices to show that the image of the
  trace on  $H^4(\Omega/G;\bZ) = H^4(SP^2(\bT^2);\bZ)$,
  $SP^2(\bT^2)$ the second symmetric product of $\bT^2$, is strictly smaller
  than the image of the trace on the fixed-point algebra
  $\cA_{p}(\Omega)^G$.  Via the analogue of \eqref{eq:forreal},
  the issue is to show that the top degree Chern character on
  $K^0(SP^2(\bT^2))$ has image which is strictly bigger than $H^4(SP^2(\bT^2);\bZ)$.
  For this we can apply \cite[(6.30) and (7.1)]{MR0151460}, which
  computes the structure of the integral cohomology ring of $SP^2(\bT^2)$.
  (In the language of \cite{MR0151460}, this is the case where $X$ is
  a curve of genus $g=1$ and $n=2$.) Macdonald shows that the
  cohomology ring is torsion-free and is generated over $\bZ$ by
  generators $\xi_1,\xi_2$ of degree $1$ (which anticommute with each other
  and with themselves)
  and a generator $\eta$ of degree $2$, commuting with $\xi_1$ and $\xi_2$,
  with the one relation
  $\eta^2=\eta\xi_1\xi_2=\text{generator of }H^4(SP^2(\bT^2);\bZ)$.
  If $L$ is the complex line bundle with $c_1(L)=\eta$, then
  $\ch([L])=1+\eta+\frac{\eta^2}{2}$, which is \textbf{not} integral,
  and its projection into $H^4(SP^2(\bT^2);\bQ)$ generates an infinite
  cyclic group which properly contains $H^4(SP^2(\bT^2);\bZ)$ (with index
  $2$).
\end{proof}

An example of an aperiodic tiling space with a $\mathbb{Z}_2$-action which is also a counterexample is described in \S \ref{subsec:equivCounter}. We then have the following corollary as an immediate consequence of this counterexample.

\begin{corollary}
The gap-labeling conjecture \textbf{fails} for rotational tiling spaces when the dimension \(d\) is \(\geq 4\).
\end{corollary}

\noindent While patch frequencies are well-defined for rotational tiling spaces, thus one can compute the associated frequency module, it is unclear if there is a Schr\"odinger operator whose gaps in the spectra are labeled by the \(K_0\)-group of the associated \(C^\ast\)-algebra.

\begin{remark}
  We should mention that there seem to be different ways to formulate
  the equivariant conjecture, and at first sight they don't seem to be
  the same.  Let illustrate with the case of Example \ref{ex:Z2action},
  a $\bZ_2$-action.  If $X$ is, say, a Cantor set with a $G$-action,
  where $G=\bZ_2$, equipped with a $G$-invariant measure $\mu$,
  then $A=C(X)$ is generated by projections $\chi_Y$ (corresponding to
  clopen subsets $Y\subseteq X$), whose classes also generate $K^0(X)$,
  and the trace $\tau$ defined by $\mu$
  sends $\chi_Y$ to $\mu(Y)$.  In Example \ref{ex:Z2action}, we said
  that the equivariant conjecture amounts to looking at $\tau$ on $A$
  along with the induced trace on $A^G=C(X/G)$.  If $X^G=Z$, then
  $X$ is the union of $Z$ and the free $G$-space $X\smallsetminus Z$,
  so the (unrenormalized) trace on $C(X/G)$ has total mass
  $\mu(Z) + \frac12 \mu(X\smallsetminus Z)$.  The formulation 
  using \eqref{eq:equivChern} amounts to looking instead at the
  pair consisting of $\mu$ on $X$ and $\mu$ restricted to $X^G=Z$.
  Knowing these, one can compute $\mu(X\smallsetminus Z)$, and thus
  get the same set of invariants.  However, because of the factor of
  $2$ that appears in the calculation, one has to be careful
  in looking at the \emph{integral} range of the trace.  The reduction
  to fixed-point sets in \eqref{eq:equivChern} appeals to the
  Segal Localization Theorem, which requires inverting elements
  of $R(G)$ not in the augmentation ideal $I$.  This involves
  inverting $2$, so there is no contradiction in relating the
  two versions of the equivariant conjecture, the one using the
  quotient space and the one using the fixed-point set.
\end{remark}